\documentclass[11pt,reqno]{amsart}
\usepackage{amsfonts}
\usepackage{amsmath}
\usepackage{mathrsfs}
\usepackage{amssymb}
\usepackage{cite}
\usepackage{needspace}

\newtheorem{theorem}{Theorem}[section]
\newtheorem{lemma}[theorem]{Lemma}
\newtheorem{corollary}[theorem]{Corollary}
\newtheorem{example}[theorem]{Example}

\newtheorem{proposition}[theorem]{Proposition}
\numberwithin{equation}{section}
\makeatletter
\renewenvironment{proof}[1][\proofname]{\par\pushQED{\qed}\normalfont
\topsep6\p@\@plus6\p@\trivlist
\item[\hskip\labelsep\normalfont #1\@addpunct{.}]\ignorespaces}
{\popQED\endtrivlist\@endpefalse}
\makeatother

\allowdisplaybreaks
\AtBeginDocument{%
  \setlength{\abovedisplayskip}{6.8pt plus 2pt minus 3pt}%
  \setlength{\belowdisplayskip}{6.8pt plus 2pt minus 3pt}%
  \setlength{\abovedisplayshortskip}{2pt plus 2pt}%
  \setlength{\belowdisplayshortskip}{5pt plus 2pt minus 2pt}%
}%

\begin{document}

\title[Generalized Principal Eigenvalues of Periodic Nonlocal Operators]{Generalized Principal Eigenvalues of Spatially Periodic Nonlocal Operators with Applications to Noncooperative KPP Systems}
\author[H.-H. Vo]{Hoang-Hung Vo$^{1,2}$}
\author[M.-Z. Xin]{Ming-Zhen Xin$^{3,*}$}
\thanks{\hspace{-.6cm}
  $^1$Department of Mathematics, International University, Ho Chi Minh City, Vietnam.
  \\
  $^2$Vietnam National University Ho Chi Minh City, Ho Chi Minh City, Vietnam.
  \\
  $^3$School of Mathematical Sciences, Shanghai Jiao Tong University, Shanghai 200240, China.
  \\
  E-mails: vhhung@hcmiu.edu.vn (H.-H. Vo); xinmingzhen@sjtu.edu.cn (M.-Z. Xin, $^*$corresponding author).}
\date{}

\begin{abstract}
In this paper, we study generalized principal eigenvalues of the
spatially periodic cooperative nonlocal operator
$\mathcal{L}\boldsymbol{\varphi}=\mathcal{J}*\boldsymbol{\varphi}
-\boldsymbol{\varphi}+A(x)\boldsymbol{\varphi}$
and their applications to noncooperative KPP systems. In the periodic
setting, we characterize the existence of a classical principal
eigenvalue and explicitly construct and classify Radon measure solutions
to the associated eigenproblem when a positive continuous principal
eigenfunction does not exist. On the whole space, we consider two notions
of generalized principal eigenvalues defined through positive test
functions and characterize them through a family of weighted periodic
spectral values and Dirichlet approximations. For large spatial periods,
we identify the spectral limits through a scalar Hamilton--Jacobi cell
problem and, under additional assumptions, derive leading-order
corrections and eigenfunction concentration results. For small spatial
periods, we characterize the limits by an averaged operator and obtain
refined expansions for specific classes of kernels, including formulas
for the difference between the two whole-space generalized principal
eigenvalues. Finally, we apply the spectral theory to establish extinction
and persistence criteria for the associated noncooperative KPP systems,
characterize spreading speeds along rays by a Freidlin--G\"artner-type
formula, and study their asymptotic behavior as the spatial period tends
to infinity or zero.

\medskip
Keywords: Nonlocal dispersal; generalized principal eigenvalue;
periodic coefficients; KPP systems; asymptotic expansion; spreading speed.

MSC 2020: 35P15, 35K57, 45C05, 92D25.
\end{abstract}

\maketitle
	
	\section{Introduction}

We study the following nonlocal Kolmogorov--Petrovskii--Piskunov (KPP)
system with $m$ components in spatially periodic habitats:
\begin{equation}
\label{eq:general-nonlocal-KPP}
\partial_t\boldsymbol u
=\mathcal J*\boldsymbol u-\boldsymbol u+A(x)\boldsymbol u
-\bigl(Q(x)\boldsymbol u\bigr)\circ\boldsymbol u,
\qquad (t,x)\in(0,\infty)\times\mathbb R^N,
\end{equation}
where $\boldsymbol u=(u_1,\ldots,u_m)^{\mathsf T}$,
$\mathcal J=\operatorname{diag}(J_1,\ldots,J_m)$, $\circ$ denotes the
componentwise product, and $A(x)=(a_{ij}(x))_{m\times m}$ and
$Q(x)=(q_{ij}(x))_{m\times m}$ are coefficient matrices.

System \eqref{eq:general-nonlocal-KPP} describes the spatial dynamics
of a population divided into discrete types, such as phenotypes.
A concrete example of such population models is the reaction--diffusion
model of Elliott and Cornell \cite{ElliottCornell2012}, in which two
phenotypes, one with faster growth and the other with faster dispersal,
are linked by mutation and subject to competition within and between
phenotypes. Here $u_i(t,x)$ is the density of type $i$,
$A(x)$ describes growth and linear coupling between types, and $Q(x)$
describes competition within and between types. For example, in mutation models,
the off-diagonal entries of $A$ give transfer rates and its diagonal
entries include the corresponding losses. Although the linearization is cooperative
under the assumptions below, competition can prevent the full
nonlinear system from preserving order.

To account for movement over long distances, we use integral dispersal
as an alternative to local diffusion; see Kot et~al.\
\cite{KotLewisVanDenDriessche1996} and Hutson et~al.\
\cite{HutsonMartinezMischaikowVickers2003}.
The convolution acts componentwise in the spatial variable:
\[
(\mathcal J*\boldsymbol u)_i(t,x)=(J_i*u_i)(t,x)
:=\int_{\mathbb R^N}J_i(x-y)u_i(t,y)\,\mathrm dy,
\qquad i=1,\ldots,m.
\]
Here $J_i(x-y)$ is the dispersal probability density from $y$ to $x$:
$J_i*u_i$ represents immigration at $x$, whereas $-u_i$ represents
emigration. The kernels allow type-dependent dispersal distances and
directional bias.

We further incorporate environmental heterogeneity through
time-independent, spatially periodic coefficients $A(x)$ and $Q(x)$,
which describe repeated patterns in growth and competition; see
Berestycki et~al.\ \cite{BHN2005}. The periods specify the spatial scale
of these environmental variations.
Fix spatial periods $L_1,\ldots,L_N>0$ and set $\Lambda:=\{1,\ldots,m\}$.
We use the following assumptions on the dispersal kernels and
coefficient matrices:

\noindent
(J) $J_i\in C(\mathbb{R}^N)$, $J_i(0)>0$, $J_i(x)\geq0$, and
$\int_{\mathbb{R}^N}J_i(x)\,\mathrm{d}x=1$ for all
$i\in\Lambda$. Moreover, there exist constants
$C>0$ and $\delta>0$ such that
\begin{equation*}
	J_i(x)\leq C(1+|x|)^{-N-\delta},
	\qquad x\in\mathbb R^N,\quad i\in\Lambda.
\end{equation*}

The polynomial upper decay imposed in (J) does not imply the existence
of nonzero exponential moments. Accordingly, formulas involving
exponential weights below retain the actual moment domain rather
than silently assuming it contains a neighborhood of the origin.
In particular, every uniform expansion in the weight is confined to
a compact subset of the interior of that domain, unless explicitly stated
otherwise.

\noindent
(A) The entries of $A$ are continuous and $L_j$-periodic in $x_j$ for
$j=1,\ldots,N$. For every $x\in\mathbb R^N$, $A(x)$ is cooperative and
irreducible.

\noindent
(Q) The entries of $Q$ are continuous, strictly positive, and
$L_j$-periodic in $x_j$ for $j=1,\ldots,N$.

A square matrix is cooperative if its off-diagonal entries are
nonnegative, and irreducible if no simultaneous permutation of its
rows and columns puts it in block upper triangular form with two
nonempty diagonal blocks.
Assumptions (J) and (A) are imposed throughout the paper; (Q) is
required only for the nonlinear results in
Section~\ref{sec:nonlinear-dynamics}.
In components, system \eqref{eq:general-nonlocal-KPP} reads
\begin{equation*}
\partial_t u_i
=(J_i*u_i)(t,x)-u_i
+\sum_{j=1}^m a_{ij}(x)u_j
-u_i\sum_{j=1}^m q_{ij}(x)u_j.
\end{equation*}
To place this problem in context, we review the relevant literature
on spectral theory and propagation in heterogeneous media.

\subsection{Related work}\label{subsec:related-work}

Principal eigenvalues of linearized operators provide spectral
criteria for persistence and extinction and enter the characterization
of spreading speeds in nonlinear reaction--diffusion equations.
In the classical differential setting, consider the scalar elliptic
operator
\begin{equation*}
\mathcal P\phi
=\sum_{i,j=1}^N b_{ij}(x)\partial_{ij}\phi
+\sum_{i=1}^N\beta_i(x)\partial_i\phi+c(x)\phi
\qquad\text{in }\mathbb R^N,
\end{equation*}
with bounded H\"older-continuous coefficients and a symmetric,
uniformly positive definite matrix $B(x)=(b_{ij}(x))$.
On the whole space, the lack of compactness prevents a direct use
of classical bounded-domain principal-eigenvalue theory.
Generalized principal eigenvalues instead use inequalities for positive
test functions. In this regular setting, two classical notions studied
by Berestycki and Rossi \cite{BerestyckiRossi2015} can be written,
with the opposite sign convention adopted in this paper, as
\begin{equation*}
\begin{aligned}
\lambda^*(\mathcal P)
&=\inf\left\{\lambda\in\mathbb R:
\exists\phi\in C^2(\mathbb R^N),\quad \phi>0,\quad
\mathcal P\phi\leq\lambda\phi\right\},\\[2pt]
\lambda'(\mathcal P)
&=\sup\left\{\lambda\in\mathbb R:
\exists\phi\in C^2(\mathbb R^N)\cap L^\infty(\mathbb R^N),\quad \phi>0,\quad
\mathcal P\phi\geq\lambda\phi\right\}.
\end{aligned}
\end{equation*}
Here positivity and both inequalities hold at every point of $\mathbb R^N$.
The first definition permits unbounded test functions, whereas the
second requires boundedness; neither requires periodicity or a positive
uniform lower bound. These different conditions at infinity can lead
to different spectral values.
For scalar parabolic operators in space-time periodic media, Nadin
\cite{Nadin2009Eigenvalue} studied analogous
quantities using time-periodic test functions. He identified the
bounded-test value with the unweighted periodic principal eigenvalue
and the other value by optimization over exponentially weighted
periodic principal eigenvalues.

Nonlocal dispersal introduces a further difficulty: even under spatially
periodic boundary conditions, the operator may fail to possess a
principal eigenvalue with a positive continuous eigenfunction; see
Shen and Xie \cite{Shen2015}. Generalized spectral notions are therefore
needed for the periodic problem itself, not only for the whole-space
problem. For the nonlocal linearization $\mathcal L$ in \eqref{O1},
the corresponding whole-space values $\lambda^*$ and $\lambda'$ are
defined in \eqref{EQ2.1} and \eqref{eq:whole-space-lambda-prime}.
Their test functions are continuous positive vector functions with
finite, continuous convolutions; the definition of $\lambda'$ also
requires boundedness. Coville \cite{Coville2010,Coville2012Harnack} established
criteria and Harnack estimates for the existence of positive principal
eigenfunctions of integral operators. For scalar eigenmeasure results,
see Coville \cite{C2} and Griette \cite{Griette2019}.
For possibly unbounded domains,
Berestycki et~al.\ \cite{Berestycki2016} studied the relations among
generalized principal eigenvalues, domain approximation, and nonlocal
scaling; their persistence criteria \cite{BCV2016Persistence} connect
this theory with population dynamics in the whole space. In
one-dimensional periodic KPP media, a spectral criterion for positive
stationary solutions was obtained by Coville et~al.\
\cite{CovilleDavilaMartinez2008}, and further one-dimensional spectral
properties were studied by Wang and Wang \cite{WangWang2026}.
Du et~al.\ \cite{DuFangNi2025} studied principal eigenvalue limits on
expanding intervals for asymmetric nonlocal operators and their
applications to KPP propagation.
For systems, principal spectral theory has been developed for coupled
integral operators by Su et~al.\ \cite{S20231} and for time-periodic
cooperative operators by Feng et~al.\ \cite{Feng2024}, with time-delay
systems treated by Liang et~al.\ \cite{Liang2019}. Of particular
relevance here, Wang and Zhang
\cite{WZ-1} approximate generalized principal values of cooperative
irreducible nonlocal systems by monotone upper and lower systems
having positive eigenfunctions; see also their time-periodic study
\cite{WangZhangTime2025}.

Once spectral thresholds have been identified, a further question is
how they depend on dispersal and environmental parameters. For
nonlocal operators, dispersal-rate and dispersal-range limits have
been studied in the whole-space persistence theory of Berestycki
et~al.\ \cite{BCV2016Persistence}, scalar time-periodic problems by
Shen and Vo \cite{Shen2019}, and Neumann-type problems by Vo
\cite{Vo2022Neumann}. For scalar time-periodic nonlocal operators,
Su et~al.\ \cite{Su2020a} established equivalent characterizations of
the generalized principal eigenvalue and studied its dependence on
temporal frequency, dispersal rate, and dispersal range; further
dispersal-range asymptotics were obtained by Su et~al.\ \cite{Su2020b}.
Related asymptotic results for
inhomogeneous cooperative systems and time-periodic systems with
coupled or uncoupled integral dispersal were obtained by Thu and Vo
\cite{Thu2025} and Wu et~al.\ \cite{WuLiSunVo2026}, respectively.
Feng et~al.\ \cite{Feng2025} studied dispersal-rate monotonicity
and joint large-frequency and large-dispersal asymptotics for a weighted
time-periodic nonlocal eigenvalue problem.
Other operator and boundary settings include
nonlocal drift studied by Vo \cite{Vo2026Advection},
cooperative-advective free-boundary systems by Bentout and Vo
\cite{BentoutVo2025,BentoutVo2026}, and spectral fractional diffusion
by Trang and Vo \cite{TrangVo2026Fractional}, which differ from the
whole-space convolution considered here.

For differential operators, Lam and Lou \cite{LamLou2016} characterized
the small-diffusion limit for cooperative elliptic systems through the
pointwise Perron root of the coefficient matrix. Liu and Lou studied
advection monotonicity \cite{LiuLou2020} and diffusion--advection
asymptotics for shear flows \cite{LiuLou2024}; see also their survey
\cite{LiuLou2027}. For scalar time-periodic parabolic operators with
Neumann boundary conditions, Liu and Lou \cite{LiuLou2022} used a
time-periodic Hamilton--Jacobi equation to characterize the joint
small-diffusion and low-frequency limit, and established monotonicity
properties and a classification of spectral level sets. Liu
\cite{Liu2025} extended this analysis to time-periodic cooperative
parabolic systems, where the limiting Hamilton--Jacobi equation remains
scalar but its Hamiltonian is determined by a matrix Perron eigenvalue;
under symmetric coupling, that work also establishes frequency
monotonicity and classifies spectral level sets.
The viscosity-solution framework for Hamilton--Jacobi equations was
introduced by Crandall and Lions \cite{CrandallLions1983}; see Crandall
et~al.\ \cite{CrandallIshiiLions1992} for comparison principles and
stability under limits. In periodic homogenization, Lions et~al.\
\cite{LionsPapanicolaouVaradhan1988} characterized effective
Hamiltonians through the additive constants of periodic cell problems.
More directly related to spatial-period variation, Hamel et~al.\
\cite{HamelNadinRoques2011} used periodic Hamilton--Jacobi
equations to characterize large-period spectral and spreading-speed
limits for scalar reaction--diffusion--advection equations.
For nonlocal convolution-type operators on bounded domains with
Dirichlet conditions, understood as zero values outside the domain,
Piatnitski and Rybalko \cite{PR24} characterized the homogenization
limit of the bottom of the spectrum through an effective
Hamilton--Jacobi additive eigenvalue problem. This boundary setting
differs from the spatially periodic spectral problems considered here.
For propagation in shifting habitats, Tao et~al.\
\cite{TaoLiRuanXu2026} used the viscosity-solution theory of
Hamilton--Jacobi equations to classify spreading speeds for nonlocal
dispersal Fisher--KPP equations and characterize the asymptotic states
of their solutions.

Principal eigenvalues also enter the description of propagation in
periodic media. For scalar periodic KPP equations, consider front-like
initial data that are uniformly positive sufficiently far behind the
front and zero sufficiently far ahead. Under suitable hypotheses, the
corresponding spreading speed coincides with the minimal speed of
pulsating fronts in the same direction; see Weinberger
\cite{Weinberger2002}. Spectral variational
formulas for minimal front speeds with local diffusion were established
by Berestycki et~al.\ \cite{BHN2005}, while Shen and Zhang \cite{SZ}
obtained spectral characterizations of spreading speeds for periodic
nonlocal equations. For compactly supported initial data, however,
the spreading speed along a ray in higher dimensions need not equal
the front-like spreading speed in the same direction. G\"artner and
Freidlin \cite{GartnerFreidlin1979} first characterized asymptotic
spreading in periodic Fisher--KPP equations with local diffusion
through weighted periodic principal eigenvalues. Combined with the
minimal-front-speed formula, the Freidlin--G\"artner formula expresses
the spreading speed along each ray in terms of front-like spreading
speeds in all directions and thereby describes the asymptotic
spreading set; see Berestycki et~al.\ \cite{BHN2005} and the extension
to more general reaction terms by Rossi \cite{Rossi2017}.
For spatially homogeneous
scalar KPP equations with asymmetric kernels satisfying an
exponential-moment condition, Guo et~al.\ \cite{GuoLiLiu2026}
obtained a geometric description allowing the intrinsic spreading set
to exclude the origin, for bounded and unbounded initial supports.
Such results emphasize that directional propagation need not be
described by a single isotropic speed.

Establishing these propagation bounds requires a connection between
spectral quantities and nonlinear invasion. Abstract theories for
order-preserving evolution systems developed by Liang and Zhao
\cite{LiangZhao2010}, Ding and Liang \cite{DL2015}, and Du et~al.\
\cite{DuLiShen2022} provide relevant comparison methods. In particular,
Ding and Liang treated periodic noncompact systems under
principal-eigenvalue assumptions. Directional spreading for periodic
nonlocal scalar equations and cooperative systems was studied by
Bao et~al.\ \cite{BaoShenShen2019} and Xin et~al.\ \cite{XinLiBao2023}. For a spatially periodic
age-structured model with local diffusion, Kang and Liu
\cite{KangLiu2025} used generalized principal eigenvalues and
homogenization to characterize persistence and spreading.
For noncooperative systems, cooperative comparison is used in the
local space-time periodic theory of Girardin \cite{Girardin2025},
while periodic nonlocal propagation was studied by Li et~al.\
\cite{LiXinZhao2025}. Slowly decaying kernels can lead to infinite
spreading speeds, as shown by Garnier \cite{Garnier2011} and
Li et~al.\ \cite{LiXinZhao2025}, so finite-speed conclusions require
appropriate moment assumptions. Related spreading--vanishing results
in free-boundary settings were obtained by Li et~al.\
\cite{LiLiWang2022} and Li and Wang \cite{LiWang2025}. For wave
problems, see Coville et~al.\ \cite{CovilleDavilaMartinez2013},
Yang et~al.\ \cite{YangLiWang2020}, Du et~al.\ \cite{DuLiXin2026},
and Vo and Trang \cite{VoTrang2025}.

These developments motivate a basic question: how does the spatial
scale of a heterogeneous habitat affect persistence and invasion in
noncooperative KPP systems with nonlocal dispersal? To answer it, one
must relate whole-space spectral thresholds to periodic spectral
problems and then connect their period dependence to the nonlinear
Cauchy problem.

\subsection{Overview of results}\label{subsec:main-results}

The results are organized into four parts, following the order of the
paper. The spectral characterization, qualitative period limits, and
nonlinear threshold and spreading results hold in arbitrary spatial
dimension; the localization results and explicit Bessel/Laplace
expansions, including their speed corrections, are one-dimensional.
\begin{enumerate}
\item Spectral characterization and eigenmeasures.
We first consider the periodic problem, for which a positive continuous
principal eigenfunction need not exist. Starting from the existence
criteria recalled in Lemma~\ref{L11},
Theorem~\ref{thm:appendix-measure-classification} establishes the
existence of eigenmeasures at the spectral threshold and, under
\eqref{eq:appendix-integrable-case}, classifies classical eigenfunctions,
absolutely continuous eigenmeasures, and eigenmeasures with a singular
part, with an explicit representation in the singular case.
We then turn to the whole space, where positive test functions need
not be bounded away from zero. Theorem~\ref{TH2.1} relates the two
generalized principal eigenvalues to Dirichlet approximations and
weighted periodic spectral values:
\[
\lambda^*=\lim_{l\to\infty}\lambda(l)
=\min_{\mu\in\mathbb R^N}\lambda_L^*(\mu),
\qquad \lambda'=\lambda_L^*(0),
\]
where $\lambda(l)$ is defined in \eqref{EQ2.3}, and
$\lambda_L^*(\mu)$ is the periodic value \eqref{def:lambdaLstar} for
the weighted operator \eqref{Q2}. We also establish convexity and
attainment of the minimum, with uniqueness under the additional
condition in Theorem~\ref{TH2.1}(i). Proposition~\ref{prop:constant-coefficient-stress-test}
provides an independent exact benchmark, while
Example~\ref{ex:one-sided-moment-domain} shows why derivative
claims at moment-domain boundaries require separate hypotheses.
We identify a symmetry condition
ensuring $\lambda^*=\lambda'$ in
Proposition~\ref{prop:detailed-balance-rigidity}.

\item Large-period spectral asymptotics.
With fixed kernels and $A_L(x)=A(x/L)$ for a fixed unit-periodic
profile $A$, Theorem~\ref{thm:large-period-limits} characterizes the
spectral limits as $L\to\infty$ through the effective Hamiltonian
\eqref{eq:effective-periodic-Hamiltonian}. This scalar Hamilton--Jacobi
description retains the system coupling through a matrix Perron
eigenvalue. Under additional assumptions,
Theorem~\ref{thm:universal-kernel-large-period} gives the leading
spectral correction and eigenfunction concentration for scalar
symmetric dispersal, while Theorem~\ref{thm:reversible-system-reduction}
transfers these conclusions to cooperative systems with a common
symmetric kernel and constant diagonal detailed balance.

\item Small-period spectral asymptotics.
For the same scaling, Theorem~\ref{thm:small-period-limits} identifies
the spectral limits as $L\to0$ through the averaged spectral problem
\eqref{eq:averaged-generalized-eigenvalue}, which retains the spatial
dependence of $A$ rather than replacing it by its mean.
Neither qualitative period limit requires exponential moments in every
direction. Under the respective spectral assumptions, we obtain
superalgebraic convergence for smooth kernels
(Theorem~\ref{thm:smooth-superalgebraic}) and explicit expansions for
Bessel and Laplace kernels, including formulas for the difference
between the two whole-space generalized principal eigenvalues
(Subsections~\ref{subsec:Bessel-period} and~\ref{subsec:Laplace-period}).

\item Noncooperative dynamics.
Using cooperative comparison systems, we apply the spectral theory
to the noncooperative KPP problem. Under (Q), for bounded continuous
nonnegative initial data,
Theorem~\ref{thm:nonlocal-threshold-dynamics} gives uniform extinction
if $\lambda'<0$, uniform persistence for strictly positive periodic
data if $\lambda'>0$, and local uniform persistence for all nonzero
data if $\lambda^*>0$. Exponential initial tails are treated in
Theorem~\ref{thm:nonlocal-tail-criteria}; zero thresholds are excluded.
For nonzero compactly supported data with $\lambda^*>0$,
Theorem~\ref{thm:nonlocal-spreading} characterizes spreading speeds
along rays by the Freidlin--G\"artner-type formula
\eqref{eq:planar-and-FG-speeds}, including a criterion for infinite speed.
Under positivity of the limiting spectral minimum,
Theorem~\ref{cor:speed-period-limits} gives period-size speed limits.
Under their respective kernel, spectral, and minimizer hypotheses,
Theorem~\ref{thm:Bessel-directional-speed-expansion} and
Corollary~\ref{cor:Laplace-minimum-speed} give small-period speed corrections.
\end{enumerate}

The paper proceeds from spectral characterization and eigenmeasures in
Section~\ref{sec:whole-space-eigenvalues} to large- and small-period
asymptotics in Sections~\ref{sec:large-period-asymptotics}
and~\ref{sec:quantitative-period}, followed by nonlinear dynamics
and propagation in Section~\ref{sec:nonlinear-dynamics}.

\section{Generalized principal eigenvalues}
\label{sec:whole-space-eigenvalues}

In arbitrary spatial dimension, we establish the spectral relations
used in the nonlinear analysis, presenting the results before their proofs.

\subsection{Main results}

We begin with the periodic spectral problem and its eigenmeasures,
then turn to the corresponding whole-space quantities.

\subsubsection{Notation and preliminaries}\label{subsec:basic-assumptions}

For the spatial periods $L_1,\ldots,L_N$ fixed in the introduction,
define the period lattice, period cell, and torus by
\[
\mathbb L:=\prod_{j=1}^N L_j\mathbb Z,
\qquad
\mathcal C:=\prod_{j=1}^N(0,L_j),
\qquad
\mathbb T^N:=\mathbb R^N/\mathbb L.
\]
We write $\mathbb S^{N-1}=\{e\in\mathbb R^N:|e|=1\}$ for the unit sphere.

	For the periodic spectral problem, let $X$ be the space of continuous
periodic vector-valued functions:
\[
X=\{\boldsymbol{\varphi}\in C(\mathbb R^N,\mathbb R^m):
\boldsymbol{\varphi}(x+p)=\boldsymbol{\varphi}(x)
\ \text{for all }x\in\mathbb R^N,\ p\in\mathbb L\},
\]
	equipped with the norm
\[
\|\boldsymbol\varphi\|_X=\|\boldsymbol\varphi\|_\infty
:=\sup_{x\in\mathbb R^N}\max_{i\in\Lambda}|\varphi_i(x)|.
\]
Matrix and operator norms are the induced norms unless stated otherwise.
We define the positive cones, with inequalities understood componentwise
for every $x\in\mathbb R^N$, by
\[
X^+=\{\boldsymbol\varphi\in X:\boldsymbol\varphi(x)\geq\boldsymbol0\},\qquad
X^{++}=\{\boldsymbol\varphi\in X:\boldsymbol\varphi(x)>\boldsymbol0\}.
\]
	For $\boldsymbol{\varphi}, \boldsymbol{\psi} \in X$, we write $\boldsymbol{\varphi} \geq \boldsymbol{\psi}$ if $\boldsymbol{\varphi}(x) \geq \boldsymbol{\psi}(x)$ for all $x \in \mathbb{R}^N$; $\boldsymbol{\varphi} > \boldsymbol{\psi}$ if $\boldsymbol{\varphi} \geq \boldsymbol{\psi}$ but $\boldsymbol{\varphi} \not\equiv \boldsymbol{\psi}$; and $\boldsymbol{\varphi} \gg \boldsymbol{\psi}$ if $\boldsymbol{\varphi}(x) > \boldsymbol{\psi}(x)$ for all $x \in \mathbb{R}^N$.

The linearization of \eqref{eq:general-nonlocal-KPP} at $\boldsymbol0$ is
\begin{equation}\label{O1}
\mathcal L[\boldsymbol\varphi](x)
=\mathcal J*\boldsymbol\varphi(x)-\boldsymbol\varphi(x)+A(x)\boldsymbol\varphi(x),
\end{equation}
where $\boldsymbol\varphi=(\varphi_1,\ldots,\varphi_m)^{\mathsf T}$ and
\[
(\mathcal J*\boldsymbol\varphi)_i(x)
=\int_{\mathbb R^N}J_i(x-y)\varphi_i(y)\,\mathrm dy,
\qquad i\in\Lambda.
\]

\Needspace{7\baselineskip}
	The operator $\mathcal L$ may not admit a positive continuous principal
eigenfunction \cite{Coville2010,KLS,DL2015,Shen2015,SZ}. Following
\cite{Berestycki2016,S20231,S20232,Thu2025,Su2020a}, we define the periodic
generalized principal eigenvalues by
	\begin{equation}\label{def:lambdaLstar}
		\lambda_L^* = \inf\{\lambda \in \mathbb{R} \mid \exists \boldsymbol{\varphi} \in X^{++},\ \mathcal{L}[\boldsymbol{\varphi}] - \lambda \boldsymbol{\varphi} \leq 0\}
	\end{equation}
	and
	\begin{equation}\label{def:lambdaLprime}
		\lambda_L' = \sup\{\lambda \in \mathbb{R} \mid \exists \boldsymbol{\varphi} \in X^{++},\ \mathcal{L}[\boldsymbol{\varphi}] - \lambda \boldsymbol{\varphi} \geq 0\}.
	\end{equation}
	These infimum and supremum characterizations are called the upper and
lower Collatz--Wielandt values, respectively. Here the subscript $L$
records the fixed periods $(L_1,\ldots,L_N)$; a common scalar period
will be varied only from Section~\ref{sec:periodic-eigenvalue-limits} onward.
	If there exists a function $\boldsymbol{\varphi} \in X^{++}$ and a constant $\lambda_L \in \mathbb{R}$ such that $(\mathcal{L} - \lambda_L)[\boldsymbol{\varphi}] = 0$, then we call $\lambda_L$ and $\boldsymbol{\varphi}$ the principal eigenvalue and principal eigenfunction of $\mathcal{L}$ in $X$, respectively.
	
For a cooperative irreducible matrix $B$, write
\[
s_{\rm PF}(B):=\max\{\operatorname{Re}\zeta:\zeta\in\sigma(B)\},
\]
where $\sigma(B)$ denotes its spectrum. Set
$s(x):=s_{\rm PF}(A(x)-I)$, with $I$ the $m\times m$ identity matrix.
Consider the eigenvalue problem
\begin{equation}\label{s(x)}
	(A(x)-I) \boldsymbol{\phi}(x) = s(x) \boldsymbol{\phi}(x).
\end{equation}
Since $A(x)-I$ is cooperative and irreducible, its Perron spectral bound
$s(x)$ is a simple eigenvalue with a positive eigenvector.  The map $s$ is
continuous and periodic, and hence defines a function on
$\mathbb T^N$. Let
\begin{equation*}
	s^* := \max_{x \in \mathbb{T}^N} s(x),\qquad
\Omega_0 := \bigl\{ x^* \in \mathbb{T}^N : s(x^*) = s^* \bigr\}.
\end{equation*}

We write $\boldsymbol1$ for the vector whose components are all one.
For later use, we record the matrix Collatz--Wielandt bounds:
for every cooperative irreducible matrix $B$,
\begin{equation}\label{eq:matrix-Collatz-Wielandt}
\min_i\frac{(Bv)_i}{v_i}\leq s_{\rm PF}(B)
\leq\max_i\frac{(Bv)_i}{v_i}\qquad(v\in(0,\infty)^m).
\end{equation}
Every positive upper test function satisfies
$(A(x)-I)\boldsymbol\varphi(x)\leq\lambda\boldsymbol\varphi(x)$,
because its convolution is nonnegative. The matrix bounds
\eqref{eq:matrix-Collatz-Wielandt} therefore give
\begin{equation}\label{eq:periodic-multiplication-edge}
\lambda_L^*\geq s^*.
\end{equation}
The same bound holds for each admissible weighted operator in \eqref{Q2}.

We recall the following periodic versions of the spectral criteria
of Bao and Shen \cite{BS} and the principal eigenvalue criteria and
Collatz--Wielandt characterization of Wang and Zhang
\cite[Lemmas~2.1, A.1 and Theorem~A]{WZ-1}.

\begin{lemma}\label{L11}
	Assume that (J) and (A) hold. Then the following statements hold:
	\begin{enumerate}
		\item[(i)] $\lambda_L^*$ is a principal eigenvalue of $\mathcal{L}$ on $X$ if and only if $\lambda^*_L > s^*$.
		\item[(ii)] If $1/(s^* - s(x)) \notin L^1(\mathbb{T}^N)$, then $\lambda_L^*$ is a principal eigenvalue of $\mathcal{L}$ on $X$. Here $L^1(\mathbb T^N)$ denotes the scalar space of integrable functions.
		\item[(iii)] $\lambda_L^* = \lambda_L'$.
	\end{enumerate}
\end{lemma}

The arguments in the cited results apply on the periodic space $X$.
They also apply to the weighted operators in \eqref{Q2} for
$\mu\in\mathcal I$: the weighted kernels are integrable and positive
near zero, and their periodic convolutions are compact by translation
continuity in $L^1$. Thus all three conclusions remain valid with
$\lambda_L^*,\lambda_L'$ replaced by
$\lambda_L^*(\mu),\lambda_L'(\mu)$, with the same $s^*$.
We omit the repeated proofs.

\subsubsection{Periodic eigenmeasures}
\label{subsec:periodic-measures}

We now establish a measure representation when the periodic spectral
value coincides with $s^*$. For a finite nonnegative vector-valued Borel
measure $\boldsymbol\nu$ on $\mathbb T^N$, denote by
$\widetilde{\boldsymbol\nu}$ its locally finite periodic extension to
$\mathbb R^N$ and define
\[
\bigl(\mathcal J[\boldsymbol\nu]\bigr)_i(x)
:=\int_{\mathbb R^N}J_i(x-y)\,\mathrm d\widetilde\nu_i(y).
\]
The measure eigenvalue equation is
\begin{equation}\label{eq:1.3}
\mathcal J[\boldsymbol\nu](x)\,\mathrm dx
+\bigl(A(x)-(\lambda_L^*+1)I\bigr)\,\mathrm d\boldsymbol\nu(x)
=\boldsymbol0,
\end{equation}
understood in the sense of vector-valued distributions on $\mathbb T^N$.
We use the Lebesgue decomposition
$\mathrm d\boldsymbol\nu=\boldsymbol\nu_{\rm ac}\,\mathrm dx
+\mathrm d\boldsymbol\nu_{\rm s}$, where $\boldsymbol\nu_{\rm ac}$ is the
density of the absolutely continuous part and $\boldsymbol\nu_{\rm s}$
is singular with respect to Lebesgue measure.

By Lemma~\ref{L11}(i)--(ii), the equality $\lambda_L^*=s^*$ can occur only if
\begin{equation}\label{eq:appendix-integrable-case}
\frac1{s^*-s(\cdot)}\in L^1(\mathbb T^N).
\end{equation}
Under (J), (A), and \eqref{eq:appendix-integrable-case}, we introduce the
operator that determines the absolutely continuous and singular alternatives.
For every $x\in\mathbb T^N$, set
\begin{equation*}
C(x):=(s^*+1)I-A(x).
\end{equation*}
For $x\notin\Omega_0$, write $C(x)^{-1}=(c_{ij}(x))_{m\times m}$.
Let $\boldsymbol v(x)$ be
the Perron vector of $A(x)-I$, normalized by $\sum_i v_i(x)=1$.
The entries $c_{ij}$ are integrable under
\eqref{eq:appendix-integrable-case}, as justified in
Lemma~\ref{lem:appendix-K-principal}; their values on the null set
$\Omega_0$ may be chosen arbitrarily. Define the positive operator
\begin{equation}\label{eq:appendix-K}
(\mathcal R\boldsymbol g)_i(x)
:=\int_{\mathbb R^N}J_i(x-y)
\sum_{j=1}^m c_{ij}(y)g_j(y)\,\mathrm dy,
\qquad \boldsymbol g\in X,\quad i\in\Lambda.
\end{equation}

We write $\gamma_1=r(\mathcal R)$ for the spectral radius of $\mathcal R$.

For $x\in\Omega_0$, the Perron vector satisfies
\begin{equation}\label{eq:appendix-Perron-null-vector}
C(x)\boldsymbol v(x)=\boldsymbol0,
\qquad \sum_{i=1}^m v_i(x)=1.
\end{equation}
For a finite positive measure $\theta$ supported on $\Omega_0$, write
\begin{equation*}
\mathrm d\boldsymbol\eta_\theta(x)
:=\boldsymbol v(x)\,\mathrm d\theta(x).
\end{equation*}

\begin{theorem}
\label{thm:appendix-measure-classification}
Assume that (J) and (A) hold. If $\lambda_L^*=s^*$, then
\eqref{eq:1.3} has a nonzero finite nonnegative vector-valued Borel measure
solution on $\mathbb T^N$. The singular part of every such solution satisfies
\begin{equation}\label{eq:singular-support-Omega0}
\operatorname{supp}\nu_{{\rm s},i}\subset\Omega_0,
\qquad i=1,\ldots,m.
\end{equation}
Here $\operatorname{supp}$ denotes the support of a measure.

Under \eqref{eq:appendix-integrable-case}, the following alternatives hold:
\begin{enumerate}
\item[(i)] $\gamma_1>1$ if and only if $\lambda_L^*>s^*$. In this case
$\lambda_L^*$ is a classical principal eigenvalue.

\item[(ii)] $\gamma_1=1$ if and only if $\lambda_L^*=s^*$ and
there exists a nonnegative nontrivial function
$\boldsymbol\varphi\in L^1(\mathbb T^N,\mathbb R^m)$ such that
$\mathrm d\boldsymbol\nu=\boldsymbol\varphi\,\mathrm dx$ solves
\eqref{eq:1.3}. In this case
every nonzero nonnegative solution is absolutely continuous, and its density
is unique up to multiplication by a positive constant.

\item[(iii)] $\gamma_1<1$ if and only if $\lambda_L^*=s^*$ and
\eqref{eq:1.3} has a nonnegative measure solution with a nontrivial singular
part. In this case every nonzero nonnegative solution has a singular part
supported on $\Omega_0$. More
precisely, for every nonzero finite positive measure $\theta$ on $\Omega_0$,
set
\begin{equation}\label{eq:singular-eigenmeasure-construction}
\boldsymbol F_\theta:=\mathcal J[\boldsymbol\eta_\theta],\qquad
\boldsymbol g_\theta:=(I-\mathcal R)^{-1}\boldsymbol F_\theta,\qquad
\boldsymbol f_\theta:=C^{-1}\boldsymbol g_\theta.
\end{equation}
Then
\begin{equation}\label{eq:singular-eigenmeasure-representation}
\mathrm d\boldsymbol\nu_\theta
=\boldsymbol f_\theta(x)\,\mathrm dx
+\mathrm d\boldsymbol\eta_\theta(x)
\end{equation}
solves \eqref{eq:1.3}; conversely, every nonzero nonnegative solution with
a singular part has this form.
\end{enumerate}
\end{theorem}

\subsubsection{Whole-space eigenvalues}
\label{subsec:whole-space-main-results}

We relate whole-space spectral values to weighted periodic problems.
For unbounded test functions, the convolution must be finite. We therefore use the admissible class
\[
\mathscr A_{\mathcal J}^{+}:=\bigl\{\boldsymbol\phi\in C(\mathbb R^N,(0,\infty)^m):
J_i*\phi_i\in C(\mathbb R^N)\text{ is finite for each }i\in\Lambda\bigr\}.
\]
In particular, every positive bounded continuous vector function belongs to this class.
\begin{equation}\label{EQ2.1}
		\lambda^* = \inf\big\{\lambda \in \mathbb{R} \mid \exists \boldsymbol{\phi} \in \mathscr A_{\mathcal J}^{+},\ \mathcal{L}[\boldsymbol{\phi}] - \lambda \boldsymbol{\phi} \leq \boldsymbol0\big\},
	\end{equation}
	\begin{equation}\label{eq:whole-space-lambda-prime}
	\lambda' = \sup\big\{\lambda \in \mathbb{R} \mid \exists \boldsymbol{\phi} \in L^{\infty}(\mathbb{R}^N, \mathbb{R}^m) \cap C(\mathbb{R}^N, \mathbb{R}^m),\ \boldsymbol{\phi}(x) > \boldsymbol0,\ \mathcal{L}[\boldsymbol{\phi}] - \lambda \boldsymbol{\phi} \geq \boldsymbol0\big\}.
	\end{equation}
	For the bounded-cube approximation, define
	\begin{equation}\label{EQ2.3}
		\lambda(l) = \inf\big\{\lambda \in \mathbb{R} \mid \exists \boldsymbol{\phi} \in C([-l, l]^N, \mathbb{R}^m),\ \boldsymbol{\phi}(x) > \boldsymbol0,\ \mathcal{L}^l[\boldsymbol{\phi}] - \lambda \boldsymbol{\phi} \leq \boldsymbol0\big\},
	\end{equation}
	where 
	\[
	\mathcal{L}^l[\boldsymbol{\phi}](x) = \mathcal{J}^l[\boldsymbol{\phi}](x) - \boldsymbol{\phi}(x) + A(x)\boldsymbol{\phi}(x),
	\]
	with the integral operator $\mathcal{J}^l$ defined componentwise as
	\[
	(\mathcal{J}^l[\boldsymbol{\phi}])_i(x)
	=\int_{[-l,l]^N}J_i(x-y)\phi_i(y)\,\mathrm{d}y,
	\quad i\in\Lambda.
	\]
	
	To relate these whole-space quantities to periodic spectral problems,
we introduce the exponential moments and their common domain:
\begin{equation}\label{eq:kernel-exponential-moments}
M_i(\mu):=\int_{\mathbb R^N}J_i(y)e^{\mu\cdot y}\,\mathrm dy,
\qquad
\mathcal I:=\{\mu\in\mathbb R^N:M_i(\mu)<\infty\text{ for all }i\in\Lambda\}.
\end{equation}
	For $\mu\in\mathcal I$, let $\lambda_L^*(\mu)$ and $\lambda_L'(\mu)$
be the values in \eqref{def:lambdaLstar}--\eqref{def:lambdaLprime}
for the weighted operator
	\begin{equation}\label{Q2}
		\mathcal{L}_{\mu}[\boldsymbol{\phi}](x) = \mathcal{J}_{\mu}[\boldsymbol{\phi}](x) - \boldsymbol{\phi}(x) + A(x)\boldsymbol{\phi}(x), \quad x\in\mathbb{R}^N,\ \boldsymbol{\phi} \in X,
	\end{equation}
	where
	\[
	(\mathcal{J}_{\mu}[\boldsymbol{\phi}])_i(x) = \int_{\mathbb{R}^N} J_i(x-y) e^{\mu \cdot (x-y)} \phi_i(y) \, \mathrm{d}y, \quad i\in\Lambda.
	\]
	In particular, $\lambda_L^*(0)=\lambda_L^*$ and
	$\lambda_L'(0)=\lambda_L'$.
	For $\mu\notin\mathcal I$, we set
	$\lambda_L^*(\mu)=\lambda_L'(\mu)=+\infty$.  Thus every minimization
	over $\mathbb R^N$ below is understood in the extended-valued sense.

\begin{lemma}\label{lem:interior-moment-regularity}
Assume (J). The common moment domain $\mathcal I$ in
\eqref{eq:kernel-exponential-moments} is convex and contains the origin.
For every compact set $K\Subset\operatorname{int}\mathcal I$ and every
nonnegative integer $r$, one has
\begin{equation}\label{eq:interior-moment-derivative-control}
\max_{1\leq i\leq m}\sup_{\mu\in K}
\int_{\mathbb R^N}(1+|z|)^rJ_i(z)e^{\mu\cdot z}\,\mathrm dz
<\infty.
\end{equation}
Consequently, the maps $\mu\mapsto J_i(\cdot)e^{\mu\cdot(\cdot)}$
are real-analytic as $L^1(\mathbb R^N)$-valued functions on
$\operatorname{int}\mathcal I$; their derivatives are
$z^\alpha J_i(z)e^{\mu\cdot z}$ for multi-indices $\alpha$.
The corresponding periodic convolution operators are real-analytic in
operator norm on compact subsets of $\operatorname{int}\mathcal I$.
None of these uniform derivative assertions is made at
$\partial\mathcal I$.
\end{lemma}

\begin{proof}
Convexity follows from H\"older's inequality and $0\in\mathcal I$
follows from normalization of the kernels. Fix
$K\Subset\operatorname{int}\mathcal I$.
For every $\mu\in K$, choose $\delta_\mu>0$ such that the cube
$\mu+[-2\delta_\mu,2\delta_\mu]^N$ lies in
$\operatorname{int}\mathcal I$. Its $2^N$ vertices belong to
$\mathcal I$. When $\nu\in\mu+[-\delta_\mu,\delta_\mu]^N$,
we can bound
\[
(1+|z|)^re^{\nu\cdot z}
\leq C_{r,\delta_\mu}\sum_{\omega\in\{-1,1\}^N}
 e^{(\mu+2\delta_\mu\omega)\cdot z},
\text{ for }z\in\mathbb R^N.
\]
Indeed, $(1+|z|)^r\leq C_{r,\delta_\mu}
 e^{\delta_\mu\sum_j|z_j|}$ and the sum contains the
largest exponential in each orthant. The right-hand side is
integrable against each $J_i$ because all the vertices lie in
$\mathcal I$. A finite number of the smaller cubes cover $K$,
which proves \eqref{eq:interior-moment-derivative-control}.
The local exponential majorants permit termwise differentiation and
expansion of the exponential into its absolutely convergent Taylor
series in $L^1$. Thus the weighted kernels are $L^1$-analytic with
the stated derivatives. Finally, convolution on
$C(\mathbb T^N)$ has operator norm at most the $L^1$ norm of
its kernel, so the assertion for the weighted periodic operators
follows.
\end{proof}

\begin{lemma}\label{lem:uniform-weight-coercivity}
Assume that (J) and (A) hold. There exist constants
$b_0,b_1>0$ and $b_2\in\mathbb R$, depending only on the common
lower bound of the kernels near the origin and on the diagonal entries
of $A$, such that, for every choice of the spatial periods and every
$\mu\in\mathcal I$,
\begin{equation}\label{eq:uniform-exponential-weight-coercivity}
\lambda_L^*(\mu)\geq b_0 e^{b_1|\mu|}-b_2.
\end{equation}
The same constants can be chosen for the kernels
$J_{i,n}=\eta_nJ_i$, provided that $\eta_n\geq\tfrac12$ on a fixed
neighborhood of the origin. \end{lemma}

\begin{proof}
Choose $r>0$ and $j_0>0$ such that
$J_i(z)\geq j_0$ for all $i$ and $|z|\leq r$.
For $\mu\ne0$, the ball
$B_{r/4}((r/2)\mu/|\mu|)$ lies in $B_r(0)$, and
$\mu\cdot z\geq r|\mu|/4$ on that ball. Hence
\[
M_i(\mu)\geq j_0|B_{r/4}(0)|e^{r|\mu|/4}
\qquad(i=1,\ldots,m).
\]
For $\mu=0$ the same bound holds using $B_{r/4}(0)$.
If $\boldsymbol\varphi\in X^{++}$ satisfies
$\mathcal L_\mu\boldsymbol\varphi\leq\lambda\boldsymbol\varphi$,
choose a minimum point $x_i$ of $\varphi_i$ on the torus.
The off-diagonal terms are nonnegative, and
\[
\lambda\geq
\frac{(\mathcal L_\mu\boldsymbol\varphi)_i(x_i)}
     {\varphi_i(x_i)}
\geq M_i(\mu)-1+\min_xa_{ii}(x).
\]
Taking the infimum over admissible upper tests proves the first
assertion. One can choose
\[
 b_0=\frac{j_0}{2}|B_{r/4}(0)|,\qquad
 b_1=\frac r4,\qquad b_2=1-\min_{i,x}a_{ii}(x).
\]
For truncated kernels the same reasoning applies with $j_0/2$
in place of $j_0$, so the stated constants work simultaneously.
\end{proof}

\begin{theorem}\label{TH2.1}
		Assume that (J) and (A) hold. Then the following statements hold:
		\begin{enumerate}
			\item[(i)] The extended-valued map
			$\mu\mapsto\lambda_L^*(\mu)$ is convex on $\mathbb R^N$ and
			attains its minimum at some $\mu_0\in\mathcal I$.  If
			$1/(s^*-s(x))\notin L^1(\mathbb T^N)$, then
			it is strictly convex on $\mathcal I$ and its minimizer is unique.
			\item[(ii)] $\lambda^* = \lim_{l \to \infty} \lambda(l) = \min_{\mu \in \mathbb{R}^N} \lambda_L^*(\mu)$.
			\item[(iii)] $\lambda' = \lambda_L^*(0)$.
		\end{enumerate}
	\end{theorem}

\subsubsection{Detailed balance}

The following detailed-balance condition makes the two whole-space
spectral values coincide.

\begin{proposition}\label{prop:detailed-balance-rigidity}
Assume that (J) and (A) hold. Suppose that all dispersal kernels are even,
and there exists a fixed positive diagonal matrix $D$ such that
\[
 D A(x)=A(x)^{\mathsf T}D\qquad\text{for every }x\in\mathbb R^N.
\]
Then the weighted periodic generalized principal eigenvalues satisfy
\[
 \lambda_L^*(\mu)=\lambda_L^*(-\mu)
 \qquad(\mu\in\mathcal I).
\]
In particular, at every fixed spatial period (and without assuming
the existence of a classical principal eigenfunction),
\begin{equation}\label{eq:detailed-balance-exact-threshold}
 \lambda^*=\lambda'=\lambda_L^*(0).
\end{equation}
The condition includes symmetric $A(x)$ but also nonsymmetric
matrices satisfying one common diagonal detailed-balance relation.
\end{proposition}

The next example shows that asymmetric dispersal can separate the two
whole-space thresholds even for constant scalar coefficients.

\begin{example}\label{ex:gaussian-gap}
Let $N=m=1$, $A(x)\equiv a$, and
\[
J(z)=\frac1{\sqrt{2\pi}\sigma}
\exp\!\left(-\frac{(z-b)^2}{2\sigma^2}\right),
\qquad \sigma>0.
\]
The weighted periodic eigenvalue is
$\lambda_L^*(\mu)=a-1+
\exp(b\mu+\sigma^2\mu^2/2)$ for every period $L$.
Theorem~\ref{TH2.1} therefore yields the exact identities
\[
\lambda'=a,\qquad
\lambda^*=a-1+\exp\!\left(-\frac{b^2}{2\sigma^2}\right).
\]
For $b\ne0$ the inequality $\lambda^*<\lambda'$ is strict.
For instance, $b=2$, $\sigma=1$, and $a=1/2$ yield
$\lambda'=1/2>0$ and $\lambda^*=-1/2+e^{-2}<0$.
This example demonstrates why the two spectral thresholds cannot be
interchanged even in a homogeneous scalar problem.
\end{example}

\begin{proposition}\label{prop:constant-coefficient-stress-test}
Suppose that (J) holds and $A(x)\equiv A_0$ is a constant
cooperative irreducible matrix. For every $\mu\in\mathcal I$ and
for every choice of spatial periods,
\begin{equation}\label{eq:constant-coefficient-exact-eigenvalue}
\lambda_L^*(\mu)=\lambda_L'(\mu)
=s_{\rm PF}\!\left(A_0-I+
 \operatorname{diag}(M_1(\mu),\ldots,M_m(\mu))\right).
\end{equation}
Moreover, for the averaged operator in
\eqref{eq:averaged-operator} and the effective Hamiltonian in
\eqref{eq:effective-periodic-Hamiltonian}, one has the exact identities
\[
\widehat\lambda(\mu)
=\overline\lambda_\infty(\mu)
=\lambda_L^*(\mu)
\ (\mu\in\mathcal I).
\]
Thus neither the periodic eigenvalue nor either period-limit quantity
depends on $L$. In the scalar case $A_0=a$, their common value is
$a-1+M_1(\mu)$. These identities hold even when $\mu$ belongs to
the boundary of $\mathcal I$.
\end{proposition}

\begin{proof}[Proof of Proposition~\ref{prop:constant-coefficient-stress-test}]
Write $B_\mu=A_0-I+\operatorname{diag}(M_i(\mu))$.
Its positive Perron eigenvector $v$ is constant and satisfies
$\mathcal L_\mu v=B_\mu v=s_{\rm PF}(B_\mu)v$. This is
simultaneously an admissible periodic upper and lower test, so that
$\lambda_L^*(\mu)\leq s_{\rm PF}(B_\mu)\leq\lambda_L'(\mu)$.
Conversely, if $\phi\in X^{++}$ is a periodic upper test with
$\mathcal L_\mu\phi\leq\lambda\phi$, integrate each component
over a period cell. Periodicity and Fubini's theorem yield
$B_\mu\bar\phi\leq\lambda\bar\phi$, where
$\bar\phi=|\mathcal C|^{-1}\int_{\mathcal C}\phi(x)\,\mathrm dx\gg0$.
The matrix Collatz--Wielandt inequality gives
$s_{\rm PF}(B_\mu)\leq\lambda$. The same calculation for a
periodic lower test gives the reverse bound. Hence both periodic
values equal $s_{\rm PF}(B_\mu)$.
The identical integration argument applies to
$\widehat{\mathcal L}_\mu$, because its averaging term has
mean $\operatorname{diag}(M_i(\mu))\bar\phi$.
Finally, the Hamiltonian $H(p,x)$ is independent of $x$ and
convex in $p$ under the present hypotheses. For every
$W\in C^1_{\rm per}(\mathbb T^N)$, Jensen's inequality and
$\int_{\mathbb T^N}\nabla W=0$ give
\[
\sup_x H(\mu+\nabla W(x))
\geq\int_{\mathbb T^N}H(\mu+\nabla W(x))\,\mathrm dx
\geq H(\mu).
\]
This remains valid for the extended-valued Hamiltonian when
one of the terms is infinite. Taking $W=0$ attains the lower
bound, proving the effective-Hamiltonian identity directly.
\end{proof}

\begin{example}\label{ex:one-sided-moment-domain}
The origin need not be an interior point of $\mathcal I$, even for
strictly positive continuous kernels satisfying (J). In one dimension,
fix $a>0$ and choose a continuous nonnegative cutoff $\chi$ with
$\chi(r)=0$ for $r\leq1$ and $\chi(r)=1$ for $r\geq2$.
After normalization to unit mass, consider
\[
J(z)=c\left[e^{-z^2}
+\boldsymbol1_{\{z>0\}}\chi(z)e^{-az}(1+z)^{-2}
+\boldsymbol1_{\{z<0\}}\chi(-z)(1+|z|)^{-4}\right],
\]
where $c>0$ is the normalizing constant. Both one-sided terms
vanish near the origin, so $J$ is continuous and $J(0)>0$.
Direct inspection of the two tails gives
\[
\mathcal I=[0,a],\qquad
M(a)<\infty,\qquad
\lim_{\mu\uparrow a}M'(\mu)=+\infty.
\]
Thus the moment domain can have a finite boundary value without
first-derivative regularity there. In particular, assumptions (J)
and (A) alone cannot justify differentiation of weighted spectra
at a boundary weight. All derivative expansions in the present paper
are therefore stated under their respective interior or explicit
exponential-moment hypotheses.
\end{example}

\subsection{Proofs}\label{subsec:spectral-approximation-proofs}

After coefficient approximation, we prove the measure classification,
the whole-space characterization, and the detailed-balance criterion.

\subsubsection{Coefficient approximation}

We record the periodic weighted form of the approximation
construction in Wang and Zhang \cite[Theorem~A and its proof]{WZ-1}.

	\begin{lemma}\label{lem:approx}
		Assume that (J) and (A) hold. Then there exist families of matrices
		\[
		\underline{A}^{\varepsilon}(x) = \bigl( \underline{a}_{i k}^{\varepsilon}(x) \bigr)_{m \times m},
		\qquad
		\overline{A}^{\varepsilon}(x) = \bigl( \overline{a}_{i k}^{\varepsilon}(x) \bigr)_{m \times m},
		\]
		defined for $\varepsilon>0$ and independent of $\mu$, with the following properties:
		\begin{enumerate}
		\item[(i)] The entries are continuous and $\mathbb L$-periodic, and the matrices are cooperative and irreducible. The Perron spectral bounds of $\underline A^\varepsilon-I$ and $\overline A^\varepsilon-I$ attain their maxima on nonempty open sets. Moreover, $A-2\varepsilon I\leq\underline A^\varepsilon\leq A-\varepsilon I$ and $A+\varepsilon I\leq\overline A^\varepsilon\leq A+2\varepsilon I$. The entries of these two matrices are nonincreasing and nondecreasing in $\varepsilon$, respectively, and
		\[
		\underline{a}_{i k}^{\varepsilon} \leq a_{i k} \leq \overline{a}_{i k}^{\varepsilon} \quad \text{in } \mathbb{R}^N, \qquad \lim_{\varepsilon \rightarrow 0^{+}} \underline{a}_{i k}^{\varepsilon} = \lim_{\varepsilon \rightarrow 0^{+}} \overline{a}_{i k}^{\varepsilon} = a_{i k} \quad \text{uniformly}.
		\]
			\item[(ii)] For every $\mu\in\mathcal I$, the operators $\underline{\mathcal L}_{\mu}^{\varepsilon}$ and $\overline{\mathcal L}_{\mu}^{\varepsilon}$ (defined by \eqref{Q2} with $A$ replaced by 
			$\underline{A}^{\varepsilon}$ and $\overline A^\varepsilon$, respectively) admit principal eigenvalues $\lambda_\varepsilon^-(\mu)$ and $\lambda_\varepsilon^+(\mu)$, respectively. The fixed spatial periods are suppressed in the notation $\lambda_\varepsilon^\pm$.
			\item[(iii)] For every $\mu\in\mathcal I$, $\lambda_\varepsilon^-(\mu) \leq \lambda_\varepsilon^+(\mu)$; $\lambda_\varepsilon^-(\mu)$ and $\lambda_\varepsilon^+(\mu)$ are strictly decreasing and increasing in $\varepsilon$, respectively, and
			\[
			\lim_{\varepsilon \rightarrow 0^{+}} \lambda_\varepsilon^-(\mu) = \lim_{\varepsilon \rightarrow 0^{+}} \lambda_\varepsilon^+(\mu) = \lambda_L^*(\mu).
			\]
			More precisely, for every $\mu\in\mathcal I$,
			\[
			0\leq\lambda_L^*(\mu)-\lambda_\varepsilon^-(\mu)\leq2\varepsilon,
			\qquad
			0\leq\lambda_\varepsilon^+(\mu)-\lambda_L^*(\mu)\leq2\varepsilon.
			\]
		\end{enumerate}
		For $\mu\notin\mathcal I$, we set $\lambda_\varepsilon^\pm(\mu)=+\infty$.
	\end{lemma}

\begin{proof}[Proof of Lemma~\ref{lem:approx}]
Using the cited construction with multiplication matrix $A-I$, set
\[
h_\varepsilon(x)=(s(x)-s^*+\varepsilon)_+,
\quad
\underline A^\varepsilon=A-(\varepsilon+h_\varepsilon)I,
\quad
\overline A^\varepsilon=A+(2\varepsilon-h_\varepsilon)I.
\]
Here $r_+=\max\{r,0\}$. This construction preserves periodicity,
leaves the off-diagonal entries unchanged, and is independent of $\mu$.
The Perron maxima are attained on the nonempty open set
$\{s>s^*-\varepsilon\}$, so the weighted version of Lemma~\ref{L11}(ii) gives
principal eigenpairs for every $\mu\in\mathcal I$.
The stated matrix bounds and monotonicity properties follow as in
the cited proof. Monotonicity and the scalar-shift identity for the
spectral value then give the two $2\varepsilon$ error bounds,
uniformly in $\mu$.
\end{proof}

We also record the comparison bounds following from
Lemmas~\ref{L11} and~\ref{lem:approx}, valid for every $\mu\in\mathcal I$:
\begin{equation}\label{eq:weighted-CW-sandwich}
\lambda_\varepsilon^-(\mu)
\leq\lambda_L'(\mu)
\leq\lambda_L^*(\mu)
\leq\lambda_\varepsilon^+(\mu).
\end{equation}

\subsubsection{Periodic eigenmeasures}
\label{subsec:measure-classification-proof}

We first establish the spectral properties of $\mathcal R$ in
\eqref{eq:appendix-K}, then use them to classify the periodic eigenmeasures.

\begin{lemma}\label{lem:appendix-K-principal}
Assume that (J) and (A) hold. Suppose that
\eqref{eq:appendix-integrable-case} holds.
Then the operator $\mathcal R$ in \eqref{eq:appendix-K} has a simple principal eigenvalue
$\gamma_1=r(\mathcal R)>0$ and an eigenfunction
$\boldsymbol g_1\in X^{++}$.
\end{lemma}

\begin{proof}[Proof of Lemma~\ref{lem:appendix-K-principal}]
For $x\notin\Omega_0$, $C(x)$ is an irreducible nonsingular $M$-matrix and
$C(x)^{-1}\gg0$. Here a nonsingular $M$-matrix is a real matrix with
nonpositive off-diagonal entries and a nonnegative inverse.
Simplicity of the Perron root implies that $\boldsymbol v$ is continuous.
Thus $v_*:=\min_{i,x}v_i(x)>0$. Since
$\boldsymbol1\leq v_*^{-1}\boldsymbol v(x)$ and $C(x)^{-1}$ is positive,
the Perron eigenvector identity gives, for $x\notin\Omega_0$,
\begin{equation}\label{eq:appendix-resolvent-integrable-bound}
C(x)^{-1}\boldsymbol v(x)
=\frac{\boldsymbol v(x)}{s^*-s(x)},
\qquad
\|C(x)^{-1}\|_\infty
\leq\frac{1}{v_*(s^*-s(x))}.
\end{equation}
Assumption \eqref{eq:appendix-integrable-case} therefore
implies $c_{ij}\in L^1(\mathbb T^N)$. 

By (J) and the boundedness of $\mathcal C$, there exists
$C_{\mathcal C}>0$ such that
\begin{equation*}
\sum_{p\in\mathbb L}\sup_{z\in\overline{\mathcal C}}J_i(z+p)
\leq C_{\mathcal C}\sum_{p\in\mathbb L}(1+|p|)^{-N-\delta}<\infty.
\end{equation*}
Hence the series
\[
J_i^{\rm per}(z):=\sum_{p\in\mathbb L}J_i(z+p)
\]
converges uniformly on $\overline{\mathcal C}$ and defines a continuous
kernel on $\mathbb T^N$.
Let $\boldsymbol0\ne\boldsymbol g\in X^+$.  Since
$|\Omega_0|=0$ under \eqref{eq:appendix-integrable-case}, some component of
$\boldsymbol g$ is positive on a ball contained in
$\mathbb T^N\setminus\Omega_0$.  On that ball,
$C^{-1}\boldsymbol g\gg0$.  By (J), there are $r,j_0>0$ such that
$J_i(z)\ge j_0$ for $|z|\le r$ and every $i$.  Consequently one application
of $\mathcal R$ makes every component positive on the $r$-neighborhood of
that ball.  Iterating this argument, and using compactness of
$\mathbb T^N$, gives an integer $n_0$ such that
\begin{equation}\label{eq:appendix-eventual-strong-positivity}
\mathcal R^{n_0}\boldsymbol g\gg0
\qquad(\boldsymbol0\ne\boldsymbol g\in X^+).
\end{equation}
Moreover, convolution of the $L^1$ functions
$c_{ij}g_j$ with the continuous kernels $J_i^{\rm per}$ is
equicontinuous;
hence $\mathcal R$ is compact on $X$. Applying the Krein--Rutman
theorem to the strongly positive power in
\eqref{eq:appendix-eventual-strong-positivity}
(see \cite[Section~19]{Deimling1985}) yields the assertion.
\end{proof}

\begin{proof}[\textbf{Proof of Theorem~\ref{thm:appendix-measure-classification}}]
We first prove the classification under \eqref{eq:appendix-integrable-case}.
For $h>0$, define
\begin{equation*}
C_h(x):=(s^*+h+1)I-A(x),
\qquad
\mathcal R_h\boldsymbol g:=\mathcal J*(C_h^{-1}\boldsymbol g).
\end{equation*}
Since $C_h^{-1}$ is continuous and positive, the compactness and
eventual strong positivity argument in
Lemma~\ref{lem:appendix-K-principal} also applies to $\mathcal R_h$.
Denote its principal eigenvalue by $\gamma_1(h)$ and choose an
eigenfunction $\boldsymbol g_h\in X^{++}$.
If $0<h_1<h_2$, the resolvent identity for the irreducible
$M$-matrices $C_h(x)$ gives
\[
C_{h_1}^{-1}(x)>C_{h_2}^{-1}(x)
\quad(x\in\mathbb T^N).
\]
Hence positivity and irreducibility imply
\begin{equation}\label{eq:appendix-gamma-h-monotone}
\gamma_1>\gamma_1(h_1)>\gamma_1(h_2)
\qquad(0<h_1<h_2).
\end{equation}
For each fixed $h_0>0$, the resolvent identity
\[
C_h^{-1}-C_{h_0}^{-1}
=(h_0-h)C_h^{-1}C_{h_0}^{-1}
\]
and the uniform boundedness of these inverses for $h$ near $h_0$ give
$\|\mathcal R_h-\mathcal R_{h_0}\|\to0$ as $h\to h_0$.
Applying the Collatz--Wielandt inequalities to $\boldsymbol g_{h_0}$ yields
\[
\min_{i,x}\frac{(\mathcal R_h\boldsymbol g_{h_0})_i(x)}
{g_{h_0,i}(x)}
\leq\gamma_1(h)\leq
\max_{i,x}\frac{(\mathcal R_h\boldsymbol g_{h_0})_i(x)}
{g_{h_0,i}(x)}.
\]
Both bounds tend to $\gamma_1(h_0)$, because
$\min_{i,x}g_{h_0,i}(x)>0$. Thus $h\mapsto\gamma_1(h)$ is
continuous on $(0,\infty)$.
For almost every $x$, the matrices $C_h(x)^{-1}$ converge to $C(x)^{-1}$
as $h\downarrow0$, and $0\leq C_h(x)^{-1}\leq C(x)^{-1}$ entrywise.
Estimate~\eqref{eq:appendix-resolvent-integrable-bound},
assumption~\eqref{eq:appendix-integrable-case}, and dominated convergence give
\[
\max_i\sum_j\int_{\mathbb T^N}
\bigl|(C_h(x)^{-1}-C(x)^{-1})_{ij}\bigr|\,\mathrm dx
\longrightarrow0.
\]
Writing the convolutions with the periodized kernels gives
\[
\|\mathcal R_h-\mathcal R\|
\leq \max_i\|J_i^{\rm per}\|_\infty
\max_i\sum_j\int_{\mathbb T^N}
\bigl|(C_h(x)^{-1}-C(x)^{-1})_{ij}\bigr|\,\mathrm dx
\longrightarrow0.
\]
Since $\boldsymbol g_1\gg0$, the Collatz--Wielandt inequalities applied
to $\mathcal R_h\boldsymbol g_1$ imply
$\gamma_1(h)\to\gamma_1$. On the other hand,
$\sup_x\|C_h(x)^{-1}\|=O(h^{-1})$ as $h\to\infty$, so that
$\|\mathcal R_h\|\to0$. Consequently,
\begin{equation}\label{eq:appendix-gamma-h-limits}
\gamma_1(h)\longrightarrow\gamma_1\quad(h\downarrow0),
\qquad
\gamma_1(h)\longrightarrow0\quad(h\to\infty).
\end{equation}

If $\mathcal R_h\boldsymbol g=\boldsymbol g$ with
$\boldsymbol g\in X^{++}$ and
$\boldsymbol\phi=C_h^{-1}\boldsymbol g$, then
\[
\mathcal J*\boldsymbol\phi=C_h\boldsymbol\phi
\quad\Longleftrightarrow\quad
\mathcal L\boldsymbol\phi=(s^*+h)\boldsymbol\phi.
\]
Conversely, for a positive principal eigenfunction $\boldsymbol\phi$,
setting $\boldsymbol g=C_h\boldsymbol\phi=\mathcal J*\boldsymbol\phi$
gives $\mathcal R_h\boldsymbol g=\boldsymbol g$. We have therefore proved
\begin{equation}\label{eq:appendix-gamma-h-one-equivalence}
\gamma_1(h)=1
\quad\Longleftrightarrow\quad
s^*+h\ \text{is a principal eigenvalue of }\mathcal L.
\end{equation}
If $\lambda_L^*>s^*$, Lemma~\ref{L11}(i) and
\eqref{eq:appendix-gamma-h-one-equivalence}, with
$h=\lambda_L^*-s^*$, give $\gamma_1(h)=1$.  Equations
\eqref{eq:appendix-gamma-h-monotone}--\eqref{eq:appendix-gamma-h-limits}
then imply $\gamma_1>1$. Conversely, if $\gamma_1>1$, continuity,
strict monotonicity, and the endpoint limits give a unique $h>0$
such that $\gamma_1(h)=1$.
By \eqref{eq:appendix-gamma-h-one-equivalence},
$s^*+h$ is a principal eigenvalue, and hence
$\lambda_L^*=s^*+h>s^*$. This proves part~(i).
The lower bound \eqref{eq:periodic-multiplication-edge}, together
with part~(i), shows that $\gamma_1\leq1$ if and only if
$\lambda_L^*=s^*$.

For $\lambda_L^*=s^*$, let
$\mathrm d\boldsymbol\nu=\boldsymbol f\,\mathrm dx
+\mathrm d\boldsymbol\eta$ be a nonnegative solution.  Because
$\mathcal J[\boldsymbol\nu](x)\,\mathrm dx$ is absolutely continuous, the
singular part of \eqref{eq:1.3} is
\begin{equation}\label{eq:appendix-singular-equation}
C(x)\,\mathrm d\boldsymbol\eta(x)=\boldsymbol0.
\end{equation}
Pointwise irreducibility of $A(x)$ implies that
$\ker C(x)=\operatorname{span}\{\boldsymbol v(x)\}$ on $\Omega_0$ and is
trivial off $\Omega_0$. More explicitly, set $\theta:=\sum_i\eta_i$ and
$\boldsymbol h:=\mathrm d\boldsymbol\eta/\mathrm d\theta$, where
$h_i=\mathrm d\eta_i/\mathrm d\theta$ is the Radon--Nikodym density of
$\eta_i$ with respect to $\theta$, so that
$\eta_i(E)=\int_Eh_i\,\mathrm d\theta$ for every Borel set $E$. Then
$\boldsymbol h\geq\boldsymbol0$, $\sum_i h_i=1$, and
$C(x)\boldsymbol h(x)=\boldsymbol0$ for $\theta$-almost every $x$.
Invertibility off $\Omega_0$ implies
$\theta(\mathbb T^N\setminus\Omega_0)=0$, and normalization of the Perron
vector implies $\boldsymbol h=\boldsymbol v$ $\theta$-almost everywhere.
Since $\Omega_0$ is closed, $\theta$ is supported on $\Omega_0$. Hence
$\mathrm d\boldsymbol\eta=\boldsymbol v\,\mathrm d\theta$ for this unique
finite positive measure $\theta$. Since the periodized kernels are
continuous and $\boldsymbol\nu$ is finite,
$\mathcal J[\boldsymbol\nu]\in X$. The absolutely continuous part of
\eqref{eq:1.3} gives
\[
C\boldsymbol f=\mathcal J[\boldsymbol\nu]
=\mathcal J*\boldsymbol f+\boldsymbol F_\theta
\quad\text{a.e. on }\mathbb T^N.
\]
We henceforth use this continuous representative of $C\boldsymbol f$.
Since $|\Omega_0|=0$ and $C^{-1}(C\boldsymbol f)=\boldsymbol f$
almost everywhere, the preceding identity becomes, in $X$,
\begin{equation}\label{eq:appendix-ac-equation}
(I-\mathcal R)(C\boldsymbol f)=\boldsymbol F_\theta.
\end{equation}
If $\gamma_1=1$, let
$\mathcal R\boldsymbol g_1=\boldsymbol g_1$ and set
$\boldsymbol f_1=C^{-1}\boldsymbol g_1\in L^1(\mathbb T^N,\mathbb R^m)$. Then
\begin{equation*}
\mathcal J*\boldsymbol f_1-C\boldsymbol f_1
=\mathcal R\boldsymbol g_1-\boldsymbol g_1=\boldsymbol0,
\end{equation*}
so $\boldsymbol f_1\,\mathrm dx$ solves \eqref{eq:1.3}. Conversely, a
nonzero nonnegative absolutely continuous solution satisfies
$\mathcal R(C\boldsymbol f)=C\boldsymbol f$ in $X$ by
\eqref{eq:appendix-ac-equation} with $\theta=0$.
This fixed point is nonzero and nonnegative. Equation
\eqref{eq:appendix-eventual-strong-positivity} makes it strictly positive,
so its eigenvalue $1$ equals $\gamma_1$. The simplicity in
Lemma~\ref{lem:appendix-K-principal} gives uniqueness up to multiplication by a positive constant.
Let $\ell$ be a nonzero positive eigenfunctional of the adjoint operator
$\mathcal R^*$ at the eigenvalue $1$, and write
$\langle\ell,\boldsymbol g\rangle=\ell(\boldsymbol g)$ for its action
on $\boldsymbol g\in X$. For every $\boldsymbol0\ne\boldsymbol g\in X^+$,
\eqref{eq:appendix-eventual-strong-positivity} gives
\[
\langle\ell,\boldsymbol g\rangle
=\langle\ell,\mathcal R^{n_0}\boldsymbol g\rangle>0.
\]
Thus $\ell$ is strictly positive.
If a solution had singular part
$\boldsymbol v\,\mathrm d\theta\ne0$, pairing
\eqref{eq:appendix-ac-equation} with $\ell$ gives
$0=\langle\ell,\boldsymbol F_\theta\rangle$. On the other hand,
\[
\int_{\mathbb T^N}(\boldsymbol F_\theta)_i(x)\,\mathrm dx
=\int_{\Omega_0}v_i(x)\,\mathrm d\theta(x)>0.
\]
Thus $\boldsymbol F_\theta\geq\boldsymbol0$ is nonzero, and the strict
positivity of $\ell$ gives
$\langle\ell,\boldsymbol F_\theta\rangle>0$, a contradiction.  This proves
part~(ii) of the theorem.

Assume finally that $\gamma_1<1$. The positive resolvent series
\[
\boldsymbol g_\theta=(I-\mathcal R)^{-1}\boldsymbol F_\theta
=\sum_{n=0}^\infty\mathcal R^n\boldsymbol F_\theta
\]
converges in $X$. Since
$C\boldsymbol f_\theta=\boldsymbol g_\theta$ almost everywhere,
$(I-\mathcal R)\boldsymbol g_\theta=\boldsymbol F_\theta$, and
\eqref{eq:appendix-Perron-null-vector} holds, we obtain
\begin{align*}
\mathcal J[\boldsymbol\nu_\theta]-C\boldsymbol f_\theta
&=\boldsymbol F_\theta
+\mathcal R\boldsymbol g_\theta-\boldsymbol g_\theta
=\boldsymbol0,
\\
C(x)\,\mathrm d\boldsymbol\eta_\theta(x)&=\boldsymbol0
\end{align*}
Thus
\eqref{eq:singular-eigenmeasure-representation} solves \eqref{eq:1.3}.

Conversely, let $\mathrm d\boldsymbol\nu=\boldsymbol f\,\mathrm dx
+\mathrm d\boldsymbol\eta_\theta$ be a nonnegative solution.
Since $r(\mathcal R)<1$, equation \eqref{eq:appendix-ac-equation} has the
unique solution $C\boldsymbol f=\boldsymbol g_\theta$, and hence
$\boldsymbol f=\boldsymbol f_\theta$. If $\theta=0$, then
$\mathcal R(C\boldsymbol f)=C\boldsymbol f$, contradicting
$\gamma_1<1$ unless $\boldsymbol f=\boldsymbol0$.  Therefore every nonzero
solution has a nontrivial singular part, and the converse follows. This proves part~(iii).

Finally, suppose that $\lambda_L^*=s^*$. Lemma~\ref{L11}(i)--(ii)
gives \eqref{eq:appendix-integrable-case}, and part~(i) gives
$\gamma_1\leq1$. If $\gamma_1=1$, the measure
$\boldsymbol f_1\,\mathrm dx$ constructed in part~(ii) is a nonzero
finite nonnegative solution. If $\gamma_1<1$, choose $x_0\in\Omega_0$
and let $\theta=\delta_{x_0}$ be the unit point mass at $x_0$.
Formula~\eqref{eq:singular-eigenmeasure-representation}, with the quantities in
\eqref{eq:singular-eigenmeasure-construction}, gives a nonzero finite
nonnegative solution with singular part $\boldsymbol v(x_0)\delta_{x_0}$.
For every nonnegative solution,
\eqref{eq:appendix-singular-equation} and the measure decomposition above
give \eqref{eq:singular-support-Omega0}. This proves the existence and
support assertions and completes the proof.
\end{proof}

\subsubsection{Whole-space eigenvalues}

On the space
$C_{\rm b}(\mathbb R^N,\mathbb R^m)$ of bounded continuous functions,
the positive semigroup generated by $\mathcal L$ is not strongly
positive: a nonzero nonnegative datum need not become uniformly
positive in every component. Neither $\mathcal L$ nor its resolvent
$(\alpha I-\mathcal L)^{-1}$ for sufficiently large $\alpha$ is compact.
Unlike periodic tests, positive upper tests and bounded positive lower
tests may yield different spectral values. Wu et~al.\ \cite{WFLS}
established a Harnack inequality and compared generalized principal
eigenvalues for compactly supported kernels; here compact support is
not required.
We prove the three assertions of Theorem~\ref{TH2.1} in order.

	\begin{proof}[\textbf{Proof of Theorem~\ref{TH2.1}(i)}]
Take $\mu_1,\mu_2\in\mathcal I$, $r\in(0,1)$, and set
$\mu_r=r\mu_1+(1-r)\mu_2$. For any $\epsilon>0$, there exist
$\boldsymbol\phi,\boldsymbol\psi\in X^{++}$ such that
\[
\int_{\mathbb R^N}J_i(x-y)e^{\mu_1\cdot(x-y)}\phi_i(y)\,\mathrm dy
-\phi_i(x)+\sum_{j=1}^m a_{ij}(x)\phi_j(x)
\leq(\lambda_L^*(\mu_1)+\epsilon)\phi_i(x),
\]
and
\[
\int_{\mathbb R^N}J_i(x-y)e^{\mu_2\cdot(x-y)}\psi_i(y)\,\mathrm dy
-\psi_i(x)+\sum_{j=1}^m a_{ij}(x)\psi_j(x)
\leq(\lambda_L^*(\mu_2)+\epsilon)\psi_i(x).
\]
Define $v_i:=\phi_i^r\psi_i^{1-r}$. By H\"older's inequality,
$\mu_r\in\mathcal I$. One can compute
\[
\begin{aligned}
\frac{(\mathcal J_{\mu_r}\boldsymbol v)_i}{v_i}
&=\int_{\mathbb R^N}J_i(z)
\left[e^{\mu_1\cdot z}\frac{\phi_i(x-z)}{\phi_i(x)}\right]^r
\left[e^{\mu_2\cdot z}\frac{\psi_i(x-z)}{\psi_i(x)}\right]^{1-r}
\,\mathrm dz\\
&\leq
\left[\int_{\mathbb R^N}J_i(z)e^{\mu_1\cdot z}
\frac{\phi_i(x-z)}{\phi_i(x)}\,\mathrm dz\right]^r
\left[\int_{\mathbb R^N}J_i(z)e^{\mu_2\cdot z}
\frac{\psi_i(x-z)}{\psi_i(x)}\,\mathrm dz\right]^{1-r}\\
&\leq r\frac{(\mathcal J_{\mu_1}\boldsymbol\phi)_i}{\phi_i}
+(1-r)\frac{(\mathcal J_{\mu_2}\boldsymbol\psi)_i}{\psi_i}.
\end{aligned}
\]
The same arithmetic--geometric mean inequality bounds
$v_j/v_i$ by $r\phi_j/\phi_i+(1-r)\psi_j/\psi_i$ for $j\ne i$.
Since $a_{ij}\geq0$ for $j\ne i$ and the diagonal term is unchanged,
\[
\frac{(\mathcal L_{\mu_r}\boldsymbol v)_i}{v_i}
\leq r\frac{(\mathcal L_{\mu_1}\boldsymbol\phi)_i}{\phi_i}
+(1-r)\frac{(\mathcal L_{\mu_2}\boldsymbol\psi)_i}{\psi_i}
\leq r\lambda_L^*(\mu_1)+(1-r)\lambda_L^*(\mu_2)+\epsilon.
\]
Letting $\epsilon\to0$, we obtain
\[
\lambda_L^*(r\mu_1+(1-r)\mu_2)
\leq r\lambda_L^*(\mu_1)+(1-r)\lambda_L^*(\mu_2).
\]
Thus $\lambda_L^*$ is convex on $\mathcal I$, and hence on
$\mathbb R^N$ when extended by $+\infty$ outside $\mathcal I$.

		Note that evaluating each component of a positive periodic super-solution at a minimum point of that component gives
		\begin{equation}\label{eq:weighted-moment-lower-bound}
		\lambda_L^*(\mu)\geq
		\max_{i\in\Lambda}\{M_i(\mu)-1+\min_xa_{ii}(x)\}.
		\end{equation}
		Lemma~\ref{lem:uniform-weight-coercivity} gives the stronger
uniform exponential estimate
\eqref{eq:uniform-exponential-weight-coercivity}, and in particular
$\lambda_L^*(\mu)\to+\infty$ as $|\mu|\to\infty$ in the
extended sense.

		We also verify lower semicontinuity, including at the boundary of $\mathcal I$. Fix $\boldsymbol v\in X^{++}$. For sequences $\mu_n\to\mu\in\mathcal I$ and $x_n\to x$ on the torus, Fatou's lemma applied to
		$J_i(z)e^{\mu_n\cdot z}v_i(x_n-z)/v_i(x_n)$ shows that
		\[
		(\mu,x)\longmapsto
		\frac{(\mathcal L_\mu\boldsymbol v)_i(x)}{v_i(x)}
		\]
		is jointly lower semicontinuous on $\mathcal I\times\mathbb T^N$. Taking the minimum over the compact set $\mathbb T^N\times\Lambda$, and then the supremum over $\boldsymbol v\in X^{++}$, preserves lower semicontinuity. The equality of the weighted Collatz--Wielandt values established by \eqref{eq:weighted-CW-sandwich} therefore gives lower semicontinuity of $\lambda_L^*$ at every point of $\mathcal I$. If $\mu\notin\mathcal I$, at least one $M_i(\mu)$ is infinite; Fatou's lemma and \eqref{eq:weighted-moment-lower-bound} then give $\liminf_{n\to\infty}\lambda_L^*(\mu_n)=+\infty$ whenever $\mu_n\to\mu$. Thus the extended-valued function is lower semicontinuous on $\mathbb R^N$. Since its value at zero is finite and it is coercive, it attains its minimum at a point of $\mathcal I$.
	
We next prove uniqueness of the minimum point. Assume that
\(1/(s^*-s(x))\notin L^1(\mathbb T^N)\). By the weighted version of
Lemma~\ref{L11}(ii),
	every exponentially weighted operator with
	\(\mu\in\mathcal I\) has a positive periodic principal eigenfunction.
We first prove uniqueness when the minimum is attained in
$\operatorname{int}\mathcal I$ by computing the Hessian.
For $\mu\in\operatorname{int}\mathcal I$, the operator
family is analytic in norm: the local exponential-moment bounds
in Lemma~\ref{lem:interior-moment-regularity} control every
Taylor coefficient of the weighted kernels. Its principal eigenvalue
is isolated and simple, so analytic perturbation theory
\cite[Chapter~VII]{Kato1995} applies.
Fix $q\in\mathbb R^N\setminus\{0\}$ and choose a locally analytic
positive eigenfunction $\boldsymbol p_t$ at weight $\mu+tq$.
Write $\lambda(t)=\lambda_L^*(\mu+tq)$ and $\boldsymbol p=\boldsymbol p_0$.
Let $\boldsymbol p^*$ be a positive adjoint eigenfunction, normalized by
\[
\sum_i\int_{\mathbb T^N}p_i(x)p_i^*(x)\,\mathrm dx=1.
\]
An $L^1$ density suffices here: the positive dual eigenmeasure has such
a density because its equation is
$\boldsymbol\nu=[(\lambda(0)+1)I-A^{\mathsf T}]^{-1}
\mathcal J_\mu^*\boldsymbol\nu$, and convolution of a finite measure
with the periodized integrable kernels is in $L^1$.
Positivity near the origin and irreducibility give $p_i^*>0$ almost
everywhere. Define
\[
\chi_i=\left.\partial_t\log p_{t,i}\right|_{t=0},\qquad
\eta_i=\left.\partial_t^2\log p_{t,i}\right|_{t=0}.
\]
Multiplying the eigenfunction equation by $p_i^*\eta_i$ and the
adjoint equation by $p_i\eta_i$, then integrating and using periodicity,
we obtain
\begin{equation}\label{eq:weighted-eigenfunction-cancellation}
\begin{aligned}
0={}&\sum_i\int_{\mathbb T^N}\int_{\mathbb R^N}
p_i^*(x)J_i(z)e^{\mu\cdot z}p_i(x-z)
[\eta_i(x-z)-\eta_i(x)]\,\mathrm dz\,\mathrm dx\\
&+\sum_{i\ne j}\int_{\mathbb T^N}
p_i^*(x)a_{ij}(x)p_j(x)[\eta_j(x)-\eta_i(x)]\,\mathrm dx.
\end{aligned}
\end{equation}
Dividing the equation for $\boldsymbol p_t$ by $p_{t,i}$ and
differentiating twice at $t=0$ yields
\[
\begin{aligned}
\lambda''(0)
={}&\int_{\mathbb R^N}J_i(z)e^{\mu\cdot z}\frac{p_i(x-z)}{p_i(x)}
[\eta_i(x-z)-\eta_i(x)]\,\mathrm dz\\
&+\sum_{j\ne i}a_{ij}(x)\frac{p_j(x)}{p_i(x)}
[\eta_j(x)-\eta_i(x)]\\
&+\int_{\mathbb R^N}J_i(z)e^{\mu\cdot z}\frac{p_i(x-z)}{p_i(x)}
[q\cdot z+\chi_i(x-z)-\chi_i(x)]^2\,\mathrm dz\\
&+\sum_{j\ne i}a_{ij}(x)\frac{p_j(x)}{p_i(x)}
[\chi_j(x)-\chi_i(x)]^2.
\end{aligned}
\]
Multiplication by $p_i p_i^*$, integration and summation, followed by
\eqref{eq:weighted-eigenfunction-cancellation}, give the square identity
\begin{equation*}
\begin{aligned}
q^{\mathsf T}D_\mu^2\lambda_L^*(\mu)q
={}&\sum_i\int_{\mathbb T^N}\int_{\mathbb R^N}
p_i^*(x)J_i(z)e^{\mu\cdot z}p_i(x-z)
[q\cdot z+\chi_i(x-z)-\chi_i(x)]^2\,\mathrm dz\,\mathrm dx\\
&+\sum_{i\ne j}\int_{\mathbb T^N}
p_i^*(x)a_{ij}(x)p_j(x)[\chi_j(x)-\chi_i(x)]^2\,\mathrm dx.
\end{aligned}
\end{equation*}
Both terms are nonnegative. If their sum were zero, positivity of
$J_i$ near zero and continuity would imply
$q\cdot z+\chi_i(x-z)-\chi_i(x)=0$ for all $x$ and all sufficiently
small $z$. Integrating over the torus gives $q\cdot z=0$ for all such
$z$, contradicting $q\ne0$. Thus the Hessian is positive definite
on $\operatorname{int}\mathcal I$. In particular, any minimizer in
this interior is unique, by convexity.

To cover possible boundary minimizers, or a moment domain with empty
interior, we finish with an argument that does not differentiate in
the weight.
	Let \(\mu_1\ne\mu_2\) belong to \(\mathcal I\), let
	\(r\in(0,1)\), and choose eigenfunctions
	\(\boldsymbol\phi,\boldsymbol\psi\in X^{++}\) at \(\mu_1,\mu_2\).
	Put $\mu_r=r\mu_1+(1-r)\mu_2$, $v_i=\phi_i^r\psi_i^{1-r}$, and
	$c=r\lambda_L^*(\mu_1)+(1-r)\lambda_L^*(\mu_2)$. The preceding calculation gives the nonnegative continuous periodic deficit
	\[
	\boldsymbol D:=c\boldsymbol v-\mathcal L_{\mu_r}\boldsymbol v\geq\boldsymbol0.
	\]
	We first show that $\boldsymbol D\not\equiv\boldsymbol0$. Otherwise every nonnegative deficit in the H\"older and arithmetic--geometric mean inequalities would vanish, at every point and in every component. Equality in the convolution inequality would imply
	\[
	e^{(\mu_1-\mu_2)\cdot(x-y)}
	\frac{\phi_i(y)\psi_i(x)}
	     {\psi_i(y)\phi_i(x)}
	=\text{constant in }y
	\quad\text{for a.e. }y\text{ with }J_i(x-y)>0.
	\]
	Since $J_i>0$ near the origin, continuity and $y=x$ show that this constant is one. Thus we would have
	\[
	\log\frac{\phi_i(y)}{\psi_i(y)}
	-\log\frac{\phi_i(x)}{\psi_i(x)}
	=(\mu_1-\mu_2)\cdot(y-x)
	\]
	for every $x$ and every $y$ in a fixed neighborhood of $x$. Chaining these neighborhoods gives the same identity for all $x,y$, contrary to the periodicity of $\phi_i/\psi_i$ and $\mu_1\ne\mu_2$. Therefore $\boldsymbol D\not\equiv\boldsymbol0$.

	Convexity already gives $\lambda_L^*(\mu_r)\leq c$. Suppose equality held and let $\boldsymbol h\in X^{++}$ be a principal eigenfunction at $\mu_r$. Define
	\[
	\alpha:=\min_{i,x}\frac{v_i(x)}{h_i(x)}>0,
	\qquad\boldsymbol w:=\boldsymbol v-\alpha\boldsymbol h\geq\boldsymbol0.
	\]
	Some component of $\boldsymbol w$ vanishes at a point, and
	\[
	(\mathcal L_{\mu_r}-c)\boldsymbol w=-\boldsymbol D\leq\boldsymbol0.
	\]
	At a zero of one component, its convolution term and all off-diagonal terms must vanish. Positivity of the kernel near zero propagates this zero throughout that component; irreducibility then propagates it to all components. Hence $\boldsymbol w\equiv\boldsymbol0$, which forces $\boldsymbol D\equiv\boldsymbol0$, a contradiction. We conclude that
	\[
	\lambda_L^*(r\mu_1+(1-r)\mu_2)
	<
	r\lambda_L^*(\mu_1)+(1-r)\lambda_L^*(\mu_2).
	\]
	Therefore the map $\mu\mapsto\lambda_L^*(\mu)$ is strictly convex on \(\mathcal I\), and its minimum point is unique.
\end{proof}

\begin{proof}[\textbf{Proof of Theorem~\ref{TH2.1}(ii)}]
		We prove this statement in the following three steps. By part~(i), the spectral minima below are attained.
		
\smallskip
\noindent Step 1. Compactly supported kernels under the nonintegrability condition.
		We first assume that $1/(s^*-s)\notin L^1(\mathbb T^N)$. This temporary condition will be removed in Step~3. The nonintegrability criterion in \cite[Lemma~2.1]{WZ-1}, applied on the periodic space, gives principal eigenfunctions for the periodic weighted problems for every $\mu\in\mathbb R^N$.
		Apply Wu et~al.\ \cite[Theorems~4.1 and~4.6]{WFLS} with
$\Omega=\mathbb R^N$, $D=I$, $K_i=J_i$, and their multiplication
matrix equal to $A-I$. Their strong irreducibility condition is
precisely pointwise irreducibility, which follows from (A).
The nonintegrability condition holds on a bounded domain containing
a period cell. Accounting for the opposite eigenvalue sign convention,
these results give
\[
\lambda^*=\lim_{l\to\infty}\lambda(l)>s^*,\qquad
\mathcal L\boldsymbol\psi=\lambda^*\boldsymbol\psi
\]
for some $\boldsymbol\psi\in C(\mathbb R^N,(0,\infty)^m)$.
Here one may exhaust $\mathbb R^N$ by the open cubes $(-l,l)^N$;
continuity extends the defining inequalities to their closures.

By \cite[Theorem~3.5]{WFLS} and periodicity, the local Harnack
constants are uniform under translations. Since all entries of
$[(\lambda^*+1)I-A(x)]^{-1}$ are bounded above and below by positive
constants, the eigenfunction equation also compares every component
with $\sum_j\psi_j$. More precisely,
for every $R>0$ a finite chain of the positivity balls connects
$x$ to each $y\in B_R(x)$. All constants in the resolvent
bounds and in the kernel lower bound are periodic and independent
of the chain's base point. Consequently one obtains a constant
$C_R$, depending on $R$ but not on $x$, such that
\begin{equation}\label{eq:whole-space-local-comparison}
 \sup_{y\in B_R(x),\,j\in\Lambda}\psi_j(y)
 \leq C_R\psi_i(x)\qquad(x\in\mathbb R^N,\ i\in\Lambda).
\end{equation}
Local equicontinuity for the normalized translates follows as in
the proof of \cite[Theorem~4.1]{WFLS}.

Positivity of the kernels near zero and irreducibility also give
\begin{equation}\label{eq:nonlocal-zero-propagation}
\boldsymbol w\geq\boldsymbol0,\quad
\mathcal L\boldsymbol w=\lambda\boldsymbol w,\quad
w_i(x_0)=0
\quad\Longrightarrow\quad
\boldsymbol w\equiv\boldsymbol0.
\end{equation}
The same argument applies to every spatial translate of $A$.
		
		We first prove that there exist $z_1\in\mathbb R$ and a positive continuous function $\boldsymbol\psi^1$ satisfying $\mathcal L\boldsymbol\psi^1=\lambda^*\boldsymbol\psi^1$ such that $x\mapsto e^{-z_1x_1}\boldsymbol\psi^1(x)$ is $L_1$-periodic with respect to $x_1$.
		Define the translation $\tau: x \in \mathbb{R}^N \mapsto x + L_1 e_1$, where $e_1 = (1,0,\dots,0) \in\mathbb S^{N-1}$, and denote $\boldsymbol{\psi}^\tau(x) = \boldsymbol{\psi}(\tau(x))$ and $\boldsymbol{\phi} = (\phi_i)_{i =1,\dots,m} = (\psi_i^\tau / \psi_i)_{i=1,\dots,m}$. Estimate~\eqref{eq:whole-space-local-comparison} shows that $\boldsymbol{\phi}$ is globally bounded. Let
		\[
		z_1 = L_1^{-1} \ln \left( \max_{1 \leq i \leq m} \sup_{x \in \mathbb{R}^N} \phi_i(x) \right).
		\]
		There exist $\bar{i} \in \{1,\dots,m\}$ and a sequence $(x_k)_{k \in \mathbb{N}} \subset \mathbb{R}^N$ such that $\phi_{\bar{i}}(x_k) \to \mathrm{e}^{z_1 L_1}$ as $k \to \infty$. Moreover, there exists $(y_k)_{k \in \mathbb{N}}\subset\overline{\mathcal C}$ such that for each $k$, $x_k - y_k \in L_1 \mathbb{Z} \times \cdots \times L_N \mathbb{Z}$. Up to extraction, we assume $y_k\to y_\infty\in\overline{\mathcal C}$.
		
		Define for each $k \in \mathbb{N}$,
		\[
		\begin{aligned}
\widehat{\boldsymbol{\psi}}_k(x) &= \frac{1}{\psi_{\bar{i}}(x_k)} \boldsymbol{\psi}(x + x_k), \qquad
\widehat{\boldsymbol{\psi}}_k^\tau(x) = \widehat{\boldsymbol{\psi}}_k(\tau(x)), \\
			\boldsymbol{w}_k(x) &= \mathrm{e}^{z_1 L_1} \widehat{\boldsymbol{\psi}}_k(x) - \widehat{\boldsymbol{\psi}}_k^\tau(x).
		\end{aligned}
		\]
		By \eqref{eq:whole-space-local-comparison} and periodicity, $(\widehat{\boldsymbol{\psi}}_k)_{k \in \mathbb{N}}$ is uniformly bounded on every fixed compact set and satisfies
		\[
		\mathcal J*\widehat{\boldsymbol\psi}_k(x)-\widehat{\boldsymbol\psi}_k(x)
		+A(x+x_k)\widehat{\boldsymbol\psi}_k(x)
		=\lambda^*\widehat{\boldsymbol\psi}_k(x).
		\]
		Here $A(x+x_k)=A(x+y_k)$ by periodicity. The same resolvent and convolution argument gives local equicontinuity. By a diagonal argument, passing to a subsequence if necessary,
		$\widehat{\boldsymbol\psi}_k\to\widehat{\boldsymbol\psi}_\infty$ locally
		uniformly, where
		$\widehat{\boldsymbol\psi}_\infty\in C(\mathbb R^N,[0,\infty)^m)$
		satisfies
		\[
		(\mathcal J*\widehat{\boldsymbol\psi}_\infty)(x)
		-\widehat{\boldsymbol{\psi}}_\infty(x)
		+A(x+y_{\infty})\widehat{\boldsymbol{\psi}}_\infty(x)
		= \lambda^* \widehat{\boldsymbol{\psi}}_\infty(x)
		\quad \forall x \in \mathbb{R}^N,
		\]
		that is,
		\[
		(\mathcal J*\widehat{\boldsymbol\psi}_\infty)(x-y_\infty)
		-\widehat{\boldsymbol{\psi}}_\infty(x-y_\infty)
		+A(x)\widehat{\boldsymbol{\psi}}_\infty(x-y_\infty)
		= \lambda^* \widehat{\boldsymbol{\psi}}_\infty(x-y_\infty),
		\quad \forall x \in \mathbb{R}^N.
		\]
		Moreover, $\widehat{\psi}_{\bar{i}, \infty}(0) = 1$, so $\widehat{\boldsymbol{\psi}}_\infty$ is nonzero. By \eqref{eq:nonlocal-zero-propagation}, it is positive.  
		Define $\boldsymbol{w}_\infty = \mathrm{e}^{z_1 L_1} \widehat{\boldsymbol{\psi}}_\infty - \widehat{\boldsymbol{\psi}}_\infty^\tau$. Since
		\[
		\boldsymbol{w}_k(x) = \widehat{\boldsymbol{\psi}}_k(x) \circ\left( \mathrm{e}^{z_1 L_1} \boldsymbol1 - \boldsymbol{\phi}(x + x_k) \right),
		\]
		we have $\boldsymbol{w}_\infty \geq \boldsymbol0$ with $w_{\bar{i}, \infty}(0) = 0$. As $\boldsymbol{w}_\infty$ satisfies the same equation as $\widehat{\boldsymbol{\psi}}_\infty$, implication~\eqref{eq:nonlocal-zero-propagation} gives $\boldsymbol{w}_\infty = \boldsymbol0$. Hence $\mathrm{e}^{z_1 L_1} \widehat{\boldsymbol{\psi}}_\infty = \widehat{\boldsymbol{\psi}}_\infty^\tau$.  
		Now, $\boldsymbol{\psi}^1(x) = \widehat{\boldsymbol{\psi}}_\infty(x - y_\infty)$ is positive, satisfies $\mathcal{L} \boldsymbol{\psi}^1 = \lambda^* \boldsymbol{\psi}^1$, and $x \mapsto \mathrm{e}^{-z_1 x_1} \boldsymbol{\psi}^1(x)$ is $L_1$-periodic in $x_1$:
		\[
		\begin{aligned}
			\mathrm{e}^{-z_1(x_1 + L_1)} \boldsymbol{\psi}^1(x + L_1 e_1) 
			&= \mathrm{e}^{-z_1 x_1} \mathrm{e}^{-z_1 L_1} \boldsymbol{\psi}^1(\tau(x)) \\
			&= \mathrm{e}^{-z_1 x_1} \mathrm{e}^{-z_1 L_1} \widehat{\boldsymbol{\psi}}_\infty^\tau(x - y_\infty) \\
&= \mathrm{e}^{-z_1 x_1} \widehat{\boldsymbol{\psi}}_\infty(x - y_\infty)
= \mathrm{e}^{-z_1 x_1} \boldsymbol{\psi}^1(x).
		\end{aligned}
		\]	
		We verify the induction over the coordinate directions. Suppose
$\boldsymbol\psi^{k-1}>0$ satisfies the eigenvalue equation and
\[
 \boldsymbol\psi^{k-1}(x+L_j e_j)
   =e^{z_jL_j}\boldsymbol\psi^{k-1}(x),\qquad 1\le j<k.
\]
The local comparison estimates just proved remain valid for this
solution, since they depend only on the coefficients and the spectral
gap. Apply the translation-ratio construction to the $k$th coordinate.
For every translation $a$ and every $j<k$, the functions
$\boldsymbol\psi^{k-1}(\cdot+a)$ and their positive scalar
normalizations satisfy precisely the same $j$th quasi-periodicity
identity. Local uniform convergence preserves that identity. Thus the
limiting eigenfunction obtained at the $k$th step satisfies all the
old identities and the new one
$\boldsymbol\psi^k(x+L_ke_k)=e^{z_kL_k}\boldsymbol\psi^k(x)$.
The shifts used to return the translated coefficients to the original
periodic phase do not affect any of these identities, because spatial
translations commute. Induction gives $z=(z_1,\ldots,z_N)$ and a
strictly positive eigenfunction $\boldsymbol\psi^N$ satisfying
\[
 \mathcal L\boldsymbol\psi^N=\lambda^*\boldsymbol\psi^N,
 \qquad e^{-z\cdot x}\boldsymbol\psi^N(x)\in X^{++}.
\]
In particular, the argument does not require the maximizing
translation ratios to attain their suprema on the noncompact space.
		For $\mu\in\mathcal I$, let $\boldsymbol\phi_\mu$ be a positive periodic principal eigenfunction of $\mathcal L_\mu$. By direct substitution, the periodic factor solves the eigenvalue problem with exponential weight $-z$. Hence uniqueness of the principal eigenpair gives
		\[
		\lambda^*=\lambda_L^*(-z),\qquad
		e^{-z\cdot x}\boldsymbol{\psi}^N
		\in\operatorname{span}(\boldsymbol{\phi}_{-z}).
		\]
		
		Since $e^{-z\cdot x}\boldsymbol\psi^N\in\operatorname{span}(\boldsymbol\phi_{-z})$ and both functions are strictly positive, there exists $C>0$ such that
		$\boldsymbol{\psi}^N(x)
		=C\boldsymbol{\phi}_{-z}(x)e^{z\cdot x}$.  Consequently,
		\[
		\lambda^*=\lambda_L^*(-z)
		\geq\min_{\alpha\in\mathbb R^N}\lambda_L^*(\alpha).
		\]
		
		On the other hand, for any $\mu\in\mathcal I$,
		$\boldsymbol{\psi}_{\mu}(x)
		=\boldsymbol{\phi}_{\mu}(x)e^{-\mu\cdot x}$ satisfies
		$\mathcal L\boldsymbol{\psi}_{\mu}
		=\lambda_L^*(\mu)\boldsymbol{\psi}_{\mu}$.
		Definition~\eqref{EQ2.1} therefore gives
		\[
		\lambda^*\leq\min_{\alpha\in\mathbb R^N}
		\lambda_L^*(\alpha).
		\]
		Thus $\lambda^* = \min_{\alpha \in \mathbb{R}^N} \lambda_L^*(\alpha)$ for compactly supported kernels.
		
\smallskip
\noindent Step 2. General kernels under the nonintegrability condition.
		For each $i=1,\dots,m$, let $J_i$ be a general kernel satisfying (J). Define
\[
\varrho(x)=\max\{\min\{2-|x|,1\},0\},\qquad
J_{i,n}(x)=J_i(x)\varrho(x/n).
\]
Write $\lambda_n(l)$ and $\lambda_{L,n}^*(\mu)$ for the bounded-domain
and weighted periodic spectral quantities with $J_i$ replaced by
$J_{i,n}$. The loss term remains $-I$. The proof of Step~1 applies
to these compactly supported kernels without requiring their masses
to equal one: the convolution stays positive and the multiplication
matrix is still $A-I$.
		Restriction of whole-space super-solutions and monotonicity with respect to the kernels give
		\[
		\min_{\mu\in\mathbb R^N}\lambda_L^*(\mu)
\geq\lambda^*\geq\lim_{l\to\infty}\lambda(l)
\geq\lim_{l\to\infty}\lambda_n(l)
=\min_{\mu\in\mathbb R^N}\lambda_{L,n}^*(\mu).
		\]
		Thus it suffices to prove that 
		\[
		\lim_{n\to\infty}\min_{\mu\in\mathbb{R}^N}\lambda_{L,n}^*(\mu) = \min_{\mu\in\mathbb{R}^N}\lambda_{L}^{*}(\mu).
		\]
		Let $\mu_n$ be a minimizer of $\lambda_{L,n}^*(\cdot)$, i.e.,
		$\lambda_{L,n}^*(\mu_n) = \min_{\mu\in\mathbb{R}^N}\lambda_{L,n}^*(\mu)$.
		
\noindent Claim 1. The minimizing sequence $(\mu_n)$ is bounded.
		Choose $r_0,j_0>0$ such that
		$J_i(y)\geq j_0$ for every $i$ and $|y|\leq r_0$.  For all large
		$n$, the same bound holds for $J_{i,n}$. For a positive periodic super-solution, evaluate the $i$th inequality at a minimum point of its $i$th component. The convolution is at least its mass times this minimum, and the off-diagonal terms are nonnegative. Taking the infimum over super-solutions gives
		\[
		\lambda_{L,n}^*(\mu)
		\geq\max_{i\in\Lambda}
		\left\{\int_{\mathbb R^N}J_{i,n}(y)e^{\mu\cdot y}\,\mathrm dy
		-1+\min_xa_{ii}(x)\right\}.
		\]
		By symmetry of the ball,
		\[
		\int_{\mathbb R^N}J_{i,n}(y)e^{\mu\cdot y}\,\mathrm dy
		\geq j_0\int_{|y|\leq r_0}\cosh(\mu\cdot y)\,\mathrm dy
		\geq C_0+C_1|\mu|^2.
		\]
		Thus the functions $\lambda_{L,n}^*$ are equicoercive: for each
		fixed upper bound on their values, the corresponding sublevel sets
		lie in one bounded set independent of $n$. Since their
		minimum values are bounded above by $\lambda_{L,n}^*(0)\leq\lambda_L^*(0)$, the
		minimizers $\mu_n$ remain in a fixed compact set.  Thus Claim~1 holds.
		
We use a monotone-minima argument, which also treats
weights on the boundary of the moment domain. Set the finite continuous functions
$f_n(\mu)=\lambda_{L,n}^*(\mu)$ on $\mathbb R^N$,
and $f(\mu)=\lambda_L^*(\mu)$ on $\mathcal I$,
with $f=+\infty$ outside $\mathcal I$.
Kernel ordering implies $f_n\le f_{n+1}\le f$. The fixed-weight comparison from
\eqref{eq:weighted-CW-sandwich}, applied to a strictly positive test
function with positive minimum, shows $f_n(\mu)\uparrow f(\mu)$
for $\mu\in\mathcal I$. If $\mu\notin\mathcal I$, at least one
moment diverges. The componentwise minimum test used in Claim~1,
followed by monotone convergence, then shows
$f_n(\mu)\uparrow+\infty=f(\mu)$.
Choose minimizers $\mu_n$. By Claim~1, a subsequence converges to
some $\mu_0$. Fix an arbitrary integer $k$. For $n\ge k$,
\[
 \min_\mu f_n(\mu)=f_n(\mu_n)\ge f_k(\mu_n).
\]
The compact-kernel spectral function $f_k$ is continuous, so
\[
 \liminf_{n\to\infty}\min_\mu f_n(\mu)
 \ge f_k(\mu_0).
\]
As $k$ is arbitrary, monotone convergence gives the lower bound
$f(\mu_0)$; in particular $\mu_0\in\mathcal I$. Conversely, for
every fixed $\mu\in\mathcal I$,
\[
 \min_\nu f_n(\nu)\le f_n(\mu)\le f(\mu).
\]
Taking the infimum over $\mu$ and combining both inequalities proves
\[
 \lim_{n\to\infty}\min_{\mu\in\mathbb R^N}f_n(\mu)
 =\min_{\mu\in\mathbb R^N}f(\mu).
\]
The same reasoning applies to every convergent subsequence of
minimizers; hence no uniform moment bound near
$\partial\mathcal I$ or convergence of minimizers in the interior
is being assumed. This completes Step~2.

\smallskip
\noindent Step 3. General matrices $A$.
By Lemma~\ref{lem:approx}, for every $\varepsilon>0$ there is a
cooperative irreducible periodic matrix $\underline A^\varepsilon$
such that
\[
A-2\varepsilon I\leq\underline A^\varepsilon\leq A-\varepsilon I,
\qquad\underline A^\varepsilon\longrightarrow A
\quad\text{uniformly as }\varepsilon\to0^+.
\]
These matrices satisfy the nonintegrability condition and do not
depend on $\mu$. The corresponding operators admit principal
eigenvalues $\lambda_\varepsilon^-(\mu)$ for every $\mu\in\mathcal I$,
and Steps~1--2 apply. Let $\underline\lambda^\varepsilon(l)$ denote
the corresponding bounded-domain quantity.

For $\mu\in\mathcal I$ and $\delta>0$, choose
$\boldsymbol\phi\in X^{++}$ such that
$\mathcal L_\mu[\boldsymbol\phi]\leq
(\lambda_L^*(\mu)+\delta)\boldsymbol\phi$.
Then $\boldsymbol\psi(x)=e^{-\mu\cdot x}\boldsymbol\phi(x)$ belongs to
$\mathscr A_{\mathcal J}^{+}$, since
$J_i*\psi_i=e^{-\mu\cdot x}(\mathcal J_\mu[\boldsymbol\phi])_i$
is finite and continuous. The same conjugation gives
$\mathcal L[\boldsymbol\psi]\leq
(\lambda_L^*(\mu)+\delta)\boldsymbol\psi$.
Letting $\delta\downarrow0$ in \eqref{EQ2.1} yields
$\lambda^*\leq\lambda_L^*(\mu)$ without requiring a principal eigenfunction.
Consequently, restriction of super-solutions and monotonicity in the coefficients give
		\[
		\min_\mu\lambda_L^*(\mu)
\geq\lambda^*\geq\lim_{l\to\infty}\lambda(l)
\geq\lim_{l\to\infty}\underline\lambda^\varepsilon(l)
=\min_\mu\lambda_\varepsilon^-(\mu).
		\]
Thus it suffices to prove that
\[
\lim_{\varepsilon\to0^+}\min_\mu\lambda_\varepsilon^-(\mu)
=\min_\mu\lambda_L^*(\mu).
\]
Let $\mu_\varepsilon$ be a minimizer of $\lambda_\varepsilon^-$.
As in Step~2, $(\mu_\varepsilon)$ remains bounded: by
Lemma~\ref{lem:approx}(iii),
\[
\lambda_L^*(\mu_\varepsilon)-2\varepsilon
\leq\lambda_\varepsilon^-(\mu_\varepsilon)
\leq\lambda_\varepsilon^-(0)\leq\lambda_L^*(0),
\]
and the proof of part~(i) gives coercivity. Along any sequence
$\varepsilon\to0^+$, we may therefore extract a subsequence such
that $\mu_\varepsilon\to\mu_0$. Lower semicontinuity from the proof of part~(i)
and the same uniform estimate imply
\[
\lambda_L^*(\mu_0)
\leq\liminf_{\varepsilon\to0^+}\lambda_\varepsilon^-(\mu_\varepsilon)
\leq\limsup_{\varepsilon\to0^+}\lambda_\varepsilon^-(\mu_\varepsilon)
\leq\min_\mu\lambda_L^*(\mu).
\]
In particular $\mu_0\in\mathcal I$, all these inequalities are
equalities, and $\mu_0$ minimizes $\lambda_L^*$. Since the initial
sequence was arbitrary, the desired limit follows. This proves part~(ii).
	\end{proof}
	
\begin{proof}[\textbf{Proof of Theorem~\ref{TH2.1}(iii)}]
The inequality $\lambda'\geq\lambda_L^*(0)$ follows from
Lemma~\ref{L11}(iii), since its periodic sub-solutions are admissible
in \eqref{eq:whole-space-lambda-prime}. Hence we only need to prove
that $\lambda'\leq\lambda_L^*(0)$.
Suppose, to the contrary, that there exist $\lambda>\lambda_L^*(0)$
and a bounded continuous function $\boldsymbol\phi\gg\boldsymbol0$
such that $\mathcal L\boldsymbol\phi\geq\lambda\boldsymbol\phi$.
By the definition of $\lambda_L^*(0)$, we can find
$\lambda_L^*(0)<\lambda_0<\lambda$ and
$\boldsymbol\varphi_0\in X^{++}$ satisfying
\[
\mathcal L\boldsymbol\varphi_0\leq\lambda_0\boldsymbol\varphi_0.
\]
Define
\[
\gamma=\max_{1\leq i\leq m}\sup_{x\in\mathbb R^N}
\frac{\phi_i(x)}{\varphi_{0,i}(x)},\qquad
\boldsymbol w=\gamma\boldsymbol\varphi_0-\boldsymbol\phi.
\]
Since $\boldsymbol\varphi_0$ is positive and periodic and
$\boldsymbol\phi$ is bounded, $0<\gamma<\infty$.
Then $\boldsymbol w\geq\boldsymbol0$ and
$\inf_{i,x}w_i(x)=0$, although this infimum need not be attained.
Moreover,
\[
\begin{aligned}
(\mathcal L-\lambda)\boldsymbol w
&=\gamma(\mathcal L-\lambda)\boldsymbol\varphi_0
  -(\mathcal L-\lambda)\boldsymbol\phi\\
&\leq\gamma(\lambda_0-\lambda)\boldsymbol\varphi_0
\leq-\delta\boldsymbol1,
\end{aligned}
\]
where
\[
\gamma_0=\min_{i,x}\varphi_{0,i}(x)>0,
\qquad\delta=\gamma(\lambda-\lambda_0)\gamma_0>0.
\]

Consider $\boldsymbol r(x)=\frac{|x|^2}{1+|x|^2}\boldsymbol1$.
Then $\boldsymbol r(0)=\boldsymbol0$, $0\leq\boldsymbol r\leq\boldsymbol1$,
and $\boldsymbol r(x)\to\boldsymbol1$ as $|x|\to\infty$.
Put
\[
 B=2+|\lambda|+\max_{i,x}\sum_{j=1}^m|a_{ij}(x)|.
\]
Since $\int J_i=1$, for every $y\in\mathbb R^N$ we have
\[
\|(\mathcal L-\lambda)[T_y\boldsymbol r]\|_\infty\leq B,
\qquad
\|(\mathcal L-\lambda)\boldsymbol1\|_\infty\leq B,
\]
where $T_y\boldsymbol r(x)=\boldsymbol r(x-y)$.
Choose $b>2B/\delta$. Since $\inf_{i,x}w_i(x)=0$, there exist
$i_0$ and $x_0$ such that
\[
 w_{i_0}(x_0)<\min\left\{\frac1{2b},\frac{\delta}{4B}\right\}.
\]
Define
\[
\widetilde{\boldsymbol w}
=\boldsymbol w+\frac1b T_{x_0}\boldsymbol r,
\qquad
\alpha=\min_{1\leq i\leq m}\inf_{x\in\mathbb R^N}\widetilde w_i(x).
\]
For $|x-x_0|$ sufficiently large, every component of
$b^{-1}T_{x_0}\boldsymbol r(x)$ is greater than $w_{i_0}(x_0)$.
Thus the infimum defining $\alpha$ is attained at some component
$i_1$ and point $x_1$, and
\[
0\leq\alpha\leq\widetilde w_{i_0}(x_0)
=w_{i_0}(x_0)<\frac{\delta}{4B}.
\]
Now $\widetilde{\boldsymbol w}-\alpha\boldsymbol1\geq\boldsymbol0$,
and the preceding estimates give
\[
\begin{aligned}
(\mathcal L-\lambda)[\widetilde{\boldsymbol w}-\alpha\boldsymbol1]
&=(\mathcal L-\lambda)\boldsymbol w
  +\frac1b(\mathcal L-\lambda)[T_{x_0}\boldsymbol r]
  -\alpha(\mathcal L-\lambda)\boldsymbol1\\
&\leq\left(-\delta+\frac{B}{b}+\alpha B\right)\boldsymbol1
<-\frac{\delta}{4}\boldsymbol1.
\end{aligned}
\]
On the other hand, its $i_1$th component vanishes at $x_1$.
At this point the convolution and off-diagonal terms are nonnegative,
while the diagonal term vanishes. Hence
\[
\bigl((\mathcal L-\lambda)
[\widetilde{\boldsymbol w}-\alpha\boldsymbol1]\bigr)_{i_1}(x_1)\geq0,
\]
a contradiction. Therefore $\lambda'\leq\lambda_L^*(0)$, and the
proof is complete.
\end{proof}

\subsubsection{Detailed balance}

We use the diagonal coefficient approximation to preserve detailed
balance and compare the eigenvalues at opposite exponential weights.

\begin{proof}[Proof of Proposition~\ref{prop:detailed-balance-rigidity}]
Since the kernels are even, their common moment domain is
symmetric. Fix $\mu\in\mathcal I$ and use the lower coefficient
approximations $A^\varepsilon=\underline A^\varepsilon$ from
Lemma~\ref{lem:approx}, writing
$\mathcal L_\mu^\varepsilon=\underline{\mathcal L}_\mu^\varepsilon$.
The matrix bounds in that lemma imply that
$A^\varepsilon-A$ is diagonal, so
$DA^\varepsilon=(A^\varepsilon)^{\mathsf T}D$.
Let $\boldsymbol\phi_\mu^\varepsilon,
\boldsymbol\phi_{-\mu}^\varepsilon\in X^{++}$ be the corresponding
periodic eigenfunctions. Evenness of the kernels and Fubini's theorem give
\[
\int_{\mathbb T^N}(\boldsymbol\phi_{-\mu}^\varepsilon)^{\mathsf T}
D\mathcal L_\mu^\varepsilon\boldsymbol\phi_\mu^\varepsilon
=\int_{\mathbb T^N}
(\mathcal L_{-\mu}^\varepsilon\boldsymbol\phi_{-\mu}^\varepsilon)^{\mathsf T}
D\boldsymbol\phi_\mu^\varepsilon.
\]
All integrals are finite because the eigenfunctions are bounded and
the weighted kernels are integrable. Substitution of the two
eigenvalue equations yields
\[
\bigl(\lambda_\varepsilon^-(\mu)-\lambda_\varepsilon^-(-\mu)\bigr)
\int_{\mathbb T^N}(\boldsymbol\phi_{-\mu}^\varepsilon)^{\mathsf T}
D\boldsymbol\phi_\mu^\varepsilon=0.
\]
The integral is strictly positive. Thus
$\lambda_\varepsilon^-(\mu)=\lambda_\varepsilon^-(-\mu)$, and
Lemma~\ref{lem:approx}(iii) gives
$\lambda_L^*(\mu)=\lambda_L^*(-\mu)$ on letting
$\varepsilon\downarrow0$. Convexity from
Theorem~\ref{TH2.1}(i) places the minimum at zero.
Parts~(i)--(ii) of that theorem then prove
\eqref{eq:detailed-balance-exact-threshold}.
\end{proof}

\section{Large-period spectral asymptotics}
\label{sec:large-period-asymptotics}
\label{sec:periodic-eigenvalue-limits}

We first characterize the large-period spectral limits through a
periodic Hamilton--Jacobi problem in arbitrary spatial dimension.
We then derive quantitative corrections and concentration estimates in
one dimension, for scalar equations and reversible cooperative systems.

\subsection{Main results}

Fix a continuous $1$-periodic matrix profile $A$ and set
$A_L(x)=A(x/L)$ with common spatial period $L=L_1=\cdots=L_N$.
After rescaling $x=L\tilde x$ and dropping the tildes, throughout this
section and Section~\ref{sec:quantitative-period} we use $\mathcal C=(0,1)^N$,
$\mathbb T^N=\mathbb R^N/\mathbb Z^N$, and
\[
X=\{\boldsymbol\varphi\in C(\mathbb R^N,\mathbb R^m):
\boldsymbol\varphi(x+p)=\boldsymbol\varphi(x)
\ \text{for every }x\in\mathbb R^N,\ p\in\mathbb Z^N\}.
\]
Elsewhere these spaces retain their original $\mathbb L$-periodic meanings.
The exponential weight $\mu$ is measured in the original spatial variables,
so the rescaled weighted operator is
\begin{equation}\label{eq:rescaled-weighted-operator}
(\mathcal L_{L,\mu}\boldsymbol\varphi)_i(x)
=\int_{\mathbb R^N}J_i(z)e^{\mu\cdot z}
\varphi_i(x-z/L)\,\mathrm dz-\varphi_i(x)
+\sum_{j=1}^m a_{ij}(x)\varphi_j(x).
\end{equation}
Its periodic generalized principal eigenvalue is denoted by
$\lambda_L^*(\mu)$. We write $\lambda^*(L)$ and $\lambda'(L)$ for the
corresponding whole-space quantities when the period varies.

\subsubsection{Hamilton--Jacobi limits}

	The large-period limit is described by a periodic Hamiltonian.
For the moments in \eqref{eq:kernel-exponential-moments}, set
	\[
	D(p):=\operatorname{diag}\bigl(M_1(p)-1,\ldots,M_m(p)-1\bigr).
	\]
	The local Hamiltonian is
	\begin{equation}\label{eq:local-Hamiltonian}
	H(p,x):=s_{\rm PF}\bigl(D(p)+A(x)\bigr),
	\qquad (p,x)\in\mathcal I\times\mathbb T^N,
	\end{equation}
	and $H(p,x):=+\infty$ for $p\notin\mathcal I$. For $\mu\in\mathcal I$,
define the effective
	periodic Hamiltonian
	\begin{equation}\label{eq:effective-periodic-Hamiltonian}
	\overline\lambda_\infty(\mu)
	:=\inf_{W\in C_{\rm per}^1(\mathbb T^N)}
	\sup_{x\in\mathbb T^N}H\bigl(\mu+\nabla W(x),x\bigr),
	\end{equation}
Here $C_{\rm per}^1(\mathbb T^N)$ is the space of scalar-valued,
continuously differentiable functions that are $1$-periodic in every
spatial direction.
Definition \eqref{eq:effective-periodic-Hamiltonian} does not require
a periodic viscosity solution of the corresponding equality.

Related Hamilton--Jacobi problems for period
limits of local diffusion models are studied in \cite{KangLiu2026}.
For $\mu\notin\mathcal I$, set
$\overline\lambda_\infty(\mu)=+\infty$.

\begin{theorem}\label{thm:large-period-limits}
\label{thm:period-size-asymptotics}

Assume that (J) and (A) hold, with $L=L_1=\cdots=L_N$.
Then the following statements hold as $L\to\infty$:
\begin{enumerate}
\item[(i)] For every $\mu\in\mathcal I$,
\[
\lambda_L^*(\mu)\longrightarrow\overline\lambda_\infty(\mu).
\]
The convergence is locally uniform on compact subsets of
$\operatorname{int}\mathcal I$; no uniformity at the boundary
of the moment domain is asserted.
\item[(ii)] The whole-space spectral values satisfy
\[
\lambda'(L)=\lambda_L^*(0)
\longrightarrow\overline\lambda_\infty(0),
\qquad
\lambda^*(L)\longrightarrow
\min_{\mu\in\mathbb R^N}\overline\lambda_\infty(\mu).
\]
\end{enumerate}
\end{theorem}
For nonsymmetric kernels, $\overline\lambda_\infty$ need not be even
or minimized at zero. The two large-period limits agree exactly when
$\overline\lambda_\infty(0)=\min_\mu\overline\lambda_\infty(\mu)$.

\subsubsection{Scalar equations}

\label{subsec:universal-large-period}

For scalar equations, the leading correction is determined by the
kernel's second moment and the curvatures at the global maxima.

\begin{theorem}\label{thm:universal-kernel-large-period}
Let $N=m=1$. Assume that (J) holds and that $J$ is even, with
\[
 0<D_J:=\frac12\int_{\mathbb R}z^2J(z)\,\mathrm dz<\infty.
\]
Suppose that $a\in C^2(\mathbb T)$ has finitely many global maxima
$x_1,\ldots,x_r$, all nondegenerate. Set
\[
 a_*=\max_{\mathbb T}a,
 \qquad \kappa_j=-a''(x_j)>0,
 \qquad \Omega_J=
 \min_{1\le j\le r}\sqrt{\frac{D_J\kappa_j}{2}}.
\]
Then the following statements hold:
\begin{enumerate}
\item[(i)] The periodic principal eigenvalue exists for all sufficiently
large $L$ and
\begin{equation}\label{eq:universal-large-period-asymptotic}
 \lambda^*(L)=\lambda'(L)=\lambda_L^*(0)
 =a_*-\frac{\Omega_J}{L}+o(L^{-1})
 \qquad(L\to\infty).
\end{equation}
\item[(ii)] If $u_L$ is a positive periodic eigenfunction normalized by
$\int_{\mathbb T}|u_L|^2\,\mathrm dx=1$,
then, for every open neighborhood $V$ of the set
$\{x_1,\ldots,x_r\}$ whose complement is nonempty,
\[
 \int_{\mathbb T\setminus V}|u_L(x)|^2\,\mathrm dx
 \leq\frac{\Omega_J+o(1)}{L\,\delta_V},
 \qquad \delta_V:=\inf_{\mathbb T\setminus V}(a_*-a)>0.
\]
In particular, the eigenfunction mass outside $V$ is $O(L^{-1})$.
\end{enumerate}
No fourth moment or differentiability of $J$ is required.
\end{theorem}

For the Laplace kernel $J(z)=\tfrac12e^{-|z|}$, $D_J=1$ and the
asymptotic formula becomes the following.

\begin{corollary}\label{thm:large-period-harmonic}
Let $N=m=1$, $J(z)=\tfrac12e^{-|z|}$, and
$A(x)=a(x)\in C^2(\mathbb T)$. Suppose that $a$ has finitely many
global maxima $x_1,\ldots,x_r$ modulo periods and that they are all
nondegenerate: $\kappa_j=-a''(x_j)>0$. Set
$a_*=\max a$ and $\omega_* =\min_{1\le j\le r}\sqrt{\kappa_j/2}$.
Then, as $L\to\infty$,
\[
\displaystyle
\lambda^*(L)=\lambda'(L)=\lambda_L^*(0)
=a_*-\frac{\omega_*}{L}+o(L^{-1}).
\]
The flattest global maxima determine this correction; for a unique
maximum, $\omega_*=\sqrt{\kappa_1/2}$.
In particular, for $a(x)=a_0+\delta\cos(2\pi x)$ with $\delta>0$,
\[
\lambda'(L)=a_0+\delta
-\frac{\pi\sqrt{2\delta}}{L}+o(L^{-1}).
\]
The equality of the two whole-space thresholds uses scalarity and
the symmetry of the kernel. It need not hold for asymmetric kernels.
\end{corollary}

The following example compares two global maxima of equal height.

\begin{example}\label{ex:two-well-large-period}
Take $J(z)=\tfrac12e^{-|z|}$, fix $0<\beta<1$, and consider
\[
a(x)=a_*-\sin^2(2\pi x)\bigl(1+\beta\cos(2\pi x)\bigr).
\]
The only global maxima are $0$ and $1/2$ modulo integers,
with respective curvatures
$\kappa_1=8\pi^2(1+\beta)$ and
$\kappa_2=8\pi^2(1-\beta)$. Consequently,
Corollary~\ref{thm:large-period-harmonic} gives
\[
\lambda^*(L)=\lambda'(L)
=a_*-\frac{2\pi\sqrt{1-\beta}}{L}+o(L^{-1}).
\]
The broader maximum at $x=1/2$ determines the correction.
\end{example}

\subsubsection{Reversible systems}
\label{subsec:reversible-system-reduction}

For the reversible systems considered below, let
$D_J=\frac12\int_{\mathbb R}z^2J(z)\,\mathrm dz$ for the common kernel $J$.
Given a fixed positive diagonal symmetrizing matrix $D$, set
\[
S(x)=D^{1/2}A(x)D^{-1/2},\qquad
s(x)=s_{\rm PF}(A(x)),
\]
and let $v(x)>0$ be the Euclidean-unit Perron eigenvector of $S(x)$.
In the following theorem and its proof, $s$ differs by one from
$s_{\rm PF}(A-I)$ in \eqref{s(x)}. Define
\[
\Pi(x)=v(x)v(x)^{\mathsf T},\qquad
\gamma=\min_{x\in\mathbb T}[s(x)-\lambda_2(S(x))]>0,
\]
where $\lambda_2(S(x))$ is the second-largest eigenvalue.
Let $\Lambda_L$ be the spectral bound of the periodic system operator
$\mathcal L_{L,0}$, and let $\Lambda_L^{\rm sc}(s)$ be the spectral
bound on $L^2(\mathbb T)$ of $K_L-I+s(x)$, where
\[
(K_Lf)(x)=\int_{\mathbb R}J(z)f(x-z/L)\,\mathrm dz.
\]
We write $\operatorname{Lip}(v)$ for the Lipschitz constant relative
to the torus distance $|x-y|_{\mathbb T}=\min_{k\in\mathbb Z}|x-y+k|$.

\begin{theorem}\label{thm:reversible-system-reduction}
Let $N=1$ and $m\geq2$. Assume that (J) and (A) hold.
Suppose that $J_1=\cdots=J_m=J$ is even with $0<D_J<\infty$,
and that $A\in C^1(\mathbb T;\mathbb R^{m\times m})$ satisfies
$DA(x)=A(x)^{\mathsf T}D$ for a fixed positive diagonal matrix $D$.
Then, with the notation above, the following statements hold:
\begin{enumerate}
\item[(i)] The vector field $v$ is periodic and Lipschitz, and for every $L>0$,
\begin{equation}\label{eq:reversible-finite-period-comparison}
0\leq\Lambda_L^{\rm sc}(s)-\Lambda_L
\leq\frac{D_J\operatorname{Lip}(v)^2}{L^2}.
\end{equation}

\item[(ii)] Suppose in addition that $A\in C^2(\mathbb T;\mathbb R^{m\times m})$ and that the global
maxima of $s$ are precisely $x_1,\ldots,x_r$, all nondegenerate.
Set $s_*=\max_{\mathbb T}s$, $\kappa_j=-s''(x_j)>0$, and
$\Omega_J=\min_j\sqrt{D_J\kappa_j/2}$. Then
\begin{equation}\label{eq:reversible-system-localization}
\Lambda_L=s_*-\frac{\Omega_J}{L}+o(L^{-1})\qquad(L\to\infty).
\end{equation}
For all sufficiently large $L$, there is a positive periodic principal
eigenfunction $u_L$. Scale $w_L=D^{1/2}u_L$ to unit norm in
$L^2(\mathbb T;\mathbb R^m)$.
Given an open neighborhood $V$ of all the maxima with nonempty complement, put
$\delta_V=\inf_{\mathbb T\setminus V}(s_*-s)>0$. Then
\begin{align}
\int_{\mathbb T\setminus V}|w_L(x)|^2\,\mathrm dx
&\leq\frac{\Omega_J+o(1)}{L\delta_V},
\label{eq:reversible-concentration}\\
\int_{\mathbb T}|(I-\Pi(x))w_L(x)|^2\,\mathrm dx
&\leq\frac{\Omega_J+o(1)}{\gamma L}.
\label{eq:reversible-adiabatic-alignment}
\end{align}
\end{enumerate}
\end{theorem}

The next example gives an explicit system localization coefficient.

\begin{example}\label{ex:reversible-system-curvature}
For $m=2$ let
\[
 A(x)=\begin{pmatrix}
 a_0+\delta\cos(2\pi x)&b\\ b&d_0
 \end{pmatrix},\qquad b>0,\quad\delta>0.
\]
The normalized Perron eigenvector depends on $x$ whenever
$\delta\ne0$. Nevertheless, $s$ has a unique nondegenerate maximum
at $x_*=0$ modulo periodicity. With
$h=a_0+\delta-d_0$, we obtain
\[
 \vartheta=\frac12\left(1+\frac{h}{\sqrt{h^2+4b^2}}\right),
 \quad 0<\vartheta<1,\qquad
 -s''(0)=4\pi^2\delta\vartheta.
\]
Consequently, for every even kernel satisfying the hypotheses of
Theorem~\ref{thm:reversible-system-reduction},
\[
 \Lambda_L=s_*-\frac{\pi}{L}
 \sqrt{2D_J\delta\vartheta}+o(L^{-1}),
 \qquad
 s_*=\frac{a_0+\delta+d_0+\sqrt{h^2+4b^2}}2.
\]
Thus even when the local Perron eigendirection varies with $x$,
its variation changes the system spectral bound by at most $O(L^{-2})$
relative to the scalar comparison problem, and not the first
localization coefficient.
\end{example}

\subsection{Proofs}
\label{subsec:localization-proofs}

We prove the Hamilton--Jacobi limit first, then the scalar and
reversible-system refinements.

\subsubsection{Hamilton--Jacobi limits}

The supremum in \eqref{eq:effective-periodic-Hamiltonian} allows for
possible discontinuities of the exponential moments at the boundary
of $\mathcal I$.
In this subsection, we use the unit-period variables introduced above.
The rescaled convolution operator is
\[
(\mathcal J_L*\boldsymbol\varphi)_i(x)
=L^N\int_{\mathbb R^N}J_i(L(x-y))\varphi_i(y)\,\mathrm dy,
\qquad i\in\Lambda,
\]
and $\mathcal L_L\boldsymbol\varphi
=\mathcal J_L*\boldsymbol\varphi-\boldsymbol\varphi+A(x)\boldsymbol\varphi$.
We first consider compactly supported kernels and zero exponential weight.
In this case the Hamiltonian $H$ in \eqref{eq:local-Hamiltonian} is finite
and continuous. Fix $x$, let
$\boldsymbol v^1,\boldsymbol v^2\gg0$ be Perron vectors at $p_1,p_2$,
and set $v_i=(v_i^1)^\theta(v_i^2)^{1-\theta}$, where $0<\theta<1$.
Convexity of $M_i$ and the weighted arithmetic--geometric mean inequality
give, componentwise,
\[
\frac{[(D(\theta p_1+(1-\theta)p_2)+A(x))\boldsymbol v]_i}{v_i}
\leq\theta H(p_1,x)+(1-\theta)H(p_2,x).
\]
The matrix bounds \eqref{eq:matrix-Collatz-Wielandt} prove convexity of $H$.
Positivity of every $J_i$ near the origin implies uniform coercivity
in $p$, meaning that $H(p,x)\to+\infty$ as $|p|\to\infty$, uniformly
in $x$. Consider the periodic ergodic problem
\begin{equation*}
H(\nabla W(x),x)=\overline\lambda_\infty(0)
\quad\text{on }\mathbb T^N
\end{equation*}
in the viscosity sense: a sub-solution is upper semicontinuous and
satisfies $H(\nabla\phi(x),x)\leq\overline\lambda_\infty(0)$ whenever
a smooth test function $\phi$ touches it from above at $x$; a
super-solution is lower semicontinuous and satisfies the reverse
inequality for tests touching from below. A continuous function
satisfying both conditions is a viscosity solution. For the present
compactly supported kernels, the preceding continuity, convexity and
coercivity give a periodic viscosity solution by classical periodic
cell-problem theory \cite{LionsPapanicolaouVaradhan1988},
whose unique additive constant is characterized by
\[
\overline\lambda_\infty(0)
=\inf_{W\in C_{\rm per}^1(\mathbb T^N)}
\max_{x\in\mathbb T^N}H(\nabla W(x),x).
\]

We prove the compact-support, zero-weight case of
Theorem~\ref{thm:large-period-limits}(i) first, then extend the argument
to arbitrary admissible weights and general kernels.

We first establish a spectral gap that is uniform for large $L$.
Write $B(x)=A(x)-I$, choose $x_0$ with
$s(x_0)=s^*$, and let $v\gg0$ satisfy
$B(x_0)v=s^*v$. Choose $r,j_0>0$ such that
$J_i(z)\geq j_0$ for all $|z|\leq r$ and all $i$. Put
\[
\Omega_L=B_{r/(2L)}(x_0),\qquad
\kappa_0=j_0|B_{r/2}(0)|>0,\qquad
\eta_L=\max_i\sup_{x\in\overline\Omega_L}
\frac{|((B(x)-B(x_0))v)_i|}{v_i}.
\]
Here $B_\rho(y)$ is the open Euclidean ball of radius $\rho$ centered
at $y$, and $|B_\rho(y)|$ is its Lebesgue measure.
Continuity gives $\eta_L\to0$. If $\boldsymbol\psi\gg0$ is a periodic
upper test function satisfying $\mathcal L_L\boldsymbol\psi
\leq\lambda\boldsymbol\psi$, choose $(i,x)$ attaining
$\min_{j,y\in\overline\Omega_L}\psi_j(y)/v_j=:a>0$.
For $y\in\Omega_L$ we have $|L(x-y)|\leq r$, and cooperativity gives
\[
\lambda a v_i
\geq (\mathcal L_L\boldsymbol\psi)_i(x)
\geq a v_i j_0L^N|\Omega_L|+a(B(x)v)_i
\geq a v_i(\kappa_0+s^*-\eta_L).
\]
Taking the infimum over upper test functions and using the constant
vector $\boldsymbol1=(1,\ldots,1)^{\mathsf T}$ for the upper bound, we obtain, for all sufficiently
large $L$,
\begin{equation}\label{eq:large-period-uniform-spectral-gap}
s^*+\frac{\kappa_0}{2}
\leq\lambda_L^*(0)
\leq C_A.
\end{equation}
Lemma~\ref{L11}(i) supplies the principal eigenfunction used in the
following estimate.

\begin{lemma}\label{L3.2}
Assume that (J) and (A) hold and that every $J_i$ has compact support.
For all sufficiently large $L$, let
$\boldsymbol\varphi_L=(\varphi_1,\ldots,\varphi_m)>0$ be a periodic
principal eigenfunction of $\mathcal L_L$, normalized by
$\max_i\|\varphi_i\|_\infty=1$. Set
\[
W_{i,L}(x):=-\frac{1}{L}\log \varphi_i(x),\quad i=1,\ldots,m;
\qquad\text{then}\quad
\min_{1\le i\le m}\inf_{\mathbb{R}^N}W_{i,L}(x)=0.
\]
	Then there exist $L_0\geq1$ and $C>0$, independent of $L$, such that
	$0\leq W_{i,L}(x)\leq C$ for every $L\geq L_0$,
	$i\in\Lambda$, and $x\in\mathbb R^N$.
\end{lemma}

\begin{proof}[Proof of Lemma~\ref{L3.2}]
Increase $L_0$ so that \eqref{eq:large-period-uniform-spectral-gap} holds
for $L\geq L_0$. Since $B(x)$ is cooperative and irreducible, the
resolvent matrix $R(x,\lambda)=(\lambda I-B(x))^{-1}$
has strictly positive entries whenever $\lambda>s(x)$. Its entries are
uniformly bounded above and below by positive constants for
$x\in\mathbb T^N$ and
$s^*+\kappa_0/2\leq\lambda\leq C_A$. The eigenvalue equation can be written as
\[
\boldsymbol\varphi_L(x)
=R(x,\lambda_L^*(0))(\mathcal J_L*\boldsymbol\varphi_L)(x).
\]
Choose $R_0>r$ containing the supports of all the kernels, and put
$F_L=\sum_j\varphi_j$. The positivity of the kernels on $B_r$ and their
boundedness and compact support give constants $c_0,C_0>0$, independent
of $L\geq L_0$, such that
\begin{equation}\label{eq:local-integral-Harnack-bounds}
c_0L^N\int_{B_{r/L}(x)}F_L(y)\,\mathrm dy
\leq\varphi_i(x)
\leq C_0L^N\int_{B_{R_0/L}(x)}F_L(y)\,\mathrm dy.
\end{equation}
Let $h=\boldsymbol1_{B_r(0)}$, the indicator function of $B_r(0)$, and
$h_L(z)=L^Nh(Lz)$. The lower bound gives both
$\varphi_i\geq c_0h_L*F_L$ and $F_L\geq c_0h_L*F_L$.
Choose a fixed integer $k\geq2$ with $kr>2R_0$. The $k$-fold convolution
$h^{*k}=h*\cdots*h$ is continuous and strictly positive in $B_{kr}(0)$, so its minimum on
$\overline B_{2R_0}(0)$ is positive. Iterating the lower bound $k$ times
therefore gives a constant $c_k>0$, independent of $L$, such that
\[
\varphi_i(x)\geq c_0^k(h_L^{*k}*F_L)(x)
\geq c_kL^N\int_{B_{2R_0/L}(x)}F_L(y)\,\mathrm dy.
\]
If $|x-y|\leq R_0/L$, then
$B_{R_0/L}(y)\subset B_{2R_0/L}(x)$. Combining this fact with the upper
bound in \eqref{eq:local-integral-Harnack-bounds} yields
$\varphi_j(y)\leq(C_0/c_k)\varphi_i(x)$ for all $i,j$.
We thus obtain the pointwise Harnack estimate: for a fixed $C_H\geq1$,
\begin{equation}\label{eq:scaled-Harnack-chain}
 \max_{1\leq i\leq m}\sup_{B_{R_0/(2L)}(x)}\varphi_i
 \leq C_H
 \min_{1\leq i\leq m}\inf_{B_{R_0/(2L)}(x)}\varphi_i,
 \qquad x\in\mathbb T^N.
\end{equation}
Any two points of the unit torus can be joined by at most \(C_1L+1\)
overlapping balls of radius \(R_0/(2L)\).  Iterating
\eqref{eq:scaled-Harnack-chain} along such a chain and using the coupling
between the components yields
\[
 \frac{\max_{i,x}\varphi_i(x)}
      {\min_{i,x}\varphi_i(x)}
 \leq C_H^{\,C_1L+1}.
\]
The normalization \(\max_{i,x}\varphi_i(x)=1\) therefore gives
\[
 0\leq W_{i,L}(x)
 =-\frac1L\log\varphi_i(x)
 \leq C_1\log C_H+\frac{\log C_H}{L},
 \qquad i\in\Lambda,\quad x\in\mathbb T^N.
\]
This completes the proof.
\end{proof}
\begin{proof}[\textbf{Proof of Theorem~\ref{thm:large-period-limits}(i): compact support and zero weight}]
We work on the compact torus $\mathbb T^N$. Write
$\lambda_L=\lambda_L^*(0)$ and
$\varepsilon_L=L^{-1}$. The positive periodic eigenfunction furnished by Lemma~\ref{L11}(i) and
\eqref{eq:large-period-uniform-spectral-gap} satisfies, for each
$i\in\Lambda$ and $x\in\mathbb T^N$,
\begin{equation}\label{eq:HJ-logarithmic-identity}
\lambda_L=\int_{\mathbb R^N}J_i(z)
 e^{L[W_{i,L}(x)-W_{i,L}(x-z/L)]}\,\mathrm dz-1
 +\sum_{j=1}^m a_{ij}(x)e^{L[W_{i,L}(x)-W_{j,L}(x)]}.
\end{equation}
By Lemma~\ref{L3.2}, $0\leq W_{i,L}\leq C$ uniformly in $L$ and
$i$. By \eqref{eq:large-period-uniform-spectral-gap}, the eigenvalues
$\lambda_L$ are bounded. Fix an arbitrary sequence $L_k\to\infty$ along
which $\lambda_{L_k}\to\lambda_\infty$.

\smallskip
\noindent Step 1. The lower half-relaxed limit.
Following the half-relaxed-limit method of \cite{BarlesPerthame1987}
(see also \cite[Section~6]{CrandallIshiiLions1992}), set
\[
U_*(x)=\min_{1\leq i\leq m}
\liminf_{\substack{k\to\infty\\y\to x}}W_{i,L_k}(y).
\]
This function is bounded, periodic and lower semicontinuous. Let
$\phi\in C^2$ touch $U_*$ from below at $x_0$, with a strict local
minimum of $U_*-\phi$ there. Choose
$\boldsymbol v=(v_1,\ldots,v_m)\gg0$ such that
\[
\bigl[D(\nabla\phi(x_0))+A(x_0)\bigr]\boldsymbol v
=H(\nabla\phi(x_0),x_0)\boldsymbol v.
\]
A contact test is needed only near $x_0$. To control the nonlocal
integral, extend $\phi$ smoothly and periodically outside a small
neighborhood of $x_0$, then subtract a sufficiently large nonnegative
smooth cutoff which vanishes near $x_0$. We denote the resulting test
again by $\phi$. Since $0\leq U_*\leq C$, the cutoff can be chosen so
that $U_*-\phi$ has a strict global minimum at $x_0$. Adding a
constant to $\phi$, we arrange $U_*(x_0)-\phi(x_0)=0$.

For each $k$, minimize over both $i\in\Lambda$ and $x\in\mathbb T^N$
\[
F_{i,k}(x):=W_{i,L_k}(x)+\varepsilon_{L_k}\log v_i-\phi(x).
\]
Let $(i_k,x_k)$ be a minimizing pair. The definition of $U_*$
provides a recovery sequence, and we first pass to the corresponding
subsequence of $k$, without changing the eigenvalue limit. Along it,
\[
\min_{i,x}F_{i,k}(x)\longrightarrow 0,
\qquad x_k\longrightarrow x_0.
\]
Indeed, the recovery sequence bounds the limsup of these minima by
zero. Conversely, compactness and the lower half-relaxed limit bound
every subsequential liminf from below by
$U_*(x)-\phi(x)\geq0$. Strict global minimality then forces
$x=x_0$. Passing to a further subsequence, we may take $i_k=i$
constant. Global minimality, first in the spatial variable and then in
the component index, yields for every $z\in\mathbb R^N$ and every $j$,
\begin{align*}
W_{i,L_k}(x_k)-W_{i,L_k}(x_k-z/L_k)
&\leq\phi(x_k)-\phi(x_k-z/L_k),\\
e^{L_k[W_{i,L_k}(x_k)-W_{j,L_k}(x_k)]}
&\leq\frac{v_j}{v_i}.
\end{align*}
The first inequality is global; in particular it is valid on
the whole support of $J_i$. Since $a_{ij}(x)\geq0$ for $j\ne i$, the
logarithmic identity \eqref{eq:HJ-logarithmic-identity} implies
\[
\lambda_{L_k}\leq
 \int_{\mathbb R^N}J_i(z)
 e^{L_k[\phi(x_k)-\phi(x_k-z/L_k)]}\,\mathrm dz-1
 +\sum_{j=1}^m a_{ij}(x_k)\frac{v_j}{v_i}.
\]
For compactly supported $J_i$, Taylor's formula shows, uniformly on its
support,
\[
L_k[\phi(x_k)-\phi(x_k-z/L_k)]
=z\cdot\nabla\phi(x_0)+o(1).
\]
It follows that
\[
\lambda_\infty\leq
M_i(\nabla\phi(x_0))-1
+\frac{[A(x_0)\boldsymbol v]_i}{v_i}
=H(\nabla\phi(x_0),x_0).
\]
Thus $U_*$ is a periodic viscosity supersolution of
$H(\nabla U,x)=\lambda_\infty$. To justify the direction of this inequality, compare the
lower semicontinuous periodic supersolution $U_*$ at level
$\lambda_\infty$ with a periodic Lipschitz viscosity solution of
$H(\nabla V,x)=\overline\lambda_\infty(0)$. The standard
periodic doubling-of-variables argument
\cite[Section~3]{CrandallIshiiLions1992}, applied at a minimum of
$U_*-V$, uses the supersolution inequality for $U_*$ and the
subsolution inequality for $V$ at the same limiting gradient. The
continuity of $H$ on bounded gradients gives
$\lambda_\infty\leq\overline\lambda_\infty(0)$.
Here boundedness and Lipschitz regularity of $V$ follow from the
coercivity of the continuous Hamiltonian.

\smallskip
\noindent Step 2. The upper half-relaxed limit.
Define
\[
U^*(x)=\max_{1\leq i\leq m}
\limsup_{\substack{k\to\infty\\y\to x}}W_{i,L_k}(y).
\]
Let $\phi\in C^2$ touch $U^*$ from above strictly at $x_0$ and take
a Perron eigenvector $\boldsymbol v\gg0$ at
$(\nabla\phi(x_0),x_0)$. Extend $\phi$ periodically, now adding
a sufficiently large nonnegative cutoff away from $x_0$. Since
$0\leq U^*\leq C$, we may require $U^*-\phi$ to have a strict global
maximum at $x_0$, with value zero. Choose a maximizing pair
$(i_k,x_k)$ for $F_{i,k}$ over $\Lambda\times\mathbb T^N$. Pass to
a subsequence realizing $U^*(x_0)$; compactness and the upper
half-relaxed-limit inequality then imply
$\max_{i,x}F_{i,k}(x)\to0$ and $x_k\to x_0$.
Along a further subsequence, $i_k=i$ is fixed, and global maximality
implies
\begin{align*}
W_{i,L_k}(x_k)-W_{i,L_k}(x_k-z/L_k)
&\geq\phi(x_k)-\phi(x_k-z/L_k)
\quad(z\in\mathbb R^N),\\
e^{L_k[W_{i,L_k}(x_k)-W_{j,L_k}(x_k)]}
&\geq\frac{v_j}{v_i}\quad(1\leq j\leq m).
\end{align*}
Inserting these inequalities into
\eqref{eq:HJ-logarithmic-identity} and using the nonnegativity of the
off-diagonal entries of $A$, we obtain
\[
\lambda_{L_k}\geq
 \int_{\mathbb R^N}J_i(z)
 e^{L_k[\phi(x_k)-\phi(x_k-z/L_k)]}\,\mathrm dz-1
 +\sum_{j=1}^m a_{ij}(x_k)\frac{v_j}{v_i}.
\]
Compact support and the same uniform Taylor estimate now imply
$\lambda_\infty\geq H(\nabla\phi(x_0),x_0)$.
Consequently $U^*$ is a periodic viscosity subsolution of
$H(\nabla U,x)=\lambda_\infty$. Compare it with the
same periodic solution $V$ at a maximum of $U^*-V$ using the
periodic doubling-of-variables argument, now with the subsolution
inequality for $U^*$ and the supersolution inequality for $V$.
The inequalities yield
$\overline\lambda_\infty(0)\le\lambda_\infty$.
The two comparisons have opposite directions; the variational
formula involving only $\inf\max H$ is not used to infer the
supersolution inequality.

Combining the two bounds gives
$\lambda_\infty=\overline\lambda_\infty(0)$. Every convergent
subsequence of the bounded family $(\lambda_L)$ has this limit;
hence $\lambda_L^*(0)\to\overline\lambda_\infty(0)$ as
$L\to\infty$.
\end{proof}

We next extend this convergence to all $\mu\in\mathcal I$ and remove
the compact-support assumption.

\begin{proof}[\textbf{Proof of Theorem~\ref{thm:large-period-limits}(i), completion}]
After the change of variables $x=L\widetilde x$, the $i$th component of
the exponentially weighted convolution is
\[
\int_{\mathbb R^N}J_i(z)e^{\mu\cdot z}
\phi_i\!\left(x-\frac zL\right)\,\mathrm dz.
\]
Suppose first that all kernels have compact support. For fixed $\mu$,
the weighted kernels $J_i(z)e^{\mu\cdot z}$ retain compact support and
positivity near zero, and have finite positive masses $M_i(\mu)$.
The proof of Lemma~\ref{L3.2} applies without unit-mass normalization:
the multiplication part remains $A-I$, and positivity of the weighted
kernels on a fixed ball yields the analogue of
\eqref{eq:large-period-uniform-spectral-gap}.
A uniform upper bound follows by testing
\eqref{def:lambdaLstar} with the constant vector and using the masses
$M_i(\mu)$.
Thus, for all sufficiently large $L$, a positive principal eigenfunction
exists and its normalized logarithm is uniformly bounded. With
$W_{i,L}=-L^{-1}\log\phi_{i,L}$, division by $\phi_{i,L}$ changes the
preceding integral into
\[
\int_{\mathbb R^N}J_i(z)
\exp\!\left(\mu\cdot z+L\left[W_{i,L}(x)
-W_{i,L}\!\left(x-\frac zL\right)\right]\right)\,\mathrm dz.
\]
Use the global minimizing and maximizing pairs on the
torus constructed in the preceding zero-weight proof, this time for
the weighted eigenfunctions. Since each pair is extremal over both
spatial points and component indices, the inequalities for the
logarithmic increments hold for every $z$ in the support of each
weighted kernel. In particular, no local-contact inequality is
extended to a nonlocal point without justification. For a smooth test function $\phi$ and $x_L\to x$,
\[
\mu\cdot z+L\bigl[\phi(x_L)-\phi(x_L-z/L)\bigr]
\longrightarrow(\mu+\nabla\phi(x))\cdot z
\]
uniformly for $z$ in the compact kernel supports. Repeating the two
half-relaxed-limit arguments in that proof therefore
gives respectively a periodic viscosity super-solution and sub-solution of
\begin{equation}\label{eq:tilted-periodic-HJ}
H(\mu+\nabla W(x),x)=c\quad\text{on }\mathbb T^N.
\end{equation}
Here $c$ is an arbitrary subsequential limit of $\lambda_L^*(\mu)$.
The comparison characterization of the additive constant in
\eqref{eq:tilted-periodic-HJ} yields
$c=\overline\lambda_\infty(\mu)$.

We now remove the compact-support assumption. Fix $\mu\in\mathcal I$.
Choose smooth cutoffs $0\leq\eta_n\uparrow1$ such that $\eta_n=1$ on
$B_n$ and $\eta_n=0$ outside $B_{2n}$, and put $J_{i,n}=\eta_nJ_i$.
The loss term remains $-\phi_i$; the truncated kernels are not normalized.
Their exponential moments are
\[
M_{i,n}(p):=\int_{\mathbb R^N}J_{i,n}(z)e^{p\cdot z}\,\mathrm dz.
\]
Write $H_n$, $\overline\lambda_{\infty,n}$ and $\lambda_{L,n}^*$ for
the corresponding quantities. For fixed $L,\mu$, put
\[
 t_n(\mu)=\max_i\int_{\mathbb R^N}(1-\eta_n(z))
 J_i(z)e^{\mu\cdot z}\,\mathrm dz\longrightarrow0.
\]
The weighted convolution operators differ in norm by at most
$t_n(\mu)$. Monotonicity immediately implies
$\lambda_{L,n}^*(\mu)\leq\lambda_L^*(\mu)$.
To obtain the reverse limiting inequality, fix $\delta>0$ and
choose $\boldsymbol\phi\in X^{++}$ satisfying
$\mathcal L_\mu\boldsymbol\phi\geq
(\lambda_L^*(\mu)-\delta)\boldsymbol\phi$;
this is possible because the periodic upper and lower values agree,
as established by \eqref{eq:weighted-CW-sandwich}
and Lemma~\ref{lem:approx}. With
$\alpha=\min_{i,x}\phi_i(x)>0$ and
$\beta=\max_{i,x}\phi_i(x)$, we have
\[
 (\mathcal L_{\mu,n}\boldsymbol\phi)_i
 \geq(\lambda_L^*(\mu)-\delta)\phi_i
      -\beta t_n(\mu)
 \geq\left(\lambda_L^*(\mu)-\delta
      -\frac{\beta}{\alpha}t_n(\mu)\right)\phi_i.
\]
The lower Collatz--Wielandt formula therefore shows that
$\liminf_n\lambda_{L,n}^*(\mu)
\geq\lambda_L^*(\mu)-\delta$.
Letting $\delta\downarrow0$ establishes the asserted spectral
convergence without appealing to continuity of the spectral bound
under arbitrary nonselfadjoint perturbations. For the Hamiltonians, use monotone convergence of the
moments and continuity of the Perron root; if a moment diverges, use
$H_n(p,x)\geq M_{i,n}(p)-1+a_{ii}(x)$.
Together with kernel ordering, these arguments yield
\begin{equation}\label{eq:monotone-kernel-truncation}
 \lambda_{L,n}^*(\mu)\uparrow\lambda_L^*(\mu),\qquad
 H_n(p,x)\uparrow H(p,x).
\end{equation}
Write $K\Subset\operatorname{int}\mathcal I$ when the closure of $K$
is compact and contained in $\operatorname{int}\mathcal I$.
Every compact $K\Subset\operatorname{int}\mathcal I$ is covered by
finitely many polytopes with vertices in $\operatorname{int}\mathcal I$.
The sum of the exponential functions associated with these vertices
dominates $e^{\mu\cdot z}$ for all $\mu\in K$. Hence
\begin{equation*}
 \max_i\sup_{\mu\in K}
 \int_{|z|>n}J_i(z)e^{\mu\cdot z}\,\mathrm dz\longrightarrow0.
\end{equation*}

\smallskip
\noindent Step 1. Convergence of the effective values.
Kernel ordering gives
\[
\overline\lambda_{\infty,n}(\mu)
\leq\overline\lambda_{\infty,n+1}(\mu)
\leq\overline\lambda_\infty(\mu)
\leq\max_xH(\mu,x)<\infty.
\]
Let $c_*=\lim_n\overline\lambda_{\infty,n}(\mu)$.
Choose mean-zero functions
$W_n\in C^1_{\rm per}(\mathbb T^N)$ satisfying
\[
 H_n(\mu+\nabla W_n(x),x)\leq\overline\lambda_{\infty,n}(\mu)+1/n
 \quad\hbox{for all }x\in\mathbb T^N.
\]
The retained kernels are positive on a fixed ball. Integrating each
kernel over a smaller ball centered in the direction $p/|p|$ and
using the diagonal lower bound for the Perron root gives constants
$b_1,b_2,b_3>0$, independent of all large $n$, such that
\begin{equation}\label{eq:truncated-Hamiltonian-coercivity}
 H_n(p,x)\geq b_1e^{b_2|p|}-b_3.
\end{equation}
Thus the gradients are uniformly bounded. After passage to a subsequence,
$W_n$ converges uniformly to a periodic Lipschitz function $W$.

Fix $k$ and a periodic mollifier $\rho_\delta$, a smooth nonnegative
averaging kernel of integral one supported within distance $\delta$
of zero on the torus. Convexity of
$H_k(\cdot,x)$ and its uniform continuity on the bounded gradient set
give, for $n\geq k$,
\[
 H_k\bigl(\mu+\nabla(W_n*\rho_\delta)(x),x\bigr)
 \leq\overline\lambda_{\infty,n}(\mu)+1/n+\omega_k(\delta),\qquad
 \omega_k(\delta)\longrightarrow0.
\]
First let $n\to\infty$, and then let $\delta\to0$ at a Lebesgue point
of $\nabla W$. Intersecting the full-measure sets for all $k$ and using
$H_k\uparrow H$ gives $H(\mu+\nabla W(x),x)\leq c_*$ almost everywhere.

On the sublevel set $H(p,y)\leq c_*$, the inequality
$H(p,y)\geq M_i(p)-1+a_{ii}(y)$ bounds all moments $M_i(p)$ uniformly.
Uniform continuity of $A$ and of the Perron root on the resulting
bounded set of matrices therefore gives a modulus of continuity
$\omega$, independent of $p$ and satisfying $\omega(\delta)\to0$ as
$\delta\to0$, such that
\[
 H(p,x)\leq H(p,y)+\omega(\delta)
 \quad\hbox{if }|x-y|\leq\delta,\ H(p,y)\leq c_*.
\]
The extended-valued function $H(\cdot,x)$ is lower semicontinuous by
Fatou's lemma and is convex. Jensen's inequality yields
\[
 H\bigl(\mu+\nabla(W*\rho_\delta)(x),x\bigr)
 \leq c_*+\omega(\delta),\qquad x\in\mathbb T^N.
\]
In particular, these mollified gradients are admissible. By the
definition of $\overline\lambda_\infty(\mu)$, followed by $\delta\to0$,
we obtain $\overline\lambda_\infty(\mu)\leq c_*$.
Thus
\begin{equation}\label{eq:effective-constant-truncation}
 \overline\lambda_{\infty,n}(\mu)\uparrow
 \overline\lambda_\infty(\mu).
\end{equation}
The compact-kernel case and kernel ordering now give
\[
 \liminf_{L\to\infty}\lambda_L^*(\mu)
 \geq\overline\lambda_{\infty,n}(\mu)\quad\hbox{for every }n,
\]
which proves the required lower bound.

\smallskip
\noindent Step 2. A positive test function for the upper bound.
Fix $\varepsilon>0$ and choose $V\in C^1_{\rm per}(\mathbb T^N)$ with
\[
 \sup_xH(\mu+\nabla V(x),x)\leq\overline\lambda_\infty(\mu)+\varepsilon.
\]
Let $S=\operatorname{conv}\{\mu+\nabla V(x):x\in\mathbb T^N\}$,
where $\operatorname{conv}$ denotes the convex hull, the set of all
finite convex combinations. It is a compact subset of the convex set
$\mathcal I$. Write $\operatorname{aff}S$ for the smallest affine
subspace containing $S$, and $\operatorname{relint}S$ for the interior
of $S$ relative to this subspace. Since
$\int_{\mathbb T^N}\nabla V(x)\,\mathrm dx=0$, a proper supporting
hyperplane through $\mu$ in $\operatorname{aff}S$ would force the whole
continuous range $\{\mu+\nabla V(x):x\in\mathbb T^N\}$ into that
hyperplane, contradicting its propriety. Hence
$\mu\in\operatorname{relint}S$.

Fix $0<\theta<1$, temporarily translate $\mu$ to zero, and work in
$\operatorname{aff}S$. If its dimension is positive, let $B$ be its
unit ball and choose $r>0$ with $rB\subset S$. Choose a finite
$\gamma$-net in $S$, that is, finitely many points of $S$ such that
every point of $S$ lies within distance $\gamma$ of one of them,
and let $P$ be their convex hull. The support function
$h_S(u):=\max_{q\in S}q\cdot u$ satisfies $h_S(u)\geq r$ for unit
vectors $u$. With $h_P$ defined in the same way,
$h_P(u)\geq h_S(u)-\gamma$. Thus
$h_P(u)\geq\theta h_S(u)$ when $\gamma<(1-\theta)r$, which gives
$\theta S\subset P\subset S$. In dimension zero, take $P=S$.
Translating back, we obtain a finite polytope satisfying
$\mu+\theta(S-\mu)\subset P\subset S$. Denote its vertices by
$q_1,\ldots,q_a\in\mathcal I$.
For all $L,x,z$,
\begin{align*}
 &\mu\cdot z+L\theta[V(x)-V(x-z/L)]=\int_0^1(\mu+\theta\nabla V(x-sz/L))\cdot z\,\mathrm ds
 \leq\max_{1\leq j\leq a}q_j\cdot z.
\end{align*}
The exponential of this expression is bounded by
$\sum_{j=1}^ae^{q_j\cdot z}$, an integrable majorant after multiplication
by any $J_i$. No interior assumption on $\mu$ is used here.

Put $p_\theta(x)=\mu+\theta\nabla V(x)$. The same majorant proves
continuity of $M_i(p_\theta(x))$. Let $\boldsymbol v_\theta(x)>0$ be
the continuous Perron vector of $D(p_\theta(x))+A(x)$, normalized by
$\sum_i v_{\theta,i}(x)=1$. Its components have a common positive lower
bound. For the positive periodic test function
\mbox{$\Phi_{L,i}(x)=e^{-L\theta V(x)}v_{\theta,i}(x)$},
the quotient of the rescaled weighted operator by $\Phi_{L,i}$ equals
\[
 \int_{\mathbb R^N}J_i(z)e^{\mu\cdot z+L\theta[V(x)-V(x-z/L)]}
 \frac{v_{\theta,i}(x-z/L)}{v_{\theta,i}(x)}\,\mathrm dz-1
 +\sum_j a_{ij}(x)
 \frac{v_{\theta,j}(x)}{v_{\theta,i}(x)}.
\]
For bounded $z$ the integrand converges uniformly in $x$, and the
preceding majorant controls the tails uniformly in $L,x$.
This quotient therefore converges uniformly to $H(p_\theta(x),x)$.
The upper characterization \eqref{def:lambdaLstar}, applied to
\eqref{eq:rescaled-weighted-operator}, and convexity of $H$ give
\[
 \limsup_{L\to\infty}\lambda_L^*(\mu)
 \leq\sup_xH(p_\theta(x),x)
 \leq\theta(\overline\lambda_\infty(\mu)+\varepsilon)
 +(1-\theta)\max_xH(\mu,x).
\]
Let $\theta\to1$ and $\varepsilon\to0$. This proves the upper bound
also for $\mu\in\partial\mathcal I$, where $\partial\mathcal I$ denotes
the boundary of $\mathcal I$, and when $\mathcal I$ has empty
interior. Finally, the spectral functions are finite and convex on
$\operatorname{int}\mathcal I$. Pointwise convergence of finite convex
functions is locally uniform there, proving the last assertion.
\end{proof}

\begin{proof}[\textbf{Proof of Theorem~\ref{thm:large-period-limits}(ii)}]
Part~(i) at zero weight and
Theorem~\ref{TH2.1}(iii) give the limit of $\lambda'(L)$.
It remains to pass to the minimum over exponential weights.
	By Theorem~\ref{TH2.1}(ii),
	\[
	\lambda^*(L)=\min_{\mu\in\mathbb R^N}\lambda_L^*(\mu).
	\]
	Let $J_{i,n}=J_i\eta_n$ be the monotone truncations used in
	the proof of part~(i). We first show that
	\begin{equation}\label{eq:effective-constant-truncation-all-weights}
	\overline\lambda_{\infty,n}(\mu)\uparrow
	\overline\lambda_\infty(\mu)
	\quad\hbox{for every }\mu\in\mathbb R^N.
	\end{equation}
	For $\mu\in\mathcal I$ this is
	\eqref{eq:effective-constant-truncation}. If $\mu\notin\mathcal I$
	and $\overline\lambda_{\infty,n}(\mu)$ were bounded, the functions $W_n$ and the mollification argument used to prove
\eqref{eq:effective-constant-truncation} would give
	a periodic Lipschitz function $W$ with
	$H(\mu+\nabla W(x),x)\leq C$ almost everywhere.
	The diagonal lower bound for the Perron root would then bound all
	$M_i(\mu+\nabla W(x))$ uniformly. Jensen's inequality and
	$\int_{\mathbb T^N}\nabla W(x)\,\mathrm dx=0$ would imply
	\[
	M_i(\mu)\leq
	\int_{\mathbb T^N}M_i(\mu+\nabla W(x))\,\mathrm dx<\infty
	\quad\hbox{for every }i,
	\]
	a contradiction. This proves the assertion also outside $\mathcal I$.

	Joint convexity of the defining expression in $\mu$ and $W$ gives
	convexity of each finite function $\overline\lambda_{\infty,n}$. Their
	pointwise supremum $\overline\lambda_\infty$ is therefore lower
	semicontinuous. The lower bound \eqref{eq:truncated-Hamiltonian-coercivity}
holds with constants independent of all large $n$. Integrating in $x$ and using Jensen's inequality
	shows that $\overline\lambda_{\infty,n}(\mu)\geq b_1e^{b_2|\mu|}-b_3$.
	The bound \eqref{eq:weighted-moment-lower-bound} and positivity of
the kernels on a fixed ball give the same type of lower bound for
$\lambda_L^*(\mu)$, independently of $L$.
	The lower semicontinuity and coercive bound imply that
	$\overline\lambda_\infty$ attains its finite minimum. The
	constant-vector upper test bounds $\lambda_L^*(0)$ independently of
	$L$, and minimality gives $\lambda_L^*(\mu_L)\leq\lambda_L^*(0)$.
	Together with the uniform lower bound above, this confines the
	minimizers $\mu_L$ of $\lambda_L^*(\cdot)$ to a common compact set.
	Consider any sequence $L_k\to\infty$. After passage to a subsequence,
	$\mu_{L_k}\to\mu_0$. For every fixed $n$, the compact-kernel limit
	and convexity imply locally uniform convergence
	$\lambda_{L,n}^*(\cdot)\to\overline\lambda_{\infty,n}$ on $\mathbb R^N$. Kernel ordering
	therefore yields
	\begin{equation*}
	\liminf_{k\to\infty}\lambda^*(L_k)
	\geq\lim_{k\to\infty}\lambda_{L_k,n}^*(\mu_{L_k})
	=\overline\lambda_{\infty,n}(\mu_0).
	\end{equation*}
	Letting $n\to\infty$ gives
	\[
	\liminf_{k\to\infty}\lambda^*(L_k)
	\geq\overline\lambda_\infty(\mu_0)
	\geq\min_\mu\overline\lambda_\infty(\mu).
	\]
	Conversely, choose a minimizer $\mu_*\in\mathcal I$ of
	$\overline\lambda_\infty$.
	Part~(i) gives
	\[
	\limsup_{L\to\infty}\lambda^*(L)
	\leq\lim_{L\to\infty}\lambda_L^*(\mu_*)
	=\overline\lambda_\infty(\mu_*).
	\]
	Since the sequence was arbitrary, this proves the large-period
	limit without requiring an interior minimizer of the moment domain.
	
\end{proof}

\subsubsection{Scalar equations}

We first prove the Laplace asymptotic independently, using the following
scalar semiclassical estimate. We then combine this asymptotic with
Fourier estimates to treat general symmetric kernels.

\begin{lemma}\label{lem:semiclassical-circle}
Let $V\in C^2(\mathbb T)$ be nonnegative, with zero set
$\{x_1,\ldots,x_r\}$ consisting of finitely many distinct points.
Suppose that $\kappa_j=V''(x_j)>0$ for every $j$, and set
$\omega_* =\min_{1\leq j\leq r}\sqrt{\kappa_j/2}$.
Let $E_\varepsilon$ denote the lowest periodic eigenvalue of
$-\varepsilon^2\partial_{xx}+V$. Then the following statements hold:
\begin{enumerate}
\item[(i)] As $\varepsilon\downarrow0$,
\[
E_\varepsilon=\varepsilon\omega_*+o(\varepsilon).
\]
\item[(ii)] Every ground state, meaning an eigenfunction for this lowest eigenvalue,
normalized to have $L^2$ norm one, has mass tending to zero
outside any fixed neighborhood of the finite zero set as $\varepsilon\downarrow0$.
\end{enumerate}
\end{lemma}

\begin{proof}[Proof of Lemma~\ref{lem:semiclassical-circle}]
Choose $x_j$ with $\sqrt{\kappa_j/2}=\omega_*$ and use local
coordinates centered at $x_j$. The normalized Gaussian
\[
g_\varepsilon(x)=\left(\frac{\omega_*}{\pi\varepsilon}\right)^{1/4}
\exp\left(-\frac{\omega_*(x-x_j)^2}{2\varepsilon}\right)
\]
has kinetic and quadratic potential energies each equal to
$\varepsilon\omega_*/2$. Multiply it by a fixed smooth cutoff
equal to one near $x_j$, and normalize again. The cutoff error is
exponentially small. Since
$V(x)=\omega_*^2(x-x_j)^2+o(|x-x_j|^2)$, rescaling
$x-x_j=\sqrt\varepsilon\,y$ shows that the potential remainder is
$o(\varepsilon)$. The resulting Rayleigh quotient is therefore
$\varepsilon\omega_*+o(\varepsilon)$, proving the upper bound.
For the converse, choose disjoint fixed neighborhoods $U_j$ of the
points $x_j$ and a smooth quadratic partition of unity
$\sum_{j=0}^r\chi_j^2=1$, where each $\chi_j$ for $j\geq1$
is supported inside $U_j$, while $\chi_0$ is supported away from
the zero set. Given $0<\eta<\frac12\min_j\kappa_j$, make the neighborhoods small enough
that $V(x)\geq\bigl(\kappa_j/2-\eta\bigr)|x-x_j|^2$
for $x\in U_j$.
For a periodic $H^1$ function $u$ (a square-integrable function with
square-integrable weak derivative), differentiation of
$\sum_j\chi_j^2=1$ gives the exact localization identity
\begin{equation}\label{eq:partition-energy-identity}
\sum_{j=0}^r\| (\chi_j u)'\|_2^2
=\|u'\|_2^2+\int_{\mathbb T}\sum_{j=0}^r|\chi_j'|^2|u|^2.
\end{equation}
The cross terms vanish because $\sum_j\chi_j\chi_j'=0$.
For a compactly supported $w\in H^1(\mathbb R)$ and $\alpha>0$,
expanding $0\leq\|\varepsilon w'+\alpha(x-x_j)w\|_2^2$ and
integrating the cross term by parts yields
\begin{equation}\label{eq:harmonic-oscillator-bound}
\varepsilon^2\|w'\|_2^2
+\alpha^2\|(x-x_j)w\|_2^2
\geq\varepsilon\alpha\|w\|_2^2.
\end{equation}
Equality holds for a Gaussian before cutoff. Applying \eqref{eq:harmonic-oscillator-bound}
to $w=\chi_j u$ in local coordinates and using
\eqref{eq:partition-energy-identity}, we obtain
\[
\begin{aligned}
\varepsilon^2\|u'\|_2^2+\int V|u|^2
&\geq \varepsilon\min_j\sqrt{\kappa_j/2-\eta}
 \sum_{j=1}^r\|\chi_j u\|_2^2
 +c\|\chi_0u\|_2^2-C\varepsilon^2\|u\|_2^2\\
&\geq \left(\varepsilon\min_j\sqrt{\kappa_j/2-\eta}
 -C\varepsilon^2\right)\|u\|_2^2
\end{aligned}
\]
for all small $\varepsilon$; here $c>0$ depends only on the
chosen neighborhoods. The lower bound follows by first taking
$\varepsilon\downarrow0$ and then $\eta\downarrow0$.
Finally, the normalized ground-state energy is $O(\varepsilon)$.
Since $V$ is bounded away from zero on any closed set disjoint
from its zeros, the concentration assertion follows.
\end{proof}

\begin{proof}[Proof of Corollary~\ref{thm:large-period-harmonic}]
Put $\varepsilon=L^{-1}$ and
$K_\varepsilon=(I-\varepsilon^2\partial_{xx})^{-1}$ on the unit
torus, and denote the periodic spectral bound $\lambda_L^*(0)$ of
$K_\varepsilon-I+a(x)$ by $\lambda_\varepsilon$.
Since $K_\varepsilon-I$ is nonpositive and self-adjoint,
$\lambda_\varepsilon\leq a_*$. The Fourier-multiplier inequality
$I-K_\varepsilon\leq-\varepsilon^2\partial_{xx}$ and the
$O(\varepsilon)$ potential loss of a normalized Gaussian of width
$\sqrt\varepsilon$ give
$\lambda_\varepsilon\geq a_*-C\varepsilon$ through the Rayleigh
quotient. Thus $\lambda_\varepsilon\to a_*$ and
$\lambda_\varepsilon>a_*-1$ for small $\varepsilon$.
Lemma~\ref{L11}(i) therefore supplies a positive periodic eigenfunction.
Let $\phi_\varepsilon>0$ be this eigenfunction, normalized in $L^2$,
and set $q_\varepsilon(x)=1+\lambda_\varepsilon-a(x)>0$
and $v_\varepsilon=q_\varepsilon\phi_\varepsilon$.
The eigenvalue equation $K_\varepsilon\phi_\varepsilon
=q_\varepsilon\phi_\varepsilon$ is equivalent to
\begin{equation}\label{eq:Laplace-resolvent-schrodinger}
-\varepsilon^2v_\varepsilon''
+\frac{\lambda_\varepsilon-a(x)}
{1+\lambda_\varepsilon-a(x)}v_\varepsilon=0.
\end{equation}
Consequently zero is the lowest eigenvalue of the periodic
Schr\"odinger operator in \eqref{eq:Laplace-resolvent-schrodinger}.
Write $V(x)=a_*-a(x)$ and
$\lambda_\varepsilon=a_*-\varepsilon c_\varepsilon$;
the preceding bounds give $0\leq c_\varepsilon\leq C$.
The potential becomes
\[
F_{\varepsilon,c}(x)
=\frac{V(x)-\varepsilon c}{1+V(x)-\varepsilon c}
=\frac{V(x)}{1+V(x)}
-\frac{\varepsilon c}{(1+V(x))(1+V(x)-\varepsilon c)}.
\]
For bounded $c$, the second term is uniformly $O(\varepsilon)$.
The first term has precisely the nondegenerate minima
$x_1,\ldots,x_r$, with second derivatives $\kappa_1,\ldots,\kappa_r$. Lemma~\ref{lem:semiclassical-circle}
therefore shows that the lowest eigenvalue of
$-\varepsilon^2\partial_{xx}+V/(1+V)$ is
$\varepsilon\omega_*+o(\varepsilon)$.
The ground states of both this operator and its bounded-$c$
perturbations concentrate near the finite set $\{x_1,\ldots,x_r\}$ of global maxima of $a$: the energy upper bound is
$O(\varepsilon)$, while the nonnegative potential $V/(1+V)$
is bounded away from zero outside any neighborhood of this set.
It follows, uniformly for bounded $c$, that the expectation of
$[(1+V)(1+V-\varepsilon c)]^{-1}$ in the corresponding normalized
ground state tends to one. The min--max principle, applied in both
directions using the two ground states, now gives
\[
\inf\sigma\bigl(-\varepsilon^2\partial_{xx}
+F_{\varepsilon,c}\bigr)
=\varepsilon\bigl(\omega_*-c\bigr)+o(\varepsilon)
\]
uniformly for bounded $c$. At $c=c_\varepsilon$ the left-hand
side is zero, so $c_\varepsilon\to\omega_*$.
The scalar problem with an even kernel is a special case of
Proposition~\ref{prop:detailed-balance-rigidity}. That proposition gives
$\lambda^*(L)=\lambda'(L)=\lambda_L^*(0)$ directly.
The cosine formula in Corollary~\ref{thm:large-period-harmonic} follows from $-a''(0)=4\pi^2\delta$.
\end{proof}

With the Laplace asymptotic established, we turn to general symmetric kernels.

\begin{proof}[\textbf{Proof of Theorem~\ref{thm:universal-kernel-large-period}}]
Put $\varepsilon=L^{-1}$ and
$K_\varepsilon f(x)=\int_{\mathbb R}J(z)f(x-\varepsilon z)\,\mathrm dz$
on the unit torus, and write
$\widehat J(\xi)=\int_{\mathbb R}J(z)e^{-\mathrm i\xi z}\,\mathrm dz$.
Since $J$ is even, $K_\varepsilon$ is
self-adjoint on $L^2(\mathbb T)$ and its multiplier at frequency
$2\pi k$ is $\widehat J(2\pi\varepsilon k)$. Its spectrum
accumulates only at zero, so $K_\varepsilon$ is compact on
$L^2(\mathbb T)$. Moreover, (J) implies that the periodized kernel
$\varepsilon^{-1}\sum_{k\in\mathbb Z}J((x+k)/\varepsilon)$
converges uniformly on $\mathbb T$ and is continuous. Thus
$K_\varepsilon$ is compact on both $C(\mathbb T)$ and
$L^2(\mathbb T)$ and maps $L^2(\mathbb T)$ into $C(\mathbb T)$.
Multiplication by $a-1$ has the same essential spectral upper edge
$a_*-1$ in both spaces, and compact perturbations preserve this
edge \cite[Chapter~IV, Theorem~5.35]{Kato1995}.
In the half-plane $\operatorname{Re}\lambda>a_*-1$, spectral points
are isolated eigenvalues, and
$u=(1+\lambda-a)^{-1}K_\varepsilon u$ makes every $L^2$
eigenfunction continuous. The upper spectral bounds therefore
coincide, including when they equal the multiplication edge; write
$\lambda_\varepsilon$ for this common value. The upper spectral bound of
$K_\varepsilon-I+a$ is characterized by
\begin{equation}\label{eq:universal-Rayleigh}
 a_*-\lambda_\varepsilon
 =\inf_{\|f\|_2=1}
 \left\{\int_{\mathbb T}(a_*-a)|f|^2\,\mathrm dx
 +\langle f,(I-K_\varepsilon)f\rangle\right\}.
\end{equation}
Here $I-K_\varepsilon$ is nonnegative, because the Fourier
transform of an even probability density is real and at most one.
The elementary inequality $1-\cos t\leq t^2/2$ yields
\begin{equation}\label{eq:kernel-second-moment-upper}
 0\leq1-\widehat J(\xi)\leq D_J\xi^2.
\end{equation}
Let $j$ realize $\Omega_J$ and test \eqref{eq:universal-Rayleigh}
with a smooth cutoff of the normalized Gaussian proportional to
$\exp[-\sqrt{\kappa_j/(2D_J)}(x-x_j)^2/(2\varepsilon)]$.
For this Gaussian the sum of the kinetic energy
$\varepsilon^2D_J\|f'\|_2^2$ and the quadratic potential energy
$\frac{\kappa_j}{2}\|(x-x_j)f\|_2^2$ equals
$\varepsilon\sqrt{D_J\kappa_j/2}$; cutoff and the $C^2$ Taylor
remainder contribute $o(\varepsilon)$. Thus
\eqref{eq:kernel-second-moment-upper} gives
\begin{equation}\label{eq:universal-upper-energy}
 a_*-\lambda_\varepsilon
 \leq\varepsilon\Omega_J+o(\varepsilon).
\end{equation}
In particular $\lambda_\varepsilon>a_*-1$ for small
$\varepsilon$. The compact perturbation of multiplication by
$a-1$ therefore has an isolated top eigenvalue. Positivity of the
convolution and irreducibility on the torus supply a positive
continuous eigenfunction. The minimum-ratio comparison with positive
periodic test functions then identifies
$\lambda_\varepsilon=\lambda_L^*(0)$.

We next establish the converse bound without a local Taylor
expansion of the convolution operator. For each $\eta\in(0,1)$
there exists $C_\eta>0$ such that, for every $\xi\in\mathbb R$,
\begin{equation}\label{eq:kernel-Pade-comparison}
 \frac{(1-\eta)D_J\xi^2}{1+C_\eta\xi^2}
 \leq1-\widehat J(\xi)\leq D_J\xi^2.
\end{equation}
Dominated convergence and the second moment give
$(1-\widehat J(\xi))/\xi^2\to D_J$ as $\xi\to0$.
Since $J$ is a nonnegative probability density,
$\widehat J(\xi)<1$ whenever $\xi\ne0$; otherwise
$\cos(\xi z)=1$ on a set of full $J(z)\,\mathrm dz$ measure,
which is impossible for a density. Also
$\widehat J(\xi)\to0$ as $|\xi|\to\infty$ by the
Riemann--Lebesgue lemma. Thus $1-\widehat J$ has a positive
minimum outside any fixed neighborhood of zero. Taking
$C_\eta$ sufficiently large establishes
\eqref{eq:kernel-Pade-comparison} there, while the second-moment
limit supplies it near zero.

Set $d_\eta=(1-\eta)D_J/C_\eta$ and let
$G_{\varepsilon,\eta}=(I-C_\eta\varepsilon^2\partial_{xx})^{-1}$.
The lower estimate in \eqref{eq:kernel-Pade-comparison} gives
the quadratic-form inequality
\[
 I-K_\varepsilon\ \geq\
 d_\eta(I-G_{\varepsilon,\eta}).
\]
Consequently, if $\lambda_{\varepsilon,\eta}$ is the top
spectral value of
$d_\eta(G_{\varepsilon,\eta}-I)+a$, then
$\lambda_\varepsilon\leq\lambda_{\varepsilon,\eta}$.
This auxiliary operator is a scaled Laplace-dispersal operator.
Applying the Laplace asymptotic of
Corollary~\ref{thm:large-period-harmonic}, proved independently above,
to the potential
$a/d_\eta$ at period
$(\sqrt{C_\eta}\varepsilon)^{-1}$, and multiplying
its eigenvalue by $d_\eta$, we obtain
\[
 \lambda_{\varepsilon,\eta}
 =a_*-\varepsilon\min_j
 \sqrt{\frac{d_\eta C_\eta\kappa_j}{2}}
 +o_\eta(\varepsilon)
 =a_*-\varepsilon\sqrt{1-\eta}\,\Omega_J
 +o_\eta(\varepsilon).
\]
Combining this with \eqref{eq:universal-upper-energy} yields
\[
 \sqrt{1-\eta}\,\Omega_J
 \leq\liminf_{\varepsilon\downarrow0}
 \frac{a_*-\lambda_\varepsilon}{\varepsilon}
 \leq\limsup_{\varepsilon\downarrow0}
 \frac{a_*-\lambda_\varepsilon}{\varepsilon}
 \leq\Omega_J.
\]
Letting $\eta\downarrow0$ proves the spectral expansion.

For normalized eigenfunctions, both terms in
\eqref{eq:universal-Rayleigh} are nonnegative, so
$\int(a_*-a)|u_L|^2=O(L^{-1})$. For every open neighborhood $V$ of the maxima with nonempty complement, the potential
$a_*-a$ is bounded below by $\delta_V>0$ on $\mathbb T\setminus V$.
The expansion \eqref{eq:universal-large-period-asymptotic}
and the variational identity \eqref{eq:universal-Rayleigh} yield
\[
 \delta_V\int_{\mathbb T\setminus V}|u_L|^2
 \leq \int_{\mathbb T}(a_*-a)|u_L|^2
 \leq a_*-\lambda_\varepsilon
 =\varepsilon\Omega_J+o(\varepsilon),
\]
which proves the quantitative concentration estimate.
Finally, the scalar even-kernel case of
Proposition~\ref{prop:detailed-balance-rigidity} gives
$\lambda^*(L)=\lambda'(L)=\lambda_L^*(0)$.
\end{proof}

\par\medskip\noindent
For the Laplace density $\tfrac12e^{-|z|}$,
$D_J=1$, and \eqref{eq:universal-large-period-asymptotic}
recovers Corollary~\ref{thm:large-period-harmonic}. For the
Gaussian density of variance $\sigma^2$, the coefficient is
$\min_j\tfrac12\sqrt{\sigma^2\kappa_j}$.
The comparison argument depends on the full kernel only through
its second moment at leading order. It does not justify an
$O(L^{-2})$ remainder without further assumptions.
\par\medskip

\subsubsection{Reversible systems}

We now use a Rayleigh-quotient comparison to transfer the scalar
estimate to the system.

\begin{proof}[\textbf{Proof of Theorem~\ref{thm:reversible-system-reduction}}]
Since $S$ is a $C^1$ symmetric cooperative irreducible
matrix, its largest eigenvalue is simple. Compactness of the torus
and the uniform spectral separation imply that the positive normalized
eigenvector $v$ is $C^1$ and has a finite Lipschitz constant.
The constant similarity transform $D^{1/2}$ commutes with the common
componentwise convolution, so the system is similar to the
self-adjoint operator $K_L-I+S(x)$ on
$L^2(\mathbb T;\mathbb R^m)$.
The compact-perturbation argument in the scalar proof of
Theorem~\ref{thm:universal-kernel-large-period} also identifies its
spectral bound on $C$ with that on $L^2$: the multiplication edge
is now $\max_xs(x)-1$, and an $L^2$ eigenfunction above this edge
is continuous by $w=(\lambda I-S+I)^{-1}K_Lw$.
If there is no spectral point above the edge, both bounds equal it.
For vector functions $w$, put
\[
 \mathcal E_L(w)=\frac12\int_{\mathbb T}\int_{\mathbb R}
 J(z)|w(x)-w(x-z/L)|^2\,\mathrm dz\,\mathrm dx.
\]
The corresponding scalar form has the same definition with absolute
values. The spectral variational principle gives
\begin{equation*}
 \Lambda_L=\sup_{\|w\|_2=1}
 \left\{\int_{\mathbb T}\langle S(x)w(x),w(x)\rangle
 \,\mathrm dx-\mathcal E_L(w)\right\}.
\end{equation*}
For any $w$ let $f(x)=|w(x)|$. The matrix Rayleigh inequality and
$|f(x)-f(y)|\leq|w(x)-w(y)|$ show that
\[
 \int\langle Sw,w\rangle-\mathcal E_L(w)
 \leq\int s(x)f(x)^2\,\mathrm dx-\mathcal E_L(f).
\]
Taking suprema gives $\Lambda_L\leq\Lambda_L^{\rm sc}(s)$.
Conversely, let $f\geq0$, $\|f\|_2=1$, and take $w(x)=f(x)v(x)$.
Then $\|w\|_2=1$ and
$\langle S(x)w(x),w(x)\rangle=s(x)f(x)^2$.
The normalization $|v|=1$ gives the exact form identity
\[
 \mathcal E_L(fv)-\mathcal E_L(f)
 =\int_{\mathbb T}\int_{\mathbb R}
 J(z)f(x)f(x-z/L)
 \bigl[1-\langle v(x),v(x-z/L)\rangle\bigr]
 \,\mathrm dz\,\mathrm dx.
\]
Since
\[
 1-\langle v(x),v(y)\rangle
 =\tfrac12|v(x)-v(y)|^2
 \leq\tfrac12\operatorname{Lip}(v)^2|x-y|_{\mathbb T}^2
 \leq\tfrac12\operatorname{Lip}(v)^2|z/L|^2
\]
when $y=x-z/L$, and since
$\int_{\mathbb T}f(x)f(x-z/L)\,\mathrm dx\leq1$,
we get
\[
 0\leq\mathcal E_L(fv)-\mathcal E_L(f)
 \leq D_J\operatorname{Lip}(v)^2L^{-2}.
\]
The scalar Rayleigh supremum may be taken over nonnegative $f$:
replacing a real-valued test by its absolute value leaves the
potential term unchanged and does not increase its energy.
Taking the supremum over such tests proves
\eqref{eq:reversible-finite-period-comparison} without assuming that
a principal eigenfunction exists at finite $L$.

Under the additional assumptions, $s\in C^2(\mathbb T)$; the
scalar Theorem~\ref{thm:universal-kernel-large-period}, applied to
the potential $s$, gives
\[
 \Lambda_L^{\rm sc}(s)=s_*-\Omega_J L^{-1}+o(L^{-1}).
\]
The error in \eqref{eq:reversible-finite-period-comparison} is
$O(L^{-2})$, which proves
\eqref{eq:reversible-system-localization}.
Since $\Lambda_L>s_*-1$ for large $L$, the compact convolution
perturbation of multiplication by $S-I$ has an isolated top
eigenvalue. Positivity and irreducibility identify its strictly
positive eigenfunction; transforming back by $D^{-1/2}$ yields $u_L$.

For such an $L$, the Rayleigh identity reads
\begin{equation}\label{eq:reversible-defect-identity}
 s_*-\Lambda_L=
 \int_{\mathbb T}(s_*-s(x))|w_L|^2\,\mathrm dx
 +\int_{\mathbb T}
 \langle(s(x)I-S(x))w_L,w_L\rangle\,\mathrm dx
 +\mathcal E_L(w_L).
\end{equation}
All three terms on the right are nonnegative.
The first term and $s_*-\Lambda_L=\Omega_J L^{-1}+o(L^{-1})$
prove \eqref{eq:reversible-concentration}.
Finally, the symmetry and simplicity of the Perron eigenvalue imply,
pointwise in $x$,
\[
 \langle(s(x)I-S(x))w,w\rangle
 \geq\bigl[s(x)-\lambda_2(S(x))\bigr]
 |(I-\Pi(x))w|^2.
\]
The continuous positive gap has a positive minimum $\gamma$.
Integration in the second term of
\eqref{eq:reversible-defect-identity} proves
\eqref{eq:reversible-adiabatic-alignment}.
\end{proof}

\par\medskip\noindent
If the normalized Perron vector $v$ is constant, then
\eqref{eq:reversible-finite-period-comparison} yields the exact
identity $\Lambda_L=\Lambda_L^{\rm sc}(s)$ for every $L>0$.
The fixed diagonal symmetrizer and the common kernel are essential to
the direct self-adjoint comparison used above. The theorem makes no
claim for unrelated component kernels or for coupling matrices that
cannot be symmetrized by the same constant diagonal matrix.
\par\medskip

\section{Small-period spectral asymptotics}
\label{sec:quantitative-period}

We first establish the small-period averaging limit in arbitrary
spatial dimension, retaining the unit-period variables of
Section~\ref{sec:large-period-asymptotics}. We then obtain quantitative
expansions and threshold separation. The smooth-kernel result holds in
arbitrary dimension; the Bessel and Laplace expansions are
one-dimensional and require the stated isolated-eigenvalue hypotheses.
The proofs and example calculations follow the main results, and speed
corrections are given in Section~\ref{subsec:quantitative-speeds}.

\subsection{Main results}

We state the averaging limit, followed by the expansions for each
kernel class and examples of nonzero corrections and threshold gaps.

\subsubsection{Averaging limit}

For the small-period limit, we introduce the averaged operator
\begin{equation}\label{eq:averaged-operator}
\widehat{\mathcal{L}}_{\mu}[\boldsymbol{\varphi}](x) = \widehat{\mathcal{J}}_{\mu}[\boldsymbol{\varphi}](x) - \boldsymbol{\varphi}(x) + A(x)\boldsymbol{\varphi}(x),
\qquad x\in\mathbb{R}^N,
\end{equation}
where
\[
\bigl(\widehat{\mathcal J}_{\mu}[\boldsymbol{\varphi}]\bigr)_i(x)
=M_i(\mu)\int_{\mathbb T^N}\varphi_i(y)\,\mathrm dy,
\qquad i\in\Lambda.
\]
Then the generalized principal eigenvalue \(\widehat{\lambda}(\mu)\) of \(\widehat{\mathcal{L}}_{\mu}\) is defined by
\begin{equation}\label{eq:averaged-generalized-eigenvalue}
\begin{aligned}
\widehat\lambda(\mu)
&=\inf\{\lambda\in\mathbb R:\exists\boldsymbol\varphi\in X^{++},\
\widehat{\mathcal L}_\mu\boldsymbol\varphi\leq\lambda\boldsymbol\varphi\}\\
&=\sup\{\lambda\in\mathbb R:\exists\boldsymbol\varphi\in X^{++},\
\widehat{\mathcal L}_\mu\boldsymbol\varphi\geq\lambda\boldsymbol\varphi\}.
\end{aligned}
\end{equation}
The equality follows from Lemma~\ref{lem:approx} and the argument
preceding Lemma~\ref{lem:conv-avg}.
For $\mu\notin\mathcal I$, set $\widehat\lambda(\mu)=+\infty$.

\begin{theorem}\label{thm:small-period-limits}
Assume that (J) and (A) hold, with $L=L_1=\cdots=L_N$.
Then, as $L\to0$,
\[
\lambda'(L)=\lambda_L^*(0)\longrightarrow\widehat\lambda(0),
\qquad
\lambda^*(L)\longrightarrow
\min_{\mu\in\mathbb R^N}\widehat\lambda(\mu).
\]
\end{theorem}

\subsubsection{Smooth kernels}

We begin with smooth weighted kernels. Fourier estimates give the
following convergence result, faster than any fixed power of the period.

\begin{theorem}\label{thm:smooth-superalgebraic}
Assume that (J) and (A) hold. Let
$K\Subset\operatorname{int}\mathcal I$ be compact.
Suppose that, for every integer $r\geq0$,
\[
\sup_{\mu\in K}\sum_{i=1}^m\sum_{|\alpha|\leq r}
\bigl\|\partial_z^\alpha(J_i(z)e^{\mu\cdot z})\bigr\|_{L^1(\mathbb R^N)}<\infty.
\]
Suppose also that $\widehat\lambda(\mu)$ is an isolated algebraically simple
principal eigenvalue of $\widehat{\mathcal L}_\mu$ for every $\mu\in K$,
with a uniform spectral separation from the rest of the spectrum.
Then, for every $q>0$,
\[
\sup_{\mu\in K}|\lambda_L^*(\mu)-\widehat\lambda(\mu)|
=O_{K,q}(L^q)\qquad(L\downarrow0).
\]
In particular, the coefficients of $L$ and $L^2$ both vanish.
If the displayed integrability condition also holds after any fixed
number of derivatives with respect to $\mu$, the same estimate holds
for those derivatives of the simple eigenvalue branch.
\end{theorem}

\subsubsection{Bessel kernels}
\label{subsec:Bessel-period}

We compute fractional-order spectral corrections and threshold
gaps. Here and below
$P$ denotes componentwise averaging over $\mathbb T=\mathbb R/\mathbb Z$.

Fix $\frac12<s<1$ and positive numbers
$\sigma_1,\ldots,\sigma_m$. For $i=1,\ldots,m$ let
\begin{equation}\label{eq:Bessel-mixture-kernel}
 J_{s,i}(z)=\frac1{\Gamma(s)}\int_0^\infty
 \frac{t^{s-1}e^{-t}}{\sqrt{4\pi\sigma_i^2t}}
 \exp\left(-\frac{z^2}{4\sigma_i^2t}\right)\,\mathrm dt.
\end{equation}
Here $\Gamma(s)=\int_0^\infty t^{s-1}e^{-t}\,\mathrm dt$ is the
Euler gamma function. These are continuous, strictly positive, even probability densities
and satisfy (J). Their Fourier transforms and exponential moments are
\begin{equation}\label{eq:Bessel-symbol}
 \int_{\mathbb R}J_{s,i}(z)e^{-\mathrm i\xi z}\,\mathrm dz
 =(1+\sigma_i^2\xi^2)^{-s},\qquad
 M_i(\mu)=(1-\sigma_i^2\mu^2)^{-s},
\end{equation}
where $|\mu|<\sigma_i^{-1}$. Introduce the bounded operators on
$C(\mathbb T,\mathbb C^m)$
\[
 H_s=(-\partial_{xx})^{-s}(I-P),\qquad
 U_s=\partial_x(-\partial_{xx})^{-s-1}(I-P),\qquad
 D_s=\operatorname{diag}(\sigma_1^{-2s},\ldots,\sigma_m^{-2s}),
\]
acting componentwise. We use the periodic Fourier coefficients
$\widehat f(k)=\int_0^1 f(x)e^{-2\pi\mathrm i kx}\,\mathrm dx$.
Their Fourier multipliers for $k\ne0$ are,
respectively, $|2\pi k|^{-2s}$ and
$\mathrm i(2\pi k)|2\pi k|^{-2s-2}$. 

We first specify the coefficients of the expansion. At weights
where $\widehat\lambda(\mu)>s^*$ is algebraically simple, let
$\phi_\mu,\psi_\mu>0$ be right and left eigenfunctions of
$T_\mu=\widehat{\mathcal L}_\mu$, normalized by
$\langle\psi_\mu,\phi_\mu\rangle=1$, where
$\langle\psi,f\rangle=\int_{\mathbb T}\psi(x)^{\mathsf T}f(x)\,\mathrm dx$
is the pairing with the real left eigenfunction, extended linearly
to complex-valued $f$.
The spectral projection and reduced resolvent are
\begin{equation}\label{eq:averaged-projection-resolvent}
 \Pi_\mu f=\phi_\mu\langle\psi_\mu,f\rangle,\qquad
 S_\mu=(I-\Pi_\mu)
 \bigl((\widehat\lambda(\mu)-T_\mu)|_{\operatorname{ran}(I-\Pi_\mu)}\bigr)^{-1}
 (I-\Pi_\mu).
\end{equation}
Here $\operatorname{ran}(I-\Pi_\mu)$ denotes the range of $I-\Pi_\mu$;
$S_\mu$ is the inverse of $\widehat\lambda(\mu)-T_\mu$ on that
complementary space, extended by zero on the spectral eigenspace.
With the bounded operator
\[
 Z_{s,\mu}=-s\,\operatorname{diag}\!\left(
 \sigma_i^{-2s}\bigl[\sigma_i^{-2}+(2s+1)\mu^2\bigr]
 \right)_{i=1}^m
 (-\partial_{xx})^{-s-1}(I-P),
\]
define
\begin{equation}\label{eq:Bessel-expansion-coefficients}
 \begin{aligned}
 b_s(\mu)&=\langle\psi_\mu,D_sH_s\phi_\mu\rangle,\qquad
 c_s(\mu)=-2s\mu\langle\psi_\mu,D_sU_s\phi_\mu\rangle,\\
 q_s(\mu)&=\langle\psi_\mu,D_sH_sS_\mu D_sH_s\phi_\mu\rangle,\qquad
 r_s(\mu)=\langle\psi_\mu,Z_{s,\mu}\phi_\mu\rangle.
 \end{aligned}
\end{equation}

\begin{theorem}\label{thm:Bessel-fractional-system}\label{thm:Bessel-refined-expansion}
Let $N=1$. Assume that (A) holds and that the kernels are given by
\eqref{eq:Bessel-mixture-kernel} for a common $\frac12<s<1$.
Fix $\mu_0\in\mathcal I$ such that $\widehat\lambda(\mu_0)>s^*$.
Suppose that this eigenvalue is algebraically simple.
Then there is an interval $U$ about $\mu_0$ such that, with the
coefficients in \eqref{eq:Bessel-expansion-coefficients}, as $L\downarrow0$,
\begin{equation}\label{eq:Bessel-four-term-eigenvalue}
 \lambda_L^*(\mu)=\widehat\lambda(\mu)
 +b_s(\mu)L^{2s}+c_s(\mu)L^{2s+1}
 +q_s(\mu)L^{4s}+r_s(\mu)L^{2s+2}
 +O_K(L^{4s+1}),
\end{equation}
uniformly on every compact subinterval $K\Subset U$, also after one
or two derivatives in $\mu$. Since $T_{-\mu}=T_\mu$, the same
statement holds near $-\mu_0$. On these reflected neighborhoods,
$\widehat\lambda$, $b_s$, $q_s$ and $r_s$ are unchanged under
$\mu\mapsto-\mu$, whereas $c_s$ changes sign.
\end{theorem}

The first possible correction has fractional order $2s\in(1,2)$.
If $\mu_0=0$, $m=1$, and $a(x)=A(x)$ is nonconstant, then
$b_s(0)>0$; see \eqref{eq:Bessel-scalar-positive}.

At zero weight, the expansion also determines the minimizing weight
and the separation of the two whole-space spectral values.
Under the hypotheses of
Theorem~\ref{thm:Bessel-fractional-system} at $\mu_0=0$, define
\begin{equation}\label{eq:Bessel-splitting-coefficients}
 \mathfrak d_s=s\big\langle\psi_0,
 \operatorname{diag}(\sigma_i^2)P\phi_0\big\rangle>0,
 \qquad
 \mathfrak c_s=\langle\psi_0,D_sU_s\phi_0\rangle.
\end{equation}

\Needspace{7\baselineskip}
\begin{theorem}\label{thm:Bessel-threshold-splitting}
Assume that the hypotheses of Theorem~\ref{thm:Bessel-fractional-system}
hold at $\mu_0=0$. Then, with the coefficients in
\eqref{eq:Bessel-splitting-coefficients}, the following statements hold:
\begin{enumerate}
\item[(i)] For all sufficiently small $L>0$, $\mu\mapsto\lambda_L^*(\mu)$
has a unique global minimum at a point $\mu_L\in\mathcal I$, with
\begin{equation}\label{eq:Bessel-weight-displacement}
 \mu_L=\frac{s\mathfrak c_s}{\mathfrak d_s}
 L^{2s+1}+o(L^{2s+1}).
\end{equation}
\item[(ii)] The whole-space spectral values satisfy
\begin{equation}\label{eq:Bessel-gap}
 \lambda^*(L)-\lambda'(L)
 =-\frac{s^2\mathfrak c_s^2}{\mathfrak d_s}
 L^{4s+2}+o(L^{4s+2}).
\end{equation}
In particular, the difference $\lambda^*(L)-\lambda'(L)$ in
\eqref{eq:Bessel-gap} can have noninteger leading order
$4s+2\in(4,6)$. If $m=1$, or $A(x)$ is pointwise symmetric,
then $\mathfrak c_s=0$; in that case the formula does not claim a
nonzero leading splitting coefficient.
\end{enumerate}
\end{theorem}

\par\medskip\noindent
If the system satisfies the detailed-balance condition of
Proposition~\ref{prop:detailed-balance-rigidity}, then
$\mathfrak c_s=0$ and the two generalized whole-space thresholds
coincide exactly for every period, not merely to order
$L^{4s+2}$. Thus a nonzero splitting coefficient detects a failure
of this constant diagonal symmetrization condition.
\par\medskip

\begin{example}\label{ex:Bessel-sharp-fractional-gap}
The leading coefficient in \eqref{eq:Bessel-gap} can be
nonzero for systems. Choose $m=2$, $\sigma_1=\sigma_2=1$,
$q=2\pi$, and
\[
 A_\delta(x)=
 \begin{pmatrix}
 \delta\cos(qx)&1+\delta\sin(qx)\\
 1&0
 \end{pmatrix},\qquad 0<|\delta|\ll1.
\]
The isolated averaged eigenvalue at $\delta=0$ is $1$ and
lies strictly above the multiplication spectrum. The shared calculation in
Subsection~\ref{subsec:threshold-gap-examples} gives
\begin{equation}\label{eq:Bessel-example-coefficients}
 \mathfrak c_s=\frac{\delta^2}{36q^{2s+1}}+O(\delta^3),
 \qquad \mathfrak d_s=s+O(\delta^2).
\end{equation}
Consequently the fractional spectral separation is sharp:
\[
 \lambda^*(L)-\lambda'(L)
 =-\left[\frac{s\delta^4}{1296(2\pi)^{4s+2}}
 +O(\delta^5)\right]L^{4s+2}
 +o_\delta(L^{4s+2}).
\]
For each fixed sufficiently small nonzero $\delta$, the coefficient
in brackets is positive.
\end{example}

\subsubsection{Laplace kernels}
\label{subsec:Laplace-period}

We now compute integer-order corrections for Laplace
dispersal. Here $N=1$ and all dispersal kernels are
$J_i(z)=\tfrac12e^{-|z|}$. Their exponential-moment domain is
$\mathcal I=(-1,1)$, and $M_i(\mu)=(1-\mu^2)^{-1}$.
For a continuous periodic scalar function $f$, set
\[
Pf=\int_0^1f(x)\,\mathrm dx,\qquad
R f=(-\partial_{xx})^{-1}(I-P)f,
\]
where $Rf$ is the unique mean-zero periodic solution of
$-(Rf)''=f-Pf$. Both operators act componentwise on vectors.

For the coefficient formulas, let
$T_\mu=\widehat{\mathcal L}_\mu$. At weights where
$\widehat\lambda(\mu)>s^*$, choose strictly positive right and left
eigenfunctions $\phi_\mu,\psi_\mu$, normalized by
$\langle\psi_\mu,\phi_\mu\rangle=1$.
Use the spectral projection $\Pi_\mu$ and reduced resolvent $S_\mu$
defined in \eqref{eq:averaged-projection-resolvent}, now for the
averaged operator for the Laplace kernels. Set
\[
B_3(\mu)=-2\mu R\partial_xR,\qquad
B_4(\mu)=4\mu^2R\partial_xR\partial_xR-(1-\mu^2)R^2.
\]
The expansion coefficients are
\begin{equation}\label{eq:Laplace-expansion-coefficients}
\begin{aligned}
\mathfrak b_2(\mu)&=\langle\psi_\mu,R\phi_\mu\rangle,\qquad
\mathfrak b_3(\mu)=-2\mu\langle\psi_\mu,R\partial_xR\phi_\mu\rangle,\\
\mathfrak b_4(\mu)&=\langle\psi_\mu,
[B_4(\mu)+RS_\mu R]\phi_\mu\rangle.
\end{aligned}
\end{equation}

\begin{theorem}\label{thm:Laplace-second-order}\label{thm:Laplace-fourth-order}
Let $N=1$. Assume that (A) holds and that
$J_i(z)=\tfrac12e^{-|z|}$ for every $i\in\Lambda$.
Fix $\mu_0\in(-1,1)$ such that $\widehat\lambda(\mu_0)>s^*$.
Then there is a neighborhood $U$ of $\mu_0$ on which the averaged principal
eigenvalue is isolated and algebraically simple. With the coefficients
in \eqref{eq:Laplace-expansion-coefficients}, as $L\downarrow0$,
\begin{equation*}
\lambda_L^*(\mu)=\widehat\lambda(\mu)
+L^2\mathfrak b_2(\mu)+L^3\mathfrak b_3(\mu)
+L^4\mathfrak b_4(\mu)+O(L^5),
\end{equation*}
uniformly on compact subsets of $U$. The same remainder estimate
holds after one or two derivatives in $\mu$.
\end{theorem}

If $\mu_0=0$, analyticity in $L^2$ gives the sharper
unweighted expansion
\begin{equation}\label{eq:Laplace-unweighted-fourth}
\lambda'(L)=\lambda_L^*(0)=\lambda_0+L^2\mathfrak b_2(0)
+L^4\mathfrak b_4(0)+O(L^6),
\end{equation}
where $\lambda_0=\widehat\lambda(0)$, $\phi=\phi_0$, $\psi=\psi_0$,
$\Pi=\Pi_0$, $S=S_0$, and
\[
\mathfrak b_4(0)=-\langle\psi,R^2\phi\rangle
+\langle\psi,RSR\phi\rangle.
\]
For $m=1$, self-adjointness allows the choice
$\psi_0=\phi_0/\int_0^1\phi_0^2$, and hence
\[
\mathfrak b_2(0)=
\frac{\displaystyle\int_0^1|(R\phi_0)'|^2\,\mathrm dx}
{\displaystyle\int_0^1\phi_0^2\,\mathrm dx}\geq0,
\]
with strict inequality precisely when $A(x)$ is nonconstant.

At zero weight, the corresponding whole-space thresholds satisfy the
following formulas.
Under the hypotheses of Theorem~\ref{thm:Laplace-second-order}
at $\mu_0=0$, set
\begin{equation}\label{eq:Laplace-splitting-coefficients}
\mathfrak d=\int_0^1\psi^{\mathsf T}P\phi\,\mathrm dx>0,
\qquad
\mathfrak c=\int_0^1\psi^{\mathsf T}R\partial_xR\phi\,\mathrm dx.
\end{equation}

\begin{corollary}\label{cor:Laplace-two-thresholds}
Assume that the hypotheses of Theorem~\ref{thm:Laplace-second-order}
hold at $\mu_0=0$. Then, with the coefficients in
\eqref{eq:Laplace-splitting-coefficients}, the following statements hold:
\begin{enumerate}
\item[(i)] The averaged spectral branch has a unique nondegenerate global
minimum at zero, with $\widehat\lambda''(0)=2\mathfrak d$.
For all sufficiently small $L>0$, the weighted periodic spectral
value has a unique global minimizer $\mu_L\in(-1,1)$, satisfying
\[
\mu_L=\frac{\mathfrak c}{\mathfrak d}L^3+O(L^5).
\]
\item[(ii)] The two whole-space thresholds satisfy the explicit
sixth-order separation formula
\begin{equation}\label{eq:sixth-order-threshold-gap}
\lambda'(L)=\lambda_0+\mathfrak b_2(0)L^2
 +\mathfrak b_4(0)L^4+O(L^6),\qquad
\lambda^*(L)-\lambda'(L)
=-\frac{\mathfrak c^2}{\mathfrak d}L^6+O(L^8).
\end{equation}
In particular, both thresholds have the same second- and fourth-order
coefficients. If $m=1$, or if $A(x)$ is symmetric for every $x$,
then $\mathfrak c=0$. No assertion that $\mathfrak c$ is nonzero
for every nonsymmetric system is intended.
\end{enumerate}
\end{corollary}

The unweighted expansion can be continued to sixth order.
Under the same zero-weight hypotheses and with the notation
above, define
\begin{equation}\label{eq:Laplace-sixth-order-coefficient}
 \mathfrak b_6=\big\langle\psi,\big[
 R^3-R^2SR-RSR^2+RS(R-\mathfrak b_2(0)I)SR
 \big]\phi\big\rangle.
\end{equation}

\begin{theorem}\label{thm:Laplace-sixth-order}
Assume that the hypotheses of Theorem~\ref{thm:Laplace-second-order}
hold at $\mu_0=0$. Then, with the coefficients in
\eqref{eq:Laplace-expansion-coefficients} and
\eqref{eq:Laplace-sixth-order-coefficient},
\begin{equation*}
\lambda_L^*(0)=\lambda_0+\mathfrak b_2(0)L^2
+\mathfrak b_4(0)L^4+\mathfrak b_6L^6+O(L^8)
\qquad (L\downarrow0).
\end{equation*}
The remainder is uniform on coefficient families for which the
spectral separation and the operator norms of the corresponding
Riesz projections and reduced resolvents are uniformly controlled.
Consequently, Corollary~\ref{cor:Laplace-two-thresholds} gives
\begin{align}
\lambda'(L)&=\lambda_0+\mathfrak b_2(0)L^2+\mathfrak b_4(0)L^4
+\mathfrak b_6L^6+O(L^8),\notag\\
\lambda^*(L)&=\lambda_0+\mathfrak b_2(0)L^2+\mathfrak b_4(0)L^4+
 \left(\mathfrak b_6-\frac{\mathfrak c^2}{\mathfrak d}\right)L^6
 +O(L^8).\label{eq:Laplace-whole-space-sixth}
\end{align}
\end{theorem}

\begin{example}\label{ex:sharp-sixth-order-gap}
The sixth-order coefficient in
\eqref{eq:sixth-order-threshold-gap} need not vanish for
non-self-adjoint cooperative systems. Let $m=2$, $q=2\pi$,
and choose
\[
A_\delta(x)=
\begin{pmatrix}
\delta\cos(qx)&1+\delta\sin(qx)\\
1&0
\end{pmatrix},\qquad 0<|\delta|\ll1.
\]
This matrix is pointwise cooperative and irreducible for
$|\delta|<1$. The strict spectral-gap hypothesis holds for
all sufficiently small $|\delta|$, since at $\delta=0$ the
averaged principal eigenvalue equals $1$ while
$\max_xs_{\rm PF}(A_0-I)=0$.
The shared calculation in
Subsection~\ref{subsec:threshold-gap-examples} yields
\begin{equation}\label{eq:Laplace-example-coefficients}
\mathfrak c=\frac{\delta^2}{36q^3}+O(\delta^3),
\qquad \mathfrak d=1+O(\delta^2).
\end{equation}
Consequently the leading gap is strictly negative:
\[
\lambda^*(L)-\lambda'(L)
=-\left(\frac{\delta^4}{1296(2\pi)^6}
+O(\delta^5)\right)L^6+O_\delta(L^8).
\]
The two limits are taken successively: first $L\downarrow0$
with $\delta$ fixed, then $\delta\to0$. This example shows
that the sixth-order splitting is not an artifact of the
perturbation argument.
\end{example}

\begin{example}\label{ex:cosine-second-order}
In the scalar case take
$A(x)=a_0+\delta\cos(2\pi x)$, with $|\delta|$ sufficiently small.
At $\mu=0$, analytic perturbation of the constant Fourier mode gives
\[
\widehat\lambda(0)=a_0+\frac{\delta^2}{2}+O(\delta^3),
\qquad
\phi_0(x)=1+\delta\cos(2\pi x)+O(\delta^2).
\]
Since $R\cos(2\pi x)=(4\pi^2)^{-1}\cos(2\pi x)$,
\[
\mathfrak b_2(0)=\frac{\delta^2}{8\pi^2}+O(\delta^3)>0
\quad(0<|\delta|\ll1).
\]
Thus the second-order coefficient need not vanish for the continuous
Laplace kernel, in contrast to Theorem~\ref{thm:smooth-superalgebraic}.
\end{example}

The new asymptotic coefficients above are obtained under
additional regularity, spectral isolation, or nondegeneracy
assumptions. They are not consequences of (J) and (A) alone. In
particular, Theorem~\ref{thm:smooth-superalgebraic} and
Theorem~\ref{thm:Laplace-second-order} show that a nonzero universal
$L^2$ coefficient is impossible across the entire kernel class (J).

\subsection{Proofs}

We prove the averaging limit first, then the smooth-kernel, Bessel,
and Laplace expansions in the order stated above. We conclude with
the shared calculation for the two threshold-gap examples.

\subsubsection{Averaging limit}

We first characterize the averaged spectral value and then prove the
convolution limit needed for the comparison.

For the averaged operator $\widehat{\mathcal L}_\mu$ in
\eqref{eq:averaged-operator}, let $C$ be either
$\underline A^\varepsilon$ or $\overline A^\varepsilon$, and set
\[
b_C=\max_x s_{\rm PF}(C(x)-I),\qquad
P_\lambda=\int_{\mathbb T^N}(\lambda I-C(x)+I)^{-1}(I+D(\mu))\,dx
\quad(\lambda>b_C).
\]
The positive matrix $P_\lambda$ depends continuously on $\lambda$. Write
$r(P_\lambda)$ for its spectral radius, the largest modulus of its
eigenvalues. Then $r(P_\lambda)\to0$ as $\lambda\to\infty$.
By Lemma~\ref{lem:approx}, the function
$x\mapsto s_{\rm PF}(C(x)-I)$ attains its maximum $b_C$ on a nonempty
open set. Positivity of the Perron projections gives
$r(P_\lambda)\to\infty$ as $\lambda\downarrow b_C$.
Choose $\lambda>b_C$ with $r(P_\lambda)=1$ and a positive vector
$c$ satisfying $P_\lambda c=c$. Then
\[
\boldsymbol\phi(x)=(\lambda I-C(x)+I)^{-1}(I+D(\mu))c,
\qquad \int_{\mathbb T^N}\boldsymbol\phi(x)\,dx=c,
\]
so $\boldsymbol\phi$ is a positive eigenfunction of the averaged operator
with coefficient $C$. Denote the two eigenvalues by
$\widehat\lambda_\varepsilon^-$ and
$\widehat\lambda_\varepsilon^+$. Since
$0\leq\overline A^\varepsilon-\underline A^\varepsilon\leq4\varepsilon I$,
comparison with either eigenfunction gives
$0\leq\widehat\lambda_\varepsilon^+
-\widehat\lambda_\varepsilon^-\leq4\varepsilon$.
Both Collatz--Wielandt values of $\widehat{\mathcal L}_\mu$ lie between
these eigenvalues, by the same test-function comparison as in
\eqref{eq:weighted-CW-sandwich}. Their difference therefore tends to zero,
which justifies the definition of $\widehat\lambda(\mu)$. Thus the same
approximation applies to the strictly positive constant kernels
$M_i(\mu)$ on the unit torus.

For a periodic function $\psi$, define its spatial average by
		\[
		\overline{\psi} = \int_{\mathbb{T}^N} \psi(y)\,\mathrm{d}y.
		\]
		
		\begin{lemma}\label{lem:conv-avg}
			Let $\phi\in C(\mathbb R^N)$ be $1$-periodic in each spatial direction, and let $J\in L^1(\mathbb{R}^N)$ be nonnegative with $\int_{\mathbb{R}^N}J(z)\,\mathrm{d}z=1$. Define $J_L(x)=L^N J(Lx)$. Then
			\[
			\lim_{L\to0} \left\| J_L * \phi - \overline{\phi} \right\|_{L^\infty(\mathbb{R}^N)} = 0.
			\]
		\end{lemma}
		
\begin{proof}[Proof of Lemma~\ref{lem:conv-avg}]
Set $\psi=\phi-\overline\phi$. Since $\psi$ has zero mean, subtracting
the constant term from a trigonometric approximation preserves
arbitrarily small uniform error. Thus, for any $\varepsilon>0$, uniform
trigonometric approximation on the torus provides
$T(x)=\sum_{k\in F}c_ke^{2\pi\mathrm i k\cdot x}$, with a finite
$F\subset\mathbb Z^N\setminus\{0\}$, such that
$\|\psi-T\|_\infty<\varepsilon$.
As $\|J_L\|_{L^1}=1$, we have
\[
\|J_L*\psi\|_\infty
\leq\varepsilon+
\sum_{k\in F}|c_k|
\left|\int_{\mathbb R^N}J(z)e^{-2\pi\mathrm i k\cdot z/L}\,\mathrm dz\right|.
\]
Each integral tends to zero by the Riemann--Lebesgue lemma. Taking
$L\to0$ and then $\varepsilon\to0$ proves the result.
\end{proof}

\begin{proof}[\textbf{Proof of Theorem~\ref{thm:small-period-limits}}]
We first establish convergence of the weighted periodic values, then
pass to the whole-space thresholds.

\smallskip
\noindent Step 1. Periodic spectral limits.
		Fix $\varepsilon>0$ and choose $\boldsymbol\phi\in X^{++}$ satisfying
		\[
		\widehat{\mathcal J}_0[\boldsymbol{\phi}](x) + A(x)\boldsymbol{\phi}(x) - \boldsymbol{\phi}(x) \le (\widehat\lambda(0) + \varepsilon)\boldsymbol{\phi}(x),\quad x\in\mathbb{R}^N.
		\]
		
By Lemma~\ref{lem:conv-avg} and $\min_{i,x}\phi_i(x)>0$, for all
sufficiently small $L>0$,
\[
\left|L^N\int_{\mathbb R^N}J_i(L(x-y))\phi_i(y)\,\mathrm dy
-\overline\phi_i\right|\leq\varepsilon\phi_i(x).
\]
Hence
\begin{align*}
(\mathcal L_L\boldsymbol\phi)_i
&=(\widehat{\mathcal L}_0\boldsymbol\phi)_i
  +(\mathcal J_L*\boldsymbol\phi)_i-\overline\phi_i\\
&\leq(\widehat\lambda(0)+2\varepsilon)\phi_i,
\end{align*}
and therefore $\lambda_L^*(0)\leq\widehat\lambda(0)+2\varepsilon$.
Similarly, choose $\boldsymbol\psi\in X^{++}$ with
\[
\widehat{\mathcal L}_0\boldsymbol\psi
\geq(\widehat\lambda(0)-\varepsilon)\boldsymbol\psi.
\]
The same estimate and Lemma~\ref{L11}(iii) yield
\[
\mathcal L_L\boldsymbol\psi
\geq(\widehat\lambda(0)-2\varepsilon)\boldsymbol\psi,
\qquad
\widehat\lambda(0)-2\varepsilon
\leq\lambda_L^*(0)\leq\widehat\lambda(0)+2\varepsilon.
\]
Letting $L\to0$ and then $\varepsilon\downarrow0$ gives
$\lambda_L^*(0)\to\widehat\lambda(0)$.

For each fixed $\mu\in\mathcal I$, the same argument applies to
$J_i(z)e^{\mu\cdot z}/M_i(\mu)$, retaining the factors $M_i(\mu)$
in the convolution terms, and gives
\begin{equation}\label{eq:small-period-tilted-local-uniform}
\lim_{L\to0}\lambda_L^*(\mu)=\widehat\lambda(\mu),
\qquad \mu\in\mathcal I.
\end{equation}
By Theorem~\ref{TH2.1}(i), these are finite convex functions on
$\operatorname{int}\mathcal I$; their pointwise convergence is therefore
uniform on every compact subset of this interior.

\smallskip
\noindent Step 2. Whole-space spectral limits.
The periodic limit proved above and Theorem~\ref{TH2.1}(iii)
give the limit of $\lambda'(L)$.
It remains to pass to the minimum over exponential weights.
	Let $\widehat\lambda_n$ correspond to $J_{i,n}=J_i\eta_n$.
The two characterizations in \eqref{eq:averaged-generalized-eigenvalue},
applied to fixed positive test functions, show that
$\widehat\lambda_n$ is continuous and
$\widehat\lambda_n(\mu)\uparrow\widehat\lambda(\mu)$ for
$\mu\in\mathcal I$. Evaluating each component at its minimum point,
as in the derivation of \eqref{eq:weighted-moment-lower-bound}, gives
\[
\widehat\lambda_n(\mu)\geq
\max_i\left\{\int_{\mathbb R^N}J_{i,n}(z)e^{\mu\cdot z}\,\mathrm dz
-1+\min_xa_{ii}(x)\right\}.
\]
Thus the same increasing convergence holds with value $+\infty$
outside $\mathcal I$, and $\widehat\lambda$ is lower semicontinuous.
Positivity on a fixed ball gives a coercive exponential lower bound,
whereas $\widehat\lambda(0)<\infty$. Its finite minimum is therefore
attained in $\mathcal I$.
By the uniform coercivity and zero-weight upper bound in the proof of
Theorem~\ref{thm:large-period-limits}, the minimizers $\mu_L$ remain in a common compact set.
Along every sequence $L_k\downarrow0$, we may therefore pass to
a subsequence with $\mu_{L_k}\to\mu_0$. The estimate \eqref{eq:weighted-moment-lower-bound}
also bounds each $M_i(\mu_{L_k})$.
Fatou's lemma yields $M_i(\mu_0)<\infty$ for every $i$, and hence
$\mu_0\in\mathcal I$.
Kernel ordering gives
\[
 \lambda^*(L_k)=\lambda_{L_k}^*(\mu_{L_k})
 \geq\lambda_{L_k,n}^*(\mu_{L_k}).
\]
For fixed $n$, Step~1 and convexity give locally uniform convergence
on $\mathbb R^N$. Thus
\[
 \lambda_{L_k,n}^*(\mu_{L_k})
 \longrightarrow\widehat\lambda_n(\mu_0),
 \qquad\text{hence}\qquad
 \liminf_{k\to\infty}\lambda^*(L_k)
 \geq\widehat\lambda_n(\mu_0).
\]
Monotone convergence of the kernels and of the generalized principal
eigenvalues gives
\(\widehat\lambda_n(\mu_0)\uparrow\widehat\lambda(\mu_0)\).  Hence
\[
 \liminf_{k\to\infty}\lambda^*(L_k)
 \geq\widehat\lambda(\mu_0)\geq\min_\mu\widehat\lambda(\mu).
\]

Conversely, choose a minimizer $\mu_*\in\mathcal I$ of
$\widehat\lambda$. By \eqref{eq:small-period-tilted-local-uniform},
\[
 \limsup_{L\downarrow0}\lambda^*(L)
 \leq\lim_{L\downarrow0}\lambda_L^*(\mu_*)
 =\widehat\lambda(\mu_*)=\min_\mu\widehat\lambda(\mu).
\]
Since the sequence $L_k$ was arbitrary, we obtain
\[
 \lambda^*(L)
 \longrightarrow\min_{\mu\in\mathbb R^N}\widehat\lambda(\mu)
 \qquad(L\downarrow0).\qedhere
\]

\end{proof}

\subsubsection{Smooth kernels}

We estimate the nonzero Fourier modes and apply perturbation theory
to the isolated principal eigenvalue.

\begin{proof}[\textbf{Proof of Theorem~\ref{thm:smooth-superalgebraic}}]
Write $P f=\int_{\mathbb T^N}f$ and let $K_{i,L,\mu}$ denote
the rescaled weighted convolution on the unit torus. Its Fourier
coefficients are
\[
\widehat K_{i,L,\mu}(k)
=\int_{\mathbb R^N}J_i(z)e^{\mu\cdot z}
 e^{-2\pi\mathrm i k\cdot z/L}\,\mathrm dz,
\qquad k\in\mathbb Z^N.
\]
The coefficient at $k=0$ equals $M_i(\mu)$. Repeated integration by
parts, for $k\ne0$, yields
$|\widehat K_{i,L,\mu}(k)|\leq C_{K,r}L^r|k|^{-r}$.
Choose $r>N$. Absolute convergence of the Fourier series for the
periodized kernel implies
\[
\|K_{i,L,\mu}-M_i(\mu)P\|_{C(\mathbb T^N)\to C(\mathbb T^N)}
\leq C_{K,r}L^r\sum_{k\ne0}|k|^{-r}.
\]
For the operators defined in \eqref{eq:rescaled-weighted-operator}
and \eqref{eq:averaged-operator}, this gives
$\mathcal L_{L,\mu}-\widehat{\mathcal L}_\mu=O(L^r)$
in operator norm, uniformly on $K$. The uniform spectral gap and
simple-eigenvalue perturbation theory
\cite[Chapter~IV, Section~3; Chapter~VII, Section~1]{Kato1995} give
\[
\sup_{\mu\in K}|\lambda_L^*(\mu)-\widehat\lambda(\mu)|=O(L^r).
\]
Here the perturbed eigenfunction remains strictly positive, so its
eigenvalue is $\lambda_L^*(\mu)$. Since $r>N$ is arbitrary, the first
assertion follows. Differentiating the Fourier multipliers and the
associated Riesz projections gives the same estimates for the
derivatives asserted in the theorem.
\end{proof}

\subsubsection{Bessel kernels}

Absolute summability of the Fourier multipliers defining $H_s$ and
$U_s$ makes these operators bounded on $C(\mathbb T)$.
The proofs use the following expansion of the weighted convolution.

\begin{lemma}\label{lem:Bessel-operator-symbol}
Let $N=1$, $\frac12<s<1$, and $\sigma_i>0$ for every $i\in\Lambda$.
Assume that the kernels are given by \eqref{eq:Bessel-mixture-kernel}.
Fix a compact subinterval $K$ of
$\mathcal I=(-\sigma_{\max}^{-1},\sigma_{\max}^{-1})$,
where $\sigma_{\max}=\max_i\sigma_i$.
Then the following statements hold:
\begin{enumerate}
\item[(i)] The rescaled weighted periodic convolution satisfies
\begin{equation}\label{eq:Bessel-operator-asymptotic}
 \mathcal J_{L,\mu}
 =\operatorname{diag}(M_i(\mu))P
 +L^{2s}D_sH_s-2s\mu L^{2s+1}D_sU_s
 +\mathcal E_{L,\mu},
\end{equation}
where there exist $L_K>0$ and $C_{K,j}>0$, $j=0,1,2$, such that
\[
 \sup_{\mu\in K}\|\partial_\mu^j\mathcal E_{L,\mu}\|_{C\to C}
 \le C_{K,j}L^{2s+2}\qquad(0<L<L_K).
\]
\item[(ii)] At zero weight, as $L\downarrow0$,
\begin{equation}\label{eq:Bessel-weight-derivative}
 \left.\partial_\mu\mathcal J_{L,\mu}\right|_{\mu=0}
 =-2sL^{2s+1}D_sU_s+O_{C\to C}(L^{2s+3}).
\end{equation}
\end{enumerate}
\end{lemma}

\begin{proof}[Proof of Lemma~\ref{lem:Bessel-operator-symbol}]
The Gaussian-mixture representation \eqref{eq:Bessel-mixture-kernel}
and Fubini's theorem yield
\eqref{eq:Bessel-symbol}. On the nonzero Fourier mode
$e^{2\pi\mathrm ikx}$, with $\xi=2\pi k/L$, the weighted
convolution multiplier equals
\[
 [1-\sigma_i^2(\mu-\mathrm i\xi)^2]^{-s}
 =\sigma_i^{-2s}|\xi|^{-2s}
 \left(1+\frac{2\mathrm i\mu}{\xi}
 +\frac{1-\sigma_i^2\mu^2}{\sigma_i^2\xi^2}\right)^{-s}.
\]
The branch is the one obtained by analytic continuation from the
positive real axis; the expression inside the original power has
strictly positive real part for $\mu\in\mathcal I$.
Uniformly on $K$, Taylor's formula with two $\mu$ derivatives gives
\[
 [1-\sigma_i^2(\mu-\mathrm i\xi)^2]^{-s}
 =\sigma_i^{-2s}|\xi|^{-2s}
 \left(1-\frac{2s\mathrm i\mu}{\xi}
 +O_{K,j}(|\xi|^{-2})\right).
\]
The $k=0$ multiplier is precisely $M_i(\mu)$. For $k\ne0$, the
error is bounded by $C_{K,j}L^{2s+2}|2\pi k|^{-2s-2}$.
The sum over $k\ne0$ converges. The resulting absolutely convergent
Fourier series is the kernel of an operator on $C(\mathbb T)$ whose
norm is at most the sum of the absolute values of its coefficients.
This proves \eqref{eq:Bessel-operator-asymptotic} and its derivative
estimates. At $\mu=0$, exact differentiation of the multiplier gives
\[
 -2s\mathrm i\sigma_i^2\xi
 (1+\sigma_i^2\xi^2)^{-s-1}
 =-2s\sigma_i^{-2s}\mathrm i\xi|\xi|^{-2s-2}
 +O(|\xi|^{-2s-3}),
\]
and absolute summation proves \eqref{eq:Bessel-weight-derivative}.
\end{proof}

\begin{proof}[\textbf{Proof of Theorem~\ref{thm:Bessel-fractional-system}}]
The strict inequality separates $\widehat\lambda(\mu_0)$
from the essential spectrum of multiplication by $A-I$. Under the
assumed simplicity, an isolating circle and normalized left and right
eigenfunctions persist for $\mu$ in a smaller interval $U$.
The use of a density for the left eigenvector is justified as follows.
If $\boldsymbol\nu_\mu$ is its positive eigenmeasure on the dual of
$C(\mathbb T,\mathbb R^m)$, the adjoint equation gives
\[
 \mathrm d\boldsymbol\nu_\mu(x)=\psi_\mu(x)\,\mathrm dx,
 \qquad
 \psi_\mu(x)=\bigl[(\widehat\lambda(\mu)+1)I-A(x)^{\mathsf T}\bigr]^{-1}
 \operatorname{diag}(M_i(\mu))\boldsymbol\nu_\mu(\mathbb T).
\]
The strict spectral gap makes the inverse matrix continuous and
entrywise positive, so $\psi_\mu$ is continuous and strictly positive.

Fix $K\Subset U$. For a nonzero periodic Fourier mode
$e^{2\pi\mathrm i kx}$, put $\xi=2\pi k/L$. The weighted kernel
symbol has the uniform expansion, together with its first two
$\mu$ derivatives,
\[
\begin{aligned}
 &[1-\sigma_i^2(\mu-\mathrm i\xi)^2]^{-s}\\
 &\quad=\sigma_i^{-2s}|\xi|^{-2s}
 \left(1-\frac{2s\mathrm i\mu}{\xi}
 -\frac{s[\sigma_i^{-2}+(2s+1)\mu^2]}{\xi^2}
 +O_K(|\xi|^{-3})\right).
\end{aligned}
\]
This follows by expanding $(1+w)^{-s}$ through $w^2$, with
$w=2\mathrm i\mu/\xi+(\sigma_i^{-2}-\mu^2)/\xi^2$.
The $k=0$ multiplier equals $M_i(\mu)$ exactly. Since
$2s>1$, absolute summation of the nonzero Fourier coefficients
and their differentiated remainders establishes on
$C(\mathbb T,\mathbb C^m)$
\[
 \mathcal L_{L,\mu}=T_\mu+L^{2s}D_sH_s
 -2s\mu L^{2s+1}D_sU_s+L^{2s+2}Z_{s,\mu}
 +O_{C\to C,K}(L^{2s+3}).
\]
The isolating contour for $\widehat\lambda(\mu)$ can be selected
uniformly on $K$. The standard simple-eigenvalue resolvent expansion
with $S_\mu$ defined by \eqref{eq:averaged-projection-resolvent} is
\[
 \lambda_L^*(\mu)=\widehat\lambda(\mu)
 +\langle\psi_\mu,\Delta_{L,\mu}\phi_\mu\rangle
 +\langle\psi_\mu,\Delta_{L,\mu}S_\mu
 \Delta_{L,\mu}\phi_\mu\rangle
 +O_K(\|\Delta_{L,\mu}\|^3),
\]
where $\Delta_{L,\mu}=\mathcal L_{L,\mu}-T_\mu$.
The first pairing supplies the coefficients of orders
$L^{2s}$, $L^{2s+1}$ and $L^{2s+2}$, and the second pairing
supplies $q_sL^{4s}$. All omitted terms are
$O_K(L^{4s+1})$ because $6s>4s+1$ and
$2s+3>4s+1$ for $1/2<s<1$.
For the derivative remainder, choose a complex
neighborhood of $K$ whose real parts lie in the common moment
interval. The Fourier multipliers, $T_\mu$, and the Riesz
projection are holomorphic there. Compactness of $K$ and the
uniform isolating spectral gap allow a smaller neighborhood on
which the resolvents are uniformly bounded. The displayed
operator-norm remainder, estimated uniformly on that neighborhood,
remains $O_K(L^{4s+1})$ after one and two real $\mu$
derivatives by the Cauchy integral formula. This justifies the
$C^2$ estimate used later to differentiate the minimizing
spectral weight.
The strictly positive eigenvector identifies the
perturbed branch with $\lambda_L^*(\mu)$.
Finally, the averaged operator satisfies $T_{-\mu}=T_\mu$.
Thus the eigendata and the reduced resolvent can be chosen even in
$\mu$; inspection of the coefficients proves the asserted parity.

Since $1/2<s<1$, the powers in
\eqref{eq:Bessel-four-term-eigenvalue} are strictly increasing. In particular,
\begin{equation*}
 \lambda_L^*(\mu)
 =\widehat\lambda(\mu)+L^{2s}b_s(\mu)
 +L^{2s+1}c_s(\mu)+O(L^{4s}),
\end{equation*}
with the same local uniformity and derivative estimates. The first
possible correction is therefore of fractional order $2s\in(1,2)$.
If $\mu_0=0$, $m=1$, and $a(x)=A(x)$ is nonconstant, then
\begin{equation}\label{eq:Bessel-scalar-positive}
 b_s(0)=\sigma_1^{-2s}\frac{\displaystyle
 \sum_{k\ne0}|2\pi k|^{-2s}|\widehat\phi_0(k)|^2}
 {\displaystyle\int_0^1\phi_0(x)^2\,\mathrm dx}>0.
\end{equation}

If $m=1$ and $\mu=0$, $T_0$ is self-adjoint, so
$\psi_0=\phi_0/\|\phi_0\|_2^2$. Parseval's identity gives
\eqref{eq:Bessel-scalar-positive}. Its numerator vanishes only
when $\phi_0$ is constant; the equation
$T_0\phi_0=\widehat\lambda(0)\phi_0$, with $T_0$ defined by
\eqref{eq:averaged-operator}, then forces $a(x)$ to be constant.
\end{proof}

\begin{proof}[\textbf{Proof of Theorem~\ref{thm:Bessel-threshold-splitting}}]
For these even kernels,
$M_i(\mu)=(1-\sigma_i^2\mu^2)^{-s}$ is strictly greater
than $M_i(0)=1$ when $\mu\ne0$. For the positive eigenfunction
$\phi_0$ of $T_0$, the residual
$T_\mu\phi_0-\widehat\lambda(0)\phi_0
=\operatorname{diag}(M_i(\mu)-1)P\phi_0$ is strictly positive
for every nonzero $\mu\in\mathcal I$. Dividing componentwise by
$\phi_0$ and using \eqref{eq:averaged-generalized-eigenvalue} gives
$\widehat\lambda(\mu)>\widehat\lambda(0)$.
Thus $\widehat\lambda$ has a strict global minimum at zero.
In a neighborhood of zero, its simplicity
gives
\[
 \widehat\lambda'(0)=0,
 \qquad \widehat\lambda''(0)
 =2s\big\langle\psi_0,
 \operatorname{diag}(\sigma_i^2)P\phi_0\big\rangle
 =2\mathfrak d_s>0.
\]
At each endpoint of $\mathcal I$, at least one $M_i(\mu)$
diverges. Comparison with the scalar diagonal subsystem gives
\[
 \lambda_L^*(\mu)\ge
 M_i(\mu)+\min_x(a_{ii}(x)-1),
\]
for each $i$. This comparison is uniform in $L$ and confines every
global minimizer to a compact subset of $\mathcal I$.
The convergence in \eqref{eq:small-period-tilted-local-uniform}
is locally uniform there, so all
minimizers converge to zero. The $C^2$ version of
Theorem~\ref{thm:Bessel-fractional-system} gives strict convexity
in a fixed neighborhood of zero for small $L$, proving uniqueness.

At $\mu=0$ the exact differentiated multiplier identity
\eqref{eq:Bessel-weight-derivative} has leading order $L^{2s+1}$.
The perturbations of the right eigenfunction and the left
eigenfunctional are $O(L^{2s})$, so their contributions to the first
derivative are $O(L^{4s+1})$. The simple-eigenvalue differentiation
formula therefore gives
\[
 \partial_\mu\lambda_L^*(0)
 =-2s\mathfrak c_sL^{2s+1}
 +O(L^{4s+1}+L^{2s+3}),
 \qquad
 \partial_{\mu\mu}\lambda_L^*(0)=2\mathfrak d_s+o(1).
\]
Taylor expansion of
$\partial_\mu\lambda_L^*(\mu_L)=0$ yields
\eqref{eq:Bessel-weight-displacement}. A second Taylor expansion
gives
\[
 \lambda_L^*(\mu_L)-\lambda_L^*(0)
 =\partial_\mu\lambda_L^*(0)\mu_L
 +\mathfrak d_s\mu_L^2+o(\mu_L^2)
 =-\frac{s^2\mathfrak c_s^2}{\mathfrak d_s}
 L^{4s+2}+o(L^{4s+2}).
\]
The identities for $\lambda^*$ and $\lambda'$ then follow from
Theorem~\ref{TH2.1}. If $m=1$, or if $A$ is pointwise symmetric,
$T_0$ is self-adjoint and $D_sU_s$ is skew-adjoint on real
$L^2$; therefore $\mathfrak c_s=0$.
\end{proof}

\par\medskip\noindent
The range $\frac12<s<1$ differs from the Laplace case in
Section~\ref{subsec:Laplace-period}. At $s=1$ and
$\sigma_i=1$, the corresponding Fourier-multiplier expansion gives the
$L^2$ and $L^3$ terms of
Lemma~\ref{lem:Laplace-operator-expansion}, and
\eqref{eq:Bessel-gap} agrees with the sixth-order formula
\eqref{eq:sixth-order-threshold-gap}. Neither a common universal
$L^2$ correction nor a universal sixth-order splitting is valid
throughout the class of admissible kernels.
\par\medskip

\subsubsection{Laplace kernels}

For the proofs, we use the following operator expansion.

\begin{lemma}\label{lem:Laplace-operator-expansion}
Let $N=1$ and $J_i(z)=\tfrac12e^{-|z|}$ for every $i\in\Lambda$.
Fix a compact set $K\Subset(-1,1)$, and let $P$ and $R$ be as defined
in Subsection~\ref{subsec:Laplace-period}.
Then the following statements hold:
\begin{enumerate}
\item[(i)] The period-one weighted convolution admits the operator-norm expansion
\begin{align}
\mathcal J_{L,\mu}={}&\frac{1}{1-\mu^2}P+L^2R
-2\mu L^3 R\partial_xR+L^4\bigl[4\mu^2R\partial_xR\partial_xR
 -(1-\mu^2)R^2\bigr]+\mathcal E_{L,\mu},
\label{eq:Laplace-operator-fourth}
\end{align}
where there exist $L_K>0$ and $C_{K,j}>0$, $j=0,1,2$, such that
\[
\sup_{\mu\in K}\|\partial_\mu^j\mathcal E_{L,\mu}\|_{C\to C}
\leq C_{K,j}L^5 \qquad(0<L<L_K).
\]
For every fixed integer $j\geq0$, expansion
\eqref{eq:Laplace-operator-fourth} remains valid after $j$ derivatives
in $\mu$, with an $O_{C\to C}(L^5)$ remainder uniform on $K$.
All compositions in \eqref{eq:Laplace-operator-fourth} are bounded
operators on $C(\mathbb T)$; in particular, the derivative $\partial_x$
always acts on the range of $R$.
\item[(ii)] At zero weight, as $L\downarrow0$,
\begin{equation}\label{eq:Laplace-zero-weight-identity}
\mathcal J_{L,0}=P+L^2R(I+L^2R)^{-1}
=P+L^2R-L^4R^2+O_{C\to C}(L^6).
\end{equation}
\end{enumerate}
\end{lemma}

\begin{proof}[Proof of Lemma~\ref{lem:Laplace-operator-expansion}]
The Fourier multiplier on $e^{2\pi\mathrm i kx}$ is
\[
\int_{\mathbb R}\frac{e^{-|z|}}2e^{\mu z}
 e^{-2\pi\mathrm i kz/L}\,\mathrm dz
=\frac{L^2}{(2\pi k)^2+2\mathrm i\mu L(2\pi k)
 +(1-\mu^2)L^2}.
\]
Its value at $k=0$ is $(1-\mu^2)^{-1}$. On mean-zero
functions write $D=-\partial_{xx}$ and
$B_{L,\mu}=2\mu L\partial_x+(1-\mu^2)L^2I$.
The multiplier formula is equivalently
\[
\mathcal J_{L,\mu}(I-P)
=L^2(D+B_{L,\mu})^{-1}(I-P)
=L^2R\bigl[I+B_{L,\mu}R\bigr]^{-1}(I-P).
\]
The last factorization is an identity of bounded operators on
$C(\mathbb T)$: $R$ maps continuous mean-zero functions to periodic
$C^2$ functions, so $\partial_xR$ is bounded. For $\mu$ in a
fixed compact subset of $(-1,1)$,
$\|B_{L,\mu}R\|_{C\to C}\leq C_KL$.
Expanding $(I+B_{L,\mu}R)^{-1}$ by a Neumann series through
order $L^2$ yields \eqref{eq:Laplace-operator-fourth}; the
remaining geometric series is $O_K(L^3)$ before multiplication
by $L^2$. Differentiating the convergent series in $\mu$
proves the stated uniform derivative bounds. At $\mu=0$,
$B_{L,0}=L^2I$, and the exact formula follows.
\end{proof}

\begin{proof}[\textbf{Proof of Theorem~\ref{thm:Laplace-second-order}}]
The strict spectral inequality places the principal eigenvalue
above the essential spectrum of multiplication by $A-I$.
There the eigenvalue equation is equivalent to the strongly positive
finite-rank problem
$\phi=[(\lambda+1)I-A]^{-1}M(\mu)P\phi$.
Positivity and irreducibility give algebraic simplicity, and the
finite-rank perturbation makes the eigenvalue isolated.
The right and left eigenfunctions are continuous and strictly positive.

Let $K\Subset U$. Lemma~\ref{lem:Laplace-operator-expansion}
gives the operator expansion, uniformly for $\mu\in K$,
\[
\mathcal L_{L,\mu}=T_\mu+L^2R+L^3B_3(\mu)
+L^4B_4(\mu)+O_{C\to C}(L^5).
\]
Its remainder and its first two $\mu$-derivatives obey the same
estimate. Since $\widehat\lambda(\mu)$ is simple and separated from the
remaining spectrum, Riesz projection perturbation supplies a unique
normalized eigenfunction $\phi_{L,\mu}$ and eigenvalue
$\lambda_{L,\mu}$, analytic in the perturbation parameters and
locally uniform in $\mu$. Write
\[
\phi_{L,\mu}=\phi_\mu+L^2\phi_{2,\mu}+O(L^3),
\qquad\langle\psi_\mu,\phi_{L,\mu}\rangle=1.
\]
At order $L^2$ we obtain
\begin{equation}\label{eq:Laplace-second-order-corrector}
(\widehat\lambda(\mu)-T_\mu)\phi_{2,\mu}
=(R-\mathfrak b_2(\mu)I)\phi_\mu,
\qquad \langle\psi_\mu,\phi_{2,\mu}\rangle=0.
\end{equation}
Applying the reduced resolvent \eqref{eq:averaged-projection-resolvent}
to \eqref{eq:Laplace-second-order-corrector} gives
$\phi_{2,\mu}=S_\mu R\phi_\mu$. Pairing the eigenvalue
equation with $\psi_\mu$ and using
$\langle\psi_\mu,\phi_{L,\mu}\rangle=1$ produces
\[
\lambda_{L,\mu}-\widehat\lambda(\mu)
=L^2\langle\psi_\mu,R\phi_{L,\mu}\rangle
+L^3\langle\psi_\mu,B_3(\mu)\phi_{L,\mu}\rangle
+L^4\langle\psi_\mu,B_4(\mu)\phi_{L,\mu}\rangle
+O(L^5),
\]
which proves the three coefficient formulas and the stated
remainder. The differentiation claim follows by differentiating the
Riesz projection formula in $\mu$; the spectral circles are uniform
on $K$. The perturbed eigenfunction remains strictly positive, so
$\lambda_{L,\mu}=\lambda_L^*(\mu)$.
At $\mu=0$, identity \eqref{eq:Laplace-zero-weight-identity}
shows that the family is analytic in $t=L^2$. Consequently only
even powers of $L$ occur, and analytic perturbation yields the
$O(L^6)$ remainder and the stated expression for
$\mathfrak b_4(0)$.

In the scalar unweighted case, $T_0$ is self-adjoint and
$\psi_0=\phi_0/\|\phi_0\|_2^2$. Integration by parts gives
$\langle\phi_0,R\phi_0\rangle=\int_0^1|(R\phi_0)'|^2$.
This vanishes exactly when $\phi_0$ is constant. The eigenvalue
equation then forces $A$ to be constant, and the converse is immediate.
\end{proof}

\begin{proof}[Proof of Corollary~\ref{cor:Laplace-two-thresholds}]
For the common Laplace kernel the averaged operator is
$T_\mu=A-I+M(\mu)P$, where
$M(\mu)=(1-\mu^2)^{-1}$.
For $\mu\ne0$, the positive eigenfunction $\phi$ of $T_0$ satisfies
\[
T_\mu\phi-\widehat\lambda(0)\phi=(M(\mu)-1)P\phi\gg0.
\]
The Collatz--Wielandt characterization yields
$\widehat\lambda(\mu)>\widehat\lambda(0)$.
Simple-eigenvalue perturbation at
$\mu=0$ gives
\[
\widehat\lambda'(0)=0,\qquad
\widehat\lambda''(0)=2\langle\psi,P\phi\rangle
=2\mathfrak d>0.
\]
The constant-vector test furnishes a lower bound
$\lambda_L^*(\mu)\geq M(\mu)-C$ with $C$ independent of
$L$ and $\mu$. Thus global minimizers exist and remain in a
fixed compact subset of $(-1,1)$ for small $L$.
The locally uniform convergence in
\eqref{eq:small-period-tilted-local-uniform} and strict minimality of
$\widehat\lambda$ force all these minimizers to tend to zero.
Lemma~\ref{lem:Laplace-operator-expansion} and simple-eigenvalue
perturbation give an analytic extension of the branch near $(L,\mu)=(0,0)$.
For some $\delta>0$,
\[
\partial_{\mu\mu}\lambda_L^*(\mu)
\longrightarrow\widehat\lambda''(\mu)>0
\quad\text{uniformly for }|\mu|\leq\delta.
\]
Thus the minimizing weight is unique for small $L$.

At $\mu=0$ the derivative of the weighted convolution is exactly
\[
\left.\partial_\mu\mathcal J_{L,\mu}\right|_{\mu=0}
=-2L^3R(I+L^2R)^{-1}\partial_x
 R(I+L^2R)^{-1}
=-2L^3R\partial_xR+O_{C\to C}(L^5).
\]
The principal left and right eigenfunctions of
$\mathcal L_{L,0}$ converge to $\psi$ and $\phi$ in the
respective dual and primal norms, with errors $O(L^2)$.
The simple-eigenvalue differentiation formula therefore yields
\[
\left.\partial_\mu\lambda_L^*(\mu)\right|_{\mu=0}
=-2\mathfrak c L^3+O(L^5),\qquad
\left.\partial_{\mu\mu}\lambda_L^*(\mu)\right|_{\mu=0}
=2\mathfrak d+O(L^2).
\]
Taylor's theorem applied to
$\partial_\mu\lambda_L^*(\mu_L)=0$ then gives
$\mu_L=(\mathfrak c/\mathfrak d)L^3+O(L^5)$.
A second Taylor expansion about zero, using the uniform analyticity
of the eigenvalue branch, shows that
\[
\lambda_L^*(\mu_L)-\lambda_L^*(0)
=(-2\mathfrak cL^3)\frac{\mathfrak c}{\mathfrak d}L^3
+\mathfrak d\left(\frac{\mathfrak c}{\mathfrak d}L^3\right)^2
+O(L^8)
=-\frac{\mathfrak c^2}{\mathfrak d}L^6+O(L^8).
\]
The fourth-order expansion of $\lambda_L^*(0)$ is
\eqref{eq:Laplace-unweighted-fourth}, and
Theorem~\ref{TH2.1} identifies the two whole-space thresholds.
Finally, in the scalar or pointwise symmetric case $T_0$ is
self-adjoint. Taking $\psi=\phi/\langle\phi,\phi\rangle$ and
using self-adjointness of $R$ and periodic integration by parts gives
$\langle\psi,R\partial_xR\phi\rangle
=\langle R\phi,\partial_xR\phi\rangle/
\langle\phi,\phi\rangle=0$.
\end{proof}

\begin{proof}[\textbf{Proof of Theorem~\ref{thm:Laplace-sixth-order}}]
At zero exponential weight the exact periodic convolution
identity in Lemma~\ref{lem:Laplace-operator-expansion} yields, with
$t=L^2$,
\[
\mathcal L_{L,0}=T_0+tR-t^2R^2+t^3R^3+O_{C\to C}(t^4).
\]
The equality holds in the operator norm of
$C(\mathbb T,\mathbb R^m)$, because $R$ is bounded and
$(I+tR)^{-1}$ is analytic for $|t|$ small.
The isolated simple eigenvalue of $T_0$ therefore admits an analytic
eigenvalue $\lambda(t)$ and an analytic right eigenfunction
$\phi(t)$, normalized by
$\langle\psi,\phi(t)\rangle=1$. Write
\[
\begin{aligned}
\lambda(t)&=\lambda_0+t\lambda_1+t^2\lambda_2
 +t^3\lambda_3+O(t^4),\\
\phi(t)&=\phi+t\phi_1+t^2\phi_2+O_{C}(t^3),
\qquad\langle\psi,\phi_j\rangle=0\quad(j=1,2).
\end{aligned}
\]
Equating first-order terms and applying $S=S_0$ from
\eqref{eq:averaged-projection-resolvent} gives
\[
\lambda_1=\langle\psi,R\phi\rangle=\mathfrak b_2(0),
\qquad \phi_1=SR\phi.
\]
At order $t^2$ one obtains
\[
\lambda_2=\langle\psi,(-R^2+RSR)\phi\rangle=\mathfrak b_4(0),
\qquad
\phi_2=S\bigl[(R-\mathfrak b_2(0)I)SR\phi-R^2\phi\bigr].
\]
Pairing the order-$t^3$ terms of
$\mathcal L_{L,0}\phi(t)=\lambda(t)\phi(t)$ with $\psi$ gives
\[
\lambda_3=\langle\psi,R^3\phi-R^2\phi_1+R\phi_2\rangle
=\big\langle\psi,\big[R^3-R^2SR-RSR^2
 +RS(R-\mathfrak b_2(0)I)SR\big]\phi\big\rangle=\mathfrak b_6.
\]
Analytic perturbation gives the $O(t^4)=O(L^8)$ remainder,
and positivity of the perturbed eigenfunction identifies this
branch with the principal generalized eigenvalue.
The first identity in \eqref{eq:Laplace-whole-space-sixth}
follows from Theorem~\ref{TH2.1}. The second follows by adding
\eqref{eq:sixth-order-threshold-gap} to the first.
All operators in the calculation are bounded on the indicated
periodic Banach space. For the coefficient families in the theorem,
the uniform bounds on $\Pi$ and $S$ give spectral circles centered
at the respective $\lambda_0$ with a common small radius, by the
resolvent Neumann series on the complementary space.
Since $R$ is fixed, the perturbation is analytic on a common
$t$-disk. Shrinking this disk uniformly ensures
$\|\Pi(t)-\Pi\|<1/2$. For this positivity argument, normalize
$\|\phi\|_C=1$. For real $t\geq0$ in the disk, choose
$\Pi(t)\phi$ as an eigenvector; its norm
is at least $1/2$. For each coefficient matrix, the set of such $t$
for which this vector is strictly positive is nonempty and open.
It is also closed: a nonnegative nonzero limiting eigenvector is
strictly positive by cooperation, pointwise irreducibility, and
strict positivity of the Laplace convolution kernel (and, at $t=0$,
by the positivity of $P$ on nonnegative nonzero components).
Consequently the analytic branch is the principal eigenvalue
throughout the common real interval. The Cauchy estimates now give
a uniform $O(t^4)$ remainder, without requiring a uniform positive
lower bound on the components of $\phi$.
\end{proof}

\subsubsection{Example calculations}\label{subsec:threshold-gap-examples}

The two threshold-gap examples use the same averaged eigenfunctions.
We compute their coefficients together on the unit torus.

\begin{proof}[Proof of the formulas in Examples~\ref{ex:Bessel-sharp-fractional-gap}
and~\ref{ex:sharp-sixth-order-gap}]
Both examples use the matrix $A_\delta$ displayed in
Example~\ref{ex:Bessel-sharp-fractional-gap}, with $q=2\pi$.
At zero weight the averaged operator is $A_\delta-I+P$ for both
kernel families. Thus the right and left averaged eigenfunctions
are the same.
Writing $v=2^{-1/2}(1,1)^{\mathsf T}$, elementary perturbation
of the isolated simple averaged eigenvalue at $\delta=0$ gives
\[
\begin{aligned}
\phi(x)&=v+\frac{\delta}{3\sqrt2}
\begin{pmatrix}2\\1\end{pmatrix}
(\cos(qx)+\sin(qx))+O(\delta^2),\\
\psi(x)&=v+\frac{\delta}{3\sqrt2}
\left[
\begin{pmatrix}2\\1\end{pmatrix}\cos(qx)
+\begin{pmatrix}1\\2\end{pmatrix}\sin(qx)
\right]+O(\delta^2).
\end{aligned}
\]
Here the remainders hold in $C(\mathbb T,\mathbb R^2)$,
with the normalization $\langle\psi,\phi\rangle=1$.
On the zero-mean Fourier modes the unperturbed
reduced eigenvalue equation is inverted by
\[
\bigl[2I-A_0\bigr]^{-1}
=\frac13\begin{pmatrix}2&1\\1&2\end{pmatrix},
\qquad A_0=\begin{pmatrix}0&1\\1&0\end{pmatrix}.
\]

For the Bessel kernels, $D_s=I$ and
\[
U_s\cos(qx)=-q^{-2s-1}\sin(qx),\qquad
U_s\sin(qx)=q^{-2s-1}\cos(qx).
\]
Since $U_s$ annihilates constants and has zero-mean range, terms
pairing a constant with either eigenfunction remainder vanish.
The functions $\cos(qx)$ and $\sin(qx)$ are orthogonal in
$L^2(0,1)$, both with squared norm $1/2$. Thus the expansions give
\[
\langle\psi,U_s\phi\rangle
=\frac{\delta^2}{36q^{2s+1}}\bigl[(4+1)-(2+2)\bigr]
+O(\delta^3).
\]
The first-order corrections have zero mean, so
$\langle\psi,P\phi\rangle=1+O(\delta^2)$.
The definitions in \eqref{eq:Bessel-splitting-coefficients} now give
\eqref{eq:Bessel-example-coefficients}; substitution into
\eqref{eq:Bessel-gap} gives the separation asserted in
Example~\ref{ex:Bessel-sharp-fractional-gap}.

For the Laplace kernel, $R\cos(qx)=q^{-2}\cos(qx)$ and
$R\sin(qx)=q^{-2}\sin(qx)$, whence
\[
R\partial_xR\cos(qx)=-q^{-3}\sin(qx),\qquad
R\partial_xR\sin(qx)=q^{-3}\cos(qx).
\]
The same integration gives
\eqref{eq:Laplace-example-coefficients}, using
\eqref{eq:Laplace-splitting-coefficients}.
Substitution into \eqref{eq:sixth-order-threshold-gap} proves the
formula in Example~\ref{ex:sharp-sixth-order-gap}.
In both examples the leading separation coefficient is nonzero
for every sufficiently small fixed $\delta\ne0$.
\end{proof}

\section{Noncooperative KPP dynamics}
\label{sec:nonlinear-dynamics}

In this section, we apply the spectral results to
\eqref{eq:general-nonlocal-KPP} through auxiliary cooperative systems.
We return to the original spatial variables and the fixed period lattice
$\mathbb L$, except when the dependence on a common period $L$ is stated.
For bounded continuous nonnegative initial data $\boldsymbol u_0$, denote
the solution of \eqref{eq:general-nonlocal-KPP} by
$\boldsymbol u(t,\cdot;\boldsymbol u_0)$.
Throughout this section, assumptions (J), (A), and (Q) from the
introduction hold.

\subsection{Main results}

We first state extinction and persistence criteria, then characterize
spreading and its dependence on the spatial period.

\subsubsection{Extinction and persistence}

The whole-space spectral values in Theorem~\ref{TH2.1} determine the
following extinction and persistence criteria.

\Needspace{6\baselineskip}
\begin{theorem}
\label{thm:nonlocal-threshold-dynamics}
Assume that (J), (A), and (Q) hold. Then the following statements hold:
\begin{enumerate}
\item[(i)] If $\lambda'<0$, then, for every
$\boldsymbol u_0\in C_{\rm b}(\mathbb R^N,[0,\infty)^m)$,
\[
\lim_{t\to\infty}
\|\boldsymbol u(t,\cdot;\boldsymbol u_0)\|_\infty=0.
\]
Consequently, $\boldsymbol0$ is the only nonnegative entire solution
bounded on $\mathbb R\times\mathbb R^N$, where an entire solution
is defined for every $t\in\mathbb R$.

\item[(ii)] If $\lambda'>0$, then there exists $\eta_{\rm per}>0$
such that, for every $\boldsymbol u_0\in X^{++}$,
\[
\liminf_{t\to\infty}
\min_{1\leq i\leq m}\min_{x\in\mathbb T^N}u_i(t,x)
\geq\eta_{\rm per}.
\]

\item[(iii)] If $\lambda^*>0$, then there exists $\eta_{\rm loc}>0$
such that, for every nonzero bounded continuous initial function
$\boldsymbol u_0\geq\boldsymbol0$
and every $R>0$,
\[
\liminf_{t\to\infty}
\min_{1\leq i\leq m}\inf_{|x|\leq R}u_i(t,x)
\geq\eta_{\rm loc}.
\]
\end{enumerate}
\end{theorem}

The next theorem shows how exponential initial tails affect local
persistence and extinction, including the case $\lambda^*<0<\lambda'$.

\begin{theorem}
\label{thm:nonlocal-tail-criteria}
Assume that (J), (A), and (Q) hold. Then the following statements hold:
\begin{enumerate}
\item[(i)] Suppose that $\lambda'>0$. Then there exists $\eta_{\rm tail}>0$
such that, for every bounded continuous initial function
$\boldsymbol u_0\geq\boldsymbol0$,
if there is $\xi\in\mathcal I\setminus\{0\}$ satisfying
the following conditions:
\begin{enumerate}
\item[(a)] $\lambda_L^*(\xi)>0$;
\item[(b)] for some $s>1$ with $s\xi\in\mathcal I$,
$\lambda_L^*(s\xi)<\lambda_L^*(\xi)$;
\item[(c)] for some $B\in\mathbb R$ and $C>0$,
\[
\min_i u_{0,i}(x)\geq C^{-1}e^{-\xi\cdot x}
\quad\hbox{whenever }\xi\cdot x\geq B,
\]
\end{enumerate}
then, for every $R>0$,
\[
\liminf_{t\to\infty}\min_i\inf_{|x|\leq R}u_i(t,x)
\geq\eta_{\rm tail}.
\]

\item[(ii)] Suppose that $\lambda^*<0$. Let
$\boldsymbol u_0\in C_{\rm b}(\mathbb R^N,[0,\infty)^m)$.
Suppose that there are $\xi\in\mathcal I$ and $C>0$ such that
\[
\lambda_L^*(\xi)<0,
\qquad
\max_i u_{0,i}(x)\leq Ce^{-\xi\cdot x}
\quad\hbox{for every }x\in\mathbb R^N.
\]
Then, for every $R>0$,
\[
\lim_{t\to\infty}\max_i\sup_{|x|\leq R}u_i(t,x)=0.
\]
\end{enumerate}
\end{theorem}

\subsubsection{Spreading}

We distinguish planar speeds from speeds along rays through the
Freidlin--G\"artner formula.

For $\lambda^*>0$ and $e\in\mathbb S^{N-1}$, define the extended planar
and ray speeds by
\begin{equation}
\label{eq:planar-and-FG-speeds}
c^*(e)
=\inf_{\substack{\rho>0\\ \rho e\in\mathcal I}}
\frac{\lambda_L^*(\rho e)}{\rho},
\quad \inf\varnothing:=+\infty,\qquad
c_{\rm FG}(e)
=\inf_{\substack{e'\in\mathbb S^{N-1}\\e\cdot e'>0}}
\frac{c^*(e')}{e\cdot e'}\in(0,+\infty].
\end{equation}
The planar speed $c^*(e)$ is finite precisely when
$\rho e\in\mathcal I$ for some $\rho>0$, whereas the ray speed
$c_{\rm FG}(e)$ is finite precisely when
$\mu\cdot e>0$ for some $\mu\in\mathcal I$.
Thus the two speeds need not coincide.

The planar and ray speeds in \eqref{eq:planar-and-FG-speeds} have the
following geometric interpretation.

\begin{proposition}\label{prop:spreading-geometry}
Assume that (J), (A), and (Q) hold and that $\lambda^*>0$.
Define the closed convex set
\[
\mathcal W=\bigcap_{e\in\mathbb S^{N-1}}
\{x\in\mathbb R^N:x\cdot e\leq c^*(e)\},
\]
where a constraint with $c^*(e)=+\infty$ is omitted.
Then, for every $e\in\mathbb S^{N-1}$,
\[
\sup\{r\geq0:re\in\mathcal W\}
=\inf_{\substack{e'\in\mathbb S^{N-1}\\e\cdot e'>0}}
\frac{c^*(e')}{e\cdot e'}=c_{\rm FG}(e).
\]
Thus the planar speed $c^*(e)$ and the radial speed
$c_{\rm FG}(e)$ have different geometric roles.
\end{proposition}

The next theorem identifies $c_{\rm FG}(e)$ as the spreading speed
along the ray in direction $e$.

\begin{theorem}
\label{thm:nonlocal-spreading}
Assume that (J), (A), and (Q) hold. Suppose that $\lambda^*>0$.
Then there exists $\eta_{\rm sp}>0$ such that,
for every nonzero compactly supported
$\boldsymbol u_0\in C_{\rm b}(\mathbb R^N,[0,\infty)^m)$, every
$e\in\mathbb S^{N-1}$, and every $R>0$,  the following statements hold:
\begin{enumerate}
\item[(i)] For every $0<c<c_{\rm FG}(e)$,
\begin{equation*}
\liminf_{t\to\infty}
\min_i\inf_{|x|\leq R}u_i(t,x+cte)\geq\eta_{\rm sp}.
\end{equation*}
\item[(ii)] If $c_{\rm FG}(e)<+\infty$, then, for every $c>c_{\rm FG}(e)$,
\begin{equation}
\lim_{t\to\infty}
\max_i\sup_{|x|\leq R}u_i(t,x+cte)=0.
\label{eq:nonlocal-spreading-outside}
\end{equation}
\end{enumerate}
\end{theorem}

\subsubsection{Period limits}

For the period-dependent statements and their proofs, fix unit-periodic profiles
$A$ and $Q$, and set $A_L(x)=A(x/L)$ and $Q_L(x)=Q(x/L)$.
All spreading speeds below are measured in the original spatial variables.
To indicate the dependence on the period, write $c_L^*(e)=c^*(e)$ and
$c_{{\rm FG},L}(e)=c_{\rm FG}(e)$.

\Needspace{5\baselineskip}
\begin{theorem}
\label{cor:speed-period-limits}
Assume that (J), (A), and (Q) hold, with $L=L_1=\cdots=L_N$.
Fix $e\in\mathbb S^{N-1}$. Then the following statements hold:
\begin{enumerate}
\item[(i)] If $\min_\mu\overline\lambda_\infty(\mu)>0$, then
\begin{gather*}
\lim_{L\to\infty}c_L^*(e)
=c_\infty^*(e):=\inf_{\substack{\rho>0\\\rho e\in\mathcal I}}
\frac{\overline\lambda_\infty(\rho e)}{\rho},\\
\lim_{L\to\infty}c_{{\rm FG},L}(e)
=\inf_{\substack{e'\in\mathbb S^{N-1}\\e\cdot e'>0}}
\frac{c_\infty^*(e')}{e\cdot e'}.
\end{gather*}
\item[(ii)] If $\min_\mu\widehat\lambda(\mu)>0$, then
\begin{gather*}
\lim_{L\to0}c_L^*(e)
=\widehat c_0^*(e):=\inf_{\substack{\rho>0\\\rho e\in\mathcal I}}
\frac{\widehat\lambda(\rho e)}{\rho},\\
\lim_{L\to0}c_{{\rm FG},L}(e)
=\inf_{\substack{e'\in\mathbb S^{N-1}\\e\cdot e'>0}}
\frac{\widehat c_0^*(e')}{e\cdot e'}.
\end{gather*}
\end{enumerate}
The limits are understood in the extended sense.
\end{theorem}

\subsubsection{Small-period corrections}

\label{subsec:quantitative-speeds}

We apply the expansions of Section~\ref{sec:quantitative-period}
to the minimizing weights and directional speeds, first for Bessel
and then for Laplace kernels.

\begin{theorem}\label{thm:Bessel-directional-speed-expansion}
Assume that the Bessel-kernel and spectral hypotheses of
Theorem~\ref{thm:Bessel-fractional-system} hold in neighborhoods of
$\rho_0$ and $-\rho_0$. Suppose that $\widehat\lambda(0)>0$ and that
\[
 F(\rho)=\frac{\widehat\lambda(\rho)}{\rho},
 \qquad 0<\rho<\sigma_{\max}^{-1},
\]
has a unique global minimizer $\rho_0$, with $F''(\rho_0)>0$.
Set $c_0=F(\rho_0)$ and
\[
 B(\rho)=\frac{b_s(\rho)}{\rho},\quad
 C(\rho)=\frac{c_s(\rho)}{\rho},\quad
 Q(\rho)=\frac{q_s(\rho)}{\rho},\quad
 E(\rho)=\frac{r_s(\rho)}{\rho}.
\]
Define the right- and left-moving planar spectral speeds by
\[
 c_L^\pm=\inf_{0<\rho<\sigma_{\max}^{-1}}
 \frac{\lambda_L^*(\pm\rho)}{\rho}.
\]
Then, as $L\downarrow0$,
\begin{equation}\label{eq:Bessel-refined-speeds}
\begin{aligned}
 c_L^\pm={}&c_0+B(\rho_0)L^{2s}
 \pm C(\rho_0)L^{2s+1}\\
 &+\left[Q(\rho_0)-\frac{[B'(\rho_0)]^2}{2F''(\rho_0)}\right]L^{4s}
 +E(\rho_0)L^{2s+2}+O(L^{4s+1}).
\end{aligned}
\end{equation}
Consequently,
\begin{equation}\label{eq:Bessel-refined-speed-gap}
 c_L^+-c_L^-=
 \frac{2c_s(\rho_0)}{\rho_0}L^{2s+1}
 +O(L^{4s+1}).
\end{equation}
If $m=1$, or if $A(x)$ is pointwise symmetric,
$c_s(\rho_0)=0$, and the displayed identity does not assert a
nonzero asymmetry. Whenever the hypotheses of
Theorem~\ref{thm:nonlocal-spreading} hold, these quantities are the
corresponding planar invasion speeds of the nonlinear system.
\end{theorem}

\begin{corollary}\label{cor:Laplace-minimum-speed}
Assume that the hypotheses of Theorem~\ref{thm:Laplace-second-order}
hold. Then the following statements hold as $L\downarrow0$:
\begin{enumerate}
\item[(i)] Suppose that the averaged spectral minimum is uniquely attained
at $\mu_0$, that $\widehat\lambda''(\mu_0)>0$, and that the minimizing
weights of $\lambda_L^*$ remain in a fixed compact subset of $U$ for
small $L$. Then
\[
\mu_L=\mu_0-L^2\frac{\mathfrak b_2'(\mu_0)}
{\widehat\lambda''(\mu_0)}+O(L^3),\qquad
\lambda^*(L)=\widehat\lambda(\mu_0)
+L^2\mathfrak b_2(\mu_0)+O(L^3).
\]
\item[(ii)] Suppose that $\widehat\lambda(0)>0$ and that
$f_0(\rho)=\widehat\lambda(\rho)/\rho$ has a unique global
minimizer $\rho_0\in(0,1)$ with $f_0''(\rho_0)>0$.
Suppose also that the eigenvalue branch is simple and separated from the
remaining spectrum on a neighborhood of $\rho_0$, and that the associated
KPP system satisfies the hypotheses of Theorem~\ref{thm:nonlocal-spreading}.
Then the planar spreading speed
in the positive direction satisfies
\[
c_L^*(1)=f_0(\rho_0)
+L^2\frac{\mathfrak b_2(\rho_0)}{\rho_0}+O(L^3),\qquad
\rho_L=\rho_0-L^2
\frac{\bigl(\mathfrak b_2(\rho)/\rho\bigr)'|_{\rho=\rho_0}}
{f_0''(\rho_0)}+O(L^3).
\]
The same statements hold for the negative direction after replacing
$\mu$ by $-\rho$.
\end{enumerate}
If $\mathfrak b_2$ is only continuously differentiable, the optimizer
remainders are $o(L^2)$ instead of $O(L^3)$.
\end{corollary}

\subsection{Proofs}

We first establish the fixed-period comparison estimates and prove the
nonlinear results, then turn to the period-dependent speeds.

\subsubsection{Auxiliary results}

We first obtain bounds for the nonlinear solution and then construct
cooperative comparison systems with the same linearization.

Set
\[
q_0:=\min_{1\leq i,j\leq m}\min_{x\in\mathbb T^N}q_{ij}(x)>0.
\]

\begin{lemma}
\label{lem:KPP-wellposed-bounded}
Assume that (J), (A), and (Q) hold.
Then there exists $M_*>0$ such that, for every
$\boldsymbol u_0\in C_{\rm b}(\mathbb R^N,[0,\infty)^m)$,
problem \eqref{eq:general-nonlocal-KPP} has a unique global nonnegative mild
solution, that is, a solution of the corresponding time-integrated
equation, satisfying
\[
\sup_{t\geq0}\|\boldsymbol u(t)\|_\infty
\leq M_*+\|\boldsymbol u_0\|_\infty,
\qquad
\limsup_{t\to\infty}\|\boldsymbol u(t)\|_\infty\leq M_*.
\]
If $\boldsymbol u_0\not\equiv\boldsymbol0$, then
$\boldsymbol u(t,x)\gg\boldsymbol0$ for every $t>0$ and
$x\in\mathbb R^N$.
\end{lemma}

\begin{proof}[Proof of Lemma~\ref{lem:KPP-wellposed-bounded}]
The convolution operators are bounded on
$C_{\rm b}(\mathbb R^N)$, and the reaction term is locally Lipschitz in
$\boldsymbol u$. Hence the integral equation associated with
\eqref{eq:general-nonlocal-KPP} has a unique local mild solution. Since the
off-diagonal entries of $A$ and the kernels are nonnegative, the
variation-of-constants formula gives positivity.

Let
\[
a_+=\max_{1\leq i\leq m}\sup_{x\in\mathbb T^N}
\sum_{j=1}^m a_{ij}^+(x).
\]
Here $r^+=r_+=\max\{r,0\}$ denotes the positive part of a real number;
for vectors or vector functions, the positive part is taken componentwise.
At a point where the $i$th component is arbitrarily close to the global
componentwise maximum $M(t):=\max_i\sup_x u_i(t,x)$, the probability-kernel property gives
$J_i*u_i-u_i\leq o(1)$, while
\[
(A(x)\boldsymbol u)_i
-u_i(Q(x)\boldsymbol u)_i
\leq a_+M(t)-q_0M(t)^2+o(1).
\]
The upper right Dini derivative
$D^+M(t):=\limsup_{h\downarrow0}(M(t+h)-M(t))/h$ therefore satisfies
$D^+M\leq a_+M-q_0M^2$. Comparison with the scalar logistic equation gives the stated bounds and the global existence of the solution.

Since $J_i(0)>0$, every $J_i$ is bounded below by a positive constant on a
small ball. Finite convolution chains propagate positivity through space.
Pointwise irreducibility of $A(x)$ propagates positivity between components.
The variation-of-constants formula then gives strict positivity for every
positive time.
\end{proof}

\begin{lemma}
\label{lem:nonlocal-present-time}
Assume that (J) and (A) hold. Let $S(t)$ be the positive semigroup
generated by $\mathcal L$ on $C_{\rm b}(\mathbb R^N,\mathbb R^m)$.
Then there exist $p\in(0,1)$ and $K_H>0$ such that the following statements hold:
\begin{enumerate}
\item[(i)] Every bounded continuous function $\boldsymbol f\geq\boldsymbol0$ satisfies
\begin{equation}
\label{eq:linear-nonlocal-Harnack}
(S(1)\boldsymbol f)_j(x)
\leq K_H\|\boldsymbol f\|_\infty^{1-p}
       (S(1)\boldsymbol f)_i(x)^p
\end{equation}
for all $i,j\in\{1,\ldots,m\}$ and $x\in\mathbb R^N$.
\item[(ii)] Suppose in addition that (Q) holds. Then, for every $M>0$,
there exists $K_M>0$ such that every solution of
\eqref{eq:general-nonlocal-KPP} with
$\boldsymbol u_0\in C_{\rm b}(\mathbb R^N,[0,\infty)^m)$ and
$\|\boldsymbol u_0\|_\infty\leq M$ satisfies
\begin{equation}
\label{eq:nonlinear-component-comparison}
u_j(t,x)\leq K_M u_i(t,x)^p,
\qquad t\geq1,\quad x\in\mathbb R^N,\quad i,j\in\{1,\ldots,m\}.
\end{equation}
\end{enumerate}
\end{lemma}

\begin{proof}[Proof of Lemma~\ref{lem:nonlocal-present-time}]
Choose $\omega>0$ so large that $B(x):=A(x)+\omega I$
is entrywise nonnegative for every $x$. It remains irreducible. Define
\[
(\mathcal B\boldsymbol v)(x)=B(x)\boldsymbol v(x),
\qquad
(\mathcal T\boldsymbol v)_i=J_i*v_i+(\mathcal B\boldsymbol v)_i.
\]
Then $\mathcal T$ is a bounded positive operator and
\[
\mathcal L=\mathcal T-(1+\omega)I,
\qquad
S(1)=e^{-(1+\omega)}
\sum_{n=0}^{\infty}\frac{\mathcal T^n}{n!}
\]
in operator norm.

For every fixed $x$, irreducibility of the nonnegative matrix $B(x)$ implies
that $\sum_{r=0}^{m-1}B(x)^r$ is entrywise positive.
Continuity and periodicity give
\begin{equation*}
\eta:=\min_{x\in\mathbb T^N}\min_{i,j}
\sum_{r=0}^{m-1}(B(x)^r)_{ij}>0.
\end{equation*}
Fix $i,j,x$ and set $h_n=(\mathcal T^n\boldsymbol f)_j(x)\geq0$.
The positive operator $\mathcal B^r\mathcal T^n$ occurs in the expansion of $\mathcal T^{n+r}$. Each operator product is counted in this way for at most $m$ values of $r$. Since $(n+r)!\leq(n+m-1)!$, we obtain
\[
(S(1)\boldsymbol f)_i(x)
\geq \frac{e^{-(1+\omega)}}{m}
\sum_{n=0}^{\infty}\sum_{r=0}^{m-1}
\frac{(\mathcal B^r\mathcal T^n\boldsymbol f)_i(x)}{(n+r)!}
\geq \frac{e^{-(1+\omega)}\eta}{m}
\sum_{n=0}^{\infty}\frac{h_n}{(n+m-1)!}
=:\sum_{n=0}^{\infty}q_n.
\]
On the other hand,
\[
(S(1)\boldsymbol f)_j(x)
=e^{-(1+\omega)}\sum_{n=0}^{\infty}\frac{h_n}{n!}.
\]
We have $0\leq h_n\leq\|\mathcal T\|^n\|\boldsymbol f\|_\infty$.
Fix $p\in(0,1)$ and define
\[
d_n=e^{-(1+\omega)(1-p)}
\left(\frac{m}{\eta}\right)^p
\frac{(n+m-1)!^p}{n!}\|\mathcal T\|^{n(1-p)}.
\]
For every $n$, we have
\[
e^{-(1+\omega)}\frac{h_n}{n!}
\leq d_n\|\boldsymbol f\|_\infty^{1-p}q_n^p.
\]
By Stirling's formula, the logarithm of the $n$th term
$d_n^{1/(1-p)}$ is $-n\log n+O(n)$. Hence
\[
\sum_{n=0}^{\infty}d_n^{1/(1-p)}<\infty,
\]
and H\"older's inequality with conjugate exponents
$1/(1-p)$ and $1/p$ therefore gives
\begin{align*}
(S(1)\boldsymbol f)_j(x)
&\leq\|\boldsymbol f\|_\infty^{1-p}
\sum_{n=0}^{\infty}d_nq_n^p
\leq\|\boldsymbol f\|_\infty^{1-p}
\left(\sum_{n=0}^{\infty}d_n^{1/(1-p)}\right)^{1-p}
\left(\sum_{n=0}^{\infty}q_n\right)^p\\
&\leq K_H\|\boldsymbol f\|_\infty^{1-p}
(S(1)\boldsymbol f)_i(x)^p.
\end{align*}
Thus \eqref{eq:linear-nonlocal-Harnack} holds.

By Lemma~\ref{lem:KPP-wellposed-bounded}, there is $d_M>0$ such that
$(Q(x)\boldsymbol u(t,x))_i\leq d_M$ for every $i,t,x$.
On every interval $[t-1,t]$, comparison of
the positive linear evolution operators gives
\[
e^{-d_M}S(1)\boldsymbol u(t-1,\cdot)
\leq\boldsymbol u(t,\cdot)
\leq S(1)\boldsymbol u(t-1,\cdot).
\]
Combining \eqref{eq:linear-nonlocal-Harnack} with the uniform bound
in Lemma~\ref{lem:KPP-wellposed-bounded}
gives \eqref{eq:nonlinear-component-comparison}. This completes the proof.
\end{proof}

For $M\geq\|\boldsymbol u_0\|_\infty$, let
\[
\bar q=\max_{i,j,x}q_{ij}(x),
\qquad
g_M=m\bar qK_M.
\]
For $t\geq1$, \eqref{eq:nonlinear-component-comparison} yields
\begin{equation}
\label{eq:competition-power-control}
u_i(Q(x)\boldsymbol u)_i
\leq g_Mu_i^{1+p}.
\end{equation}
On the other hand,
\begin{equation}
\label{eq:competition-quadratic-control}
u_i(Q(x)\boldsymbol u)_i\geq q_0u_i^2.
\end{equation}
For a nonnegative vector $\boldsymbol v$ and $r>0$, write
$\boldsymbol v^{\circ r}=(v_1^r,\ldots,v_m^r)^{\mathsf T}$ for its
componentwise power. Consider the following two auxiliary cooperative KPP systems:
\begin{align}
\partial_t\boldsymbol v^+
&=\mathcal L\boldsymbol v^+-q_0(\boldsymbol v^+)^{\circ2},
\label{eq:upper-cooperative-KPP}\\
\partial_t\boldsymbol v^-
&=\mathcal L\boldsymbol v^--g_M(\boldsymbol v^-)^{\circ(1+p)}.
\label{eq:lower-cooperative-KPP}
\end{align}
Both systems have the same linearization $\mathcal L$ as
\eqref{eq:general-nonlocal-KPP}.
Let $\boldsymbol u$ solve \eqref{eq:general-nonlocal-KPP}, and choose
$M\geq\|\boldsymbol u_0\|_\infty$ in the definition of $g_M$.
By \eqref{eq:competition-quadratic-control} and the comparison
principle for \eqref{eq:upper-cooperative-KPP}, if
$\boldsymbol v^+(0)\geq\boldsymbol u_0$, then
\[
\boldsymbol u(t)\leq\boldsymbol v^+(t),\qquad t\geq0.
\]
Similarly, by \eqref{eq:competition-power-control} and comparison
for \eqref{eq:lower-cooperative-KPP}, if $t_0\geq1$ and
$\boldsymbol0\leq\boldsymbol v^-(0)\leq\boldsymbol u(t_0)$, then
\[
\boldsymbol v^-(t-t_0)\leq\boldsymbol u(t),
\qquad t\geq t_0.
\]
Both comparisons also apply after a translation of the initial time,
with $M$ chosen to bound the translated initial function.

We next establish some properties of cooperative KPP systems by using the theory of periodic monotone semiflows. These properties apply to both
\eqref{eq:upper-cooperative-KPP} and \eqref{eq:lower-cooperative-KPP}.

\begin{lemma}
\label{lem:cooperative-KPP-package}
Assume that (J) and (A) hold. Fix $q>0$, and let $b_i(x)$ be
continuous, positive and $\mathbb L$-periodic for every $i\in\Lambda$.
Write $\boldsymbol b=(b_1,\ldots,b_m)^{\mathsf T}$ and let
$\boldsymbol z_0\in C_{\rm b}(\mathbb R^N,[0,\infty)^m)$. Consider
\begin{equation}
\label{eq:abstract-cooperative-KPP}
\partial_tz_i=(\mathcal L\boldsymbol z)_i-b_i(x)z_i^{1+q}.
\end{equation}
Then, with $\boldsymbol z(0)=\boldsymbol z_0$, the following statements hold:
\begin{enumerate}
\item[(i)] If $\lambda'<0$, every bounded nonnegative solution converges to
$\boldsymbol0$ uniformly in space. If $\lambda'>0$, there is a unique
positive periodic equilibrium $\boldsymbol p$, and every nonzero periodic
solution satisfies
\[
\liminf_{t\to\infty}\min_i\min_{x\in\mathbb T^N}z_i(t,x)
\geq\min_i\min_{x\in\mathbb T^N}p_i(x)
=: \eta_{\rm per}(q,\boldsymbol b)>0.
\]
\item[(ii)] If $\lambda^*>0$, every nonzero bounded solution satisfies
\[
\liminf_{t\to\infty}\min_i\inf_{x\in K}z_i(t,x)
\geq\eta_{\rm per}(q,\boldsymbol b)>0
\]
for every compact set $K\subset\mathbb R^N$.
\item[(iii)] Let $\xi\in\mathcal I$. If $\lambda_L^*(\xi)<0$ and
$\boldsymbol z_0(x)\leq Ce^{-\xi\cdot x}\boldsymbol1$ for all $x\in\mathbb R^N$
and some $C>0$, then
$\boldsymbol z(t,\cdot)$ converges locally uniformly to zero. If
$\lambda_L^*(\xi)>0$ and condition~(i)(b) of
Theorem~\ref{thm:nonlocal-tail-criteria} holds,
and there exist $B\in\mathbb R$ and $C>0$ such that
\[
\boldsymbol z_0(x)\geq C^{-1}e^{-\xi\cdot x}\boldsymbol1
\quad\hbox{when }\xi\cdot x\geq B,
\]
then, for every $R>0$,
\[
\liminf_{t\to\infty}\min_i\inf_{|x|\leq R}z_i(t,x)
\geq\eta_{\rm per}(q,\boldsymbol b)>0.
\]
\item[(iv)] Suppose that all $J_i$ are compactly supported and $\lambda^*>0$. Then every
solution with nonzero compactly supported initial data spreads in direction
$e\in\mathbb S^{N-1}$ at the speed $c_{\rm FG}(e)$ defined in
\eqref{eq:planar-and-FG-speeds}.
\end{enumerate}
\end{lemma}

\begin{proof}[Proof of Lemma~\ref{lem:cooperative-KPP-package}]
We divide the proof into five steps.

\smallskip
\noindent Step 1. Order, strict positivity, and coefficient
approximation.
The only nonlinear terms in \eqref{eq:abstract-cooperative-KPP} are
diagonal, so its right-hand side is cooperative. Hence its solution operator
$Q_t$ is order preserving. The variation-of-constants formula, the fact that
$J_i$ is positive on a ball about the origin, and the pointwise
irreducibility of $A(x)$ imply
\[
Q_t\boldsymbol z_0\gg\boldsymbol0,
\qquad t>0,\quad \boldsymbol z_0\geq\boldsymbol0,\quad
\boldsymbol z_0\not\equiv\boldsymbol0.
\]
Moreover, for $0<\theta<1$,
\[
\mathcal L(\theta\boldsymbol z)
-\boldsymbol b\circ(\theta\boldsymbol z)^{\circ(1+q)}
-\theta\bigl(\mathcal L\boldsymbol z
-\boldsymbol b\circ\boldsymbol z^{\circ(1+q)}\bigr)
=
(\theta-\theta^{1+q})\boldsymbol b
\circ\boldsymbol z^{\circ(1+q)}\gg\boldsymbol0
\]
whenever $\boldsymbol z\gg\boldsymbol0$. Thus $Q_t$ is strongly
subhomogeneous: $Q_t(\theta\boldsymbol z)\gg\theta Q_t\boldsymbol z$
for $t>0$, $0<\theta<1$ and $\boldsymbol z\gg\boldsymbol0$.

Use the matrices $\underline A^\varepsilon$, $\overline A^\varepsilon$
and principal eigenpairs from Lemma~\ref{lem:approx}. In particular,
\begin{equation*}
0\leq\lambda_L^*(\mu)-\lambda_\varepsilon^-(\mu)\leq2\varepsilon,\qquad
0\leq\lambda_\varepsilon^+(\mu)-\lambda_L^*(\mu)\leq2\varepsilon,
\quad \mu\in\mathcal I.
\end{equation*}
Let $Q_t^{\varepsilon,-}$ and $Q_t^{\varepsilon,+}$ be obtained by replacing
$A$ with $\underline A^\varepsilon$ and $\overline A^\varepsilon$,
respectively. Comparison gives
\begin{equation}
\label{eq:cooperative-approximation-sandwich}
Q_t^{\varepsilon,-}\boldsymbol z_0
\leq Q_t\boldsymbol z_0
\leq Q_t^{\varepsilon,+}\boldsymbol z_0.
\end{equation}

\smallskip
\noindent Step 2. Periodic threshold dynamics.
Suppose first that $\lambda'<0$. Choose $\varepsilon$ so small that
$\lambda_\varepsilon^+(0)<0$, and let
$\overline{\boldsymbol\phi}^{\varepsilon}\gg0$ be the normalized periodic
principal eigenfunction. If
$M\overline{\boldsymbol\phi}^{\varepsilon}\geq\boldsymbol z_0$, then
$\boldsymbol W(t,x)=M e^{\lambda_\varepsilon^+(0)t}
\overline{\boldsymbol\phi}^{\varepsilon}(x)$ satisfies
\[
\partial_t\boldsymbol W-
\bigl(\overline{\mathcal L}_0^{\varepsilon}\boldsymbol W
-\boldsymbol b\circ\boldsymbol W^{\circ(1+q)}\bigr)
=\boldsymbol b\circ\boldsymbol W^{\circ(1+q)}\geq0.
\]
It is therefore a super-solution for the upper approximating system. By
\eqref{eq:cooperative-approximation-sandwich},
$0\leq\boldsymbol z(t)\leq\boldsymbol W(t)$, which proves uniform
extinction.

Suppose next that $\lambda'>0$. We briefly give the periodic equilibrium
argument, as in Bao \cite{Bao2020ZAMP}, for the present kernels and
nonlinearities. Choose $\varepsilon>0$ with
$\lambda_\varepsilon^-(0)>0$. For $a>0$ small and $M>0$ large,
$\boldsymbol\ell=a\underline{\boldsymbol\phi}^{\varepsilon}$ and
$M\boldsymbol1$ are ordered stationary sub- and super-solutions.
Thus $Q_t\boldsymbol\ell$ increases to a bounded periodic function
$\boldsymbol p\geq\boldsymbol\ell\gg0$. Dominated convergence in the
time-integrated equation on $[t,t+1]$ gives
\[
\boldsymbol h:=\mathcal J*\boldsymbol p,\qquad
\boldsymbol h+(A-I)\boldsymbol p
-\boldsymbol b\circ\boldsymbol p^{\circ(1+q)}=\boldsymbol0.
\]
The function $\boldsymbol h$ is continuous and strictly positive.
We record the maximum-ratio argument, including its use to establish
continuity. For fixed $x$, write $h=\boldsymbol h(x)\gg0$.
If $y,z\gg0$ solve
$h+(A(x)-I)y-\boldsymbol b(x)\circ y^{\circ(1+q)}=0$, let
$\theta=\max_i(y_i/z_i)$. If $\theta>1$, choose $i$ with
$y_i=\theta z_i$. Since $A(x)$ has nonnegative off-diagonal
entries and $y\leq\theta z$, we have
\[
\begin{aligned}
0&=h_i+[(A(x)-I)y]_i-b_i(x)y_i^{1+q}\\
 &\leq h_i+\theta[(A(x)-I)z]_i
 -b_i(x)\theta^{1+q}z_i^{1+q}\\
 &= (1-\theta)h_i+b_i(x)
 (\theta-\theta^{1+q})z_i^{1+q}<0,
\end{aligned}
\]
a contradiction. Interchanging $y,z$ shows equality. Thus the positive
solution of the pointwise algebraic equation is unique.
The monotone limit $\boldsymbol p$ is bounded and bounded away from
zero. If $x_n\to x$, every convergent subsequence of
$\boldsymbol p(x_n)$ solves the same limiting algebraic equation at
$x$, because $A$, $\boldsymbol b$, and $\boldsymbol h$ are continuous.
Pointwise uniqueness therefore implies
$\boldsymbol p(x_n)\to\boldsymbol p(x)$, proving continuity.
For two positive periodic equilibria $\boldsymbol p,\widetilde{\boldsymbol p}$, set
$\theta=\max_{i,x}p_i(x)/\widetilde p_i(x)$.
At a contact point $(i,x)$ the positivity of convolution and of the
off-diagonal matrix entries gives
$(\mathcal L\boldsymbol p)_i(x)\leq\theta(\mathcal L\widetilde{\boldsymbol p})_i(x)$.
If $\theta>1$, their equilibrium equations imply
$0\leq b_i(x)(\theta-\theta^{1+q})
\widetilde p_i(x)^{1+q}<0$, a contradiction. Interchanging
the equilibria proves uniqueness of the positive periodic equilibrium.
Finally, for
\[
c=\min_{i,x}b_i(x)p_i(x)^q>0,\qquad
\alpha'=c\alpha(1-\alpha^q),
\]
$\alpha(t)\boldsymbol p$ is a sub-solution when $0<\alpha(0)<1$
and a super-solution when $\alpha(0)>1$. Strict positivity at positive
times allows every nonzero nonnegative periodic solution to be
sandwiched between these functions. Since $\alpha(t)\to1$, convergence
to $\boldsymbol p$ is uniform. This argument uses only $q>0$ and
positivity of the kernels, not second derivatives of the reaction at zero.
This proves assertion~(i), with
\[
\eta_{\rm per}(q,\boldsymbol b)
=\min_i\min_{x\in\mathbb T^N}p_i(x)>0.
\]

We next show how local persistence gives the same lower bound. Let
$\chi_L\in C_{\rm c}(\mathbb R^N,[0,1])$ equal one on $B_L$ and vanish
outside $B_{2L}$; here $C_{\rm c}$ denotes the space of continuous
compactly supported functions. For $0<\theta<1$, the nonlinear diagonal
loss in the difference equation for the two ordered solutions is
nonnegative. Comparison therefore gives
\[
0\leq Q_\tau(\theta\boldsymbol p)
-Q_\tau\bigl(\theta\boldsymbol p\chi_L(\cdot-y)\bigr)
\leq S(\tau)\bigl[\theta\boldsymbol p(1-\chi_L(\cdot-y))\bigr].
\]
Put $k=\sum_iJ_i$ and choose
$a_0\geq0$ bounding the row sums of $A$. For a nonnegative vector function
$\boldsymbol f$,
\[
S(\tau)\boldsymbol f
\leq e^{a_0\tau}\sum_{n=0}^\infty\frac{\tau^n}{n!}
\bigl(k^{*n}*\max_i f_i\bigr)\boldsymbol1,
\qquad k^{*0}=\delta_0.
\]
Here $\delta_0$ is the unit Dirac measure at the origin.
The measure kernel on the right has finite total mass. Hence, for fixed
$\tau,R>0$, there are numbers $\varepsilon_L(\tau,R)\to0$ such that
\begin{equation}\label{eq:cooperative-cutoff-comparison}
Q_\tau\bigl(\theta\boldsymbol p\chi_L(\cdot-y)\bigr)(x)
\geq Q_\tau(\theta\boldsymbol p)(x)
-\varepsilon_L(\tau,R)\boldsymbol1,
\qquad |x-y|\leq R,
\end{equation}
uniformly in $y\in\mathbb R^N$.

Suppose that a solution satisfies, for some $\delta>0$, $c\geq0$ and
$e\in\mathbb S^{N-1}$,
\[
\liminf_{t\to\infty}\min_i\inf_{|x|\leq R}z_i(t,x+cte)
\geq\delta
\quad\hbox{for every }R>0.
\]
Choose $0<\theta<1$ with $\theta\|\boldsymbol p\|_\infty<\delta$.
For fixed $L,\tau$, the solution at time $t-\tau$ eventually dominates
$\theta\boldsymbol p\chi_L(\cdot-(t-\tau)ce)$. Applying \eqref{eq:cooperative-cutoff-comparison}
with radius $R+c\tau$ and then letting $L\to\infty$ gives
\[
\liminf_{t\to\infty}\min_i\inf_{|x|\leq R}
\bigl[z_i(t,x+cte)-(Q_\tau(\theta\boldsymbol p))_i(x+cte)\bigr]
\geq0.
\]
Since $Q_\tau(\theta\boldsymbol p)\to\boldsymbol p$ uniformly as
$\tau\to\infty$, we conclude that
\begin{equation*}
\liminf_{t\to\infty}\min_i\inf_{|x|\leq R}
\bigl[z_i(t,x+cte)-p_i(x+cte)\bigr]\geq0.
\end{equation*}
Thus the lower bound can be taken to be $\eta_{\rm per}(q,\boldsymbol b)$,
both on fixed balls and on moving balls.

\smallskip
\noindent Step 3. Exponential initial tails.
Assume $\lambda_L^*(\xi)<0$. Choose the upper approximation so that
$\lambda_\varepsilon^+(\xi)<0$, and let
$\overline{\boldsymbol\phi}^{\varepsilon}_\xi\gg0$ be its normalized
exponentially weighted eigenfunction. For a sufficiently large $M$,
\[
\boldsymbol W(t,x)=M
e^{\lambda_\varepsilon^+(\xi)t-\xi\cdot x}
\overline{\boldsymbol\phi}^{\varepsilon}_\xi(x)
\]
dominates the initial function. Work in the exponentially weighted space
\begin{equation}\label{eq:exponentially-weighted-space}
E_\xi:=\{\boldsymbol f:e^{\xi\cdot x}\boldsymbol f(x)
\in C_{\rm b}(\mathbb R^N,\mathbb R^m)\},\qquad
\|\boldsymbol f\|_{E_\xi}
=\|e^{\xi\cdot x}\boldsymbol f\|_\infty.
\end{equation}
Since $\xi\in\mathcal I$, convolution is bounded on $E_\xi$, with norm
at most $\max_i M_i(\xi)$. The bounded cooperative multiplication part
then makes $\mathcal L$ a bounded generator of a positive semigroup
on this space. Its exponential series agrees pointwise with the one
on $C_{\rm b}$ for data in $C_{\rm b}\cap E_\xi$, so we use the same
notation $S(t)$ for both semigroups.
The initial tail bound gives $\boldsymbol z_0\in C_{\rm b}\cap E_\xi$,
while $\boldsymbol W$ is continuously differentiable in $E_\xi$ and
satisfies
\[
\partial_t\boldsymbol W-\mathcal L\boldsymbol W
=(\overline A^\varepsilon-A)\boldsymbol W\geq\boldsymbol0.
\]
Variation of constants for $\boldsymbol z$ in $C_{\rm b}$ and
nonnegativity of its loss term give the first inequality below;
linear comparison in $E_\xi$ gives the second:
\begin{equation}\label{eq:weighted-linear-comparison}
\boldsymbol z(t)\leq S(t)\boldsymbol z_0\leq\boldsymbol W(t).
\end{equation}
Consequently the strict negativity of
$\lambda_\varepsilon^+(\xi)$ implies local uniform extinction.

For the persistence assertion, condition~(i)(b) of
Theorem~\ref{thm:nonlocal-tail-criteria} and convexity allow us to
choose $s>1$ sufficiently close to one that
$0<\lambda_L^*(s\xi)<\lambda_L^*(\xi)$.
To see this, interpolate between $1$ and a weight parameter satisfying
that condition. The convex upper bound is strictly below
$\lambda_L^*(\xi)$, while continuity of the finite convex function
at $1$ ensures positivity for parameters sufficiently close to $1$.
Convexity gives
$\lambda'=\lambda_L^*(0)\geq
(s\lambda_L^*(\xi)-\lambda_L^*(s\xi))/(s-1)>0$.
By the lower
approximation, the same strict inequalities hold for
$\lambda_\varepsilon^-(s\xi)$ and
$\lambda_\varepsilon^-(\xi)$ when $\varepsilon$ is small. Write
\[
E_r(x)=e^{-r\xi\cdot x}
\underline{\boldsymbol\phi}_{r\xi}^\varepsilon(x),
\qquad
\underline{\mathcal L}_0^\varepsilon E_r
=\lambda_\varepsilon^-(r\xi)E_r,
\quad r\in\{1,s\}.
\]
Choose $0<\kappa<\lambda_\varepsilon^-(s\xi)$ and set
\[
\alpha=\lambda_\varepsilon^-(\xi)-\kappa
>\beta=\lambda_\varepsilon^-(s\xi)-\kappa>0.
\]
For $K>0$, let $\boldsymbol h_*=(E_1-KE_s)^+$.
The periodic positive bounds on the eigenfunctions give, on
$\{E_{1,i}>KE_{s,i}\}$,
\[
e^{-(s-1)\xi\cdot x}\leq C_1/K,
\qquad 0<h_{*,i}(x)\leq C_2K^{-1/(s-1)}.
\]
Since $h_{*,i}=0$ elsewhere, we obtain
\begin{equation}
\label{eq:two-exponential-smallness}
\|\boldsymbol h_*\|_\infty\longrightarrow0
\quad\hbox{as }K\to\infty.
\end{equation}
Fix $K$ so large that $b^*\|\boldsymbol h_*\|_\infty^q\leq\kappa$.
This choice is independent of $\boldsymbol z_0$, $C$ and $B$.

For $T>0$ and $0\leq t\leq T$, define
\begin{equation}\label{eq:exponential-tail-subsolution}
\boldsymbol H_T(t,x)
=e^{-\alpha(T-t)}E_1(x)-K e^{-\beta(T-t)}E_s(x).
\end{equation}
Since $\alpha>\beta$, we have the componentwise bound and derivative identity
\[
\boldsymbol0\leq\boldsymbol H_T^+(t,x)
\leq e^{-\alpha(T-t)}\boldsymbol h_*(x)\leq\boldsymbol h_*(x),\qquad
\partial_t\boldsymbol H_T
=\underline{\mathcal L}_0^\varepsilon\boldsymbol H_T
-\kappa\boldsymbol H_T.
\]
At positive components, positivity of the kernels and cooperativity give
\[
\partial_tH_{T,i}^+
\leq(\underline{\mathcal L}_0^\varepsilon\boldsymbol H_T^+)_i
-\kappa H_{T,i}^+
\leq(\underline{\mathcal L}_0^\varepsilon\boldsymbol H_T^+)_i
-b_i(x)(H_{T,i}^+)^{1+q}.
\]
At zero components, the integral and off-diagonal terms are nonnegative
and dominate their values before taking positive parts. Hence the same
inequality holds with the upper right derivative. Thus
$\boldsymbol H_T^+$ is a bounded generalized sub-solution on $[0,T]$.
The construction \eqref{eq:exponential-tail-subsolution}, with
\eqref{eq:two-exponential-smallness}, is adapted from
\cite[Section~3.5]{Girardin2025}.

For all sufficiently large $T$,
\[
e^{-\alpha T}\max_{i,x}\underline\phi_{\xi,i}^\varepsilon(x)\leq C^{-1},
\qquad
K e^{(\alpha-\beta)T}
\min_{i,x}\frac{\underline\phi_{s\xi,i}^\varepsilon(x)}
{\underline\phi_{\xi,i}^\varepsilon(x)}
\geq e^{(s-1)B}.
\]
The second inequality makes $\boldsymbol H_T^+(0,x)=\boldsymbol0$ when
$\xi\cdot x\leq B$, and the first matches the assumed exponential tail.
Hence $\boldsymbol H_T^+(0)\leq\boldsymbol z_0$. Comparison gives
\[
Q_T^{\varepsilon,-}\boldsymbol z_0\geq\boldsymbol h_*,
\qquad
Q_{T+1}^{\varepsilon,-}\boldsymbol z_0
\geq Q_1^{\varepsilon,-}\boldsymbol h_*
\]
for all sufficiently large $T$. For each fixed $\ell\in\mathbb L$, the
same argument applies to $\boldsymbol z_0(\cdot+\ell)$, whose tail constants
are $B-\xi\cdot\ell$ and $Ce^{\xi\cdot\ell}$. Periodicity therefore gives
\[
Q_{T+1}^{\varepsilon,-}\boldsymbol z_0(x)
\geq (Q_1^{\varepsilon,-}\boldsymbol h_*)(x-\ell)
\]
for all sufficiently large $T$, with the waiting time allowed to depend
on $\ell$. Pointwise strict positivity, continuity and
\eqref{eq:cooperative-approximation-sandwich}, together with the fact that
every fixed ball is covered by finitely many cells
$\ell+\overline{\mathcal C}$, give
\[
\liminf_{t\to\infty}\min_i\inf_{|x|\leq R}z_i(t,x)
\geq\min_{i\in\Lambda,\,x\in\overline{\mathcal C}}
(Q_1^{\varepsilon,-}\boldsymbol h_*)_i(x)>0
\quad\hbox{for every }R>0.
\]
The argument in Step~2 now gives the lower bound
$\eta_{\rm per}(q,\boldsymbol b)$, proving assertion~(iii).

\smallskip
\noindent Step 4. Local persistence by bounded-domain sub-solutions.
Assume $\lambda^*>0$. By Theorem~\ref{TH2.1}(ii), choose $l$ so large
that $\overline{\mathcal C}\subset\Omega:=(-l,l)^N$ and
$\lambda(l)>0$. Since $\Omega$ contains a full closed period cell,
$\max_{x\in\overline\Omega}s(x)=s^*$. Choose $0<\delta<\lambda(l)$ and set
\[
A_\delta(x)=A(x)-\bigl(s(x)-s^*+\delta\bigr)_+I,
\qquad x\in\overline\Omega.
\]
Then $A-\delta I\leq A_\delta\leq A$, and the pointwise Perron root
of $A_\delta-I$ equals
$\min\{s(x),s^*-\delta\}$.
It therefore attains its maximum on a nonempty open subset of $\Omega$.
The bounded-domain nonintegrability criterion in \cite[Lemma~2.1]{WZ-1},
together with the compact-operator argument on cubes used in the proof of
Theorem~\ref{TH2.1}(ii), gives a positive continuous principal eigenfunction
$\boldsymbol\phi$ of
\[
(\mathcal L_\Omega^-\boldsymbol\phi)_i(x)
=\int_\Omega J_i(x-y)\phi_i(y)\,\mathrm dy
-\phi_i(x)+(A_\delta(x)\boldsymbol\phi(x))_i.
\]
In this application the kernel is only used on the bounded difference set
$\overline\Omega-\overline\Omega$; it can be cut off outside that set
without changing the operator. Monotonicity and the diagonal shift give
its eigenvalue $\lambda_\Omega\geq\lambda(l)-\delta>0$.
Normalize $\max_{i\in\Lambda,\,x\in\overline\Omega}\phi_i(x)=1$. If
\[
a'=\lambda_\Omega a-b^*a^{1+q},\qquad a(0)>0,
\]
then direct substitution gives
\[
\partial_t(a\phi_i)
-\bigl((\mathcal L_\Omega^-[a\boldsymbol\phi])_i
-b_i(x)(a\phi_i)^{1+q}\bigr)
=a^{1+q}\phi_i\bigl(b_i(x)\phi_i^q-b^*\bigr)\leq0.
\]
Thus $a(t)\boldsymbol\phi$ is a sub-solution of the cooperative problem
on $\Omega$ with zero exterior data.
Moreover $a(t)\to a_\infty=(\lambda_\Omega/b^*)^{1/q}>0$.

For any nonzero bounded initial function, fix $t_0>0$. Pointwise strict positivity
gives $\boldsymbol z(t_0,x)\gg0$. For each $p\in\mathbb L$, periodicity
allows us to translate the preceding bounded-domain sub-solution to
$\Omega+p$. Choose its initial amplitude $a_p(0)>0$ small enough that
$a_p(0)\boldsymbol\phi(x-p)\leq\boldsymbol z(t_0,x)$ on
$\overline\Omega+p$. Since the omitted exterior convolution terms are
nonnegative and $A_\delta\leq A$, the restriction of the whole-space
solution to this cube is a super-solution of the lower problem. Bounded-domain
comparison consequently gives
\[
\boldsymbol z(t_0+t,x)\geq a_p(t)\boldsymbol\phi(x-p),
\qquad x\in\overline\Omega+p,\quad t\geq0.
\]
Every compact set is covered by finitely many translates
$p+\overline{\mathcal C}$. Since every amplitude tends to the same
$a_\infty$, independently of $a_p(0)$, we obtain
\[
\liminf_{t\to\infty}\min_i\inf_{x\in K}z_i(t,x)
\geq a_\infty
\min_{i\in\Lambda,\,x\in\overline{\mathcal C}}\phi_i(x)
>0.
\]
Applying Step~2 upgrades this lower bound to
$\eta_{\rm per}(q,\boldsymbol b)$ and proves assertion~(ii).

\smallskip
\noindent Step 5. Spreading for compactly supported kernels.
Assume that all $J_i$ are compactly supported and $\lambda^*>0$.
Fix the lower coefficient approximation $\underline A^\varepsilon$
from Step~1, with $2\varepsilon<\lambda^*$, and let
$S_\varepsilon(t)$ be the semigroup generated by
$\underline{\mathcal L}_0^\varepsilon$. By Lemma~\ref{lem:approx},
\[
\inf_{\mu\in\mathbb R^N}\lambda_\varepsilon^-(\mu)
\geq\lambda^*-2\varepsilon>0.
\]
Let $\boldsymbol\pi=M_{\rm sup}\boldsymbol1$ be a common constant super-solution
for the original cooperative system and its lower approximation,
dominating the initial function. Choose
$a>\max\{0,\max_{i,x}(1-\underline a_{ii}^\varepsilon(x))\}$, so that
$\underline A^\varepsilon-I+aI$ is entrywise nonnegative and
$\mathcal P_\varepsilon:=\underline{\mathcal L}_0^\varepsilon+aI$ is positive.
Choose $0<2d<\inf_\mu\lambda_\varepsilon^-(\mu)$ and let
$C=\max_{0\leq s\leq1}\|S_\varepsilon(s)\|$.
Choose $\eta>0$ with $C\eta\leq M_{\rm sup}$ and
$b^*(C\eta)^q\leq d$. For $0<h\leq1$, define
\begin{equation}\label{eq:cooperative-discrete-comparison-maps}
 T_h=e^{-(a+d)h}(I+h\mathcal P_\varepsilon),
 \qquad
 F_h[\boldsymbol u]=T_h[\min\{\boldsymbol u,\eta\boldsymbol1\}],
 \qquad \boldsymbol u\geq\boldsymbol0.
\end{equation}
The minimum is componentwise. These maps are positive and commute with
lattice translations. Since the kernels have compact support, their
values on a compact set $K$ depend only on the input on
$K+\overline B_R=\{x+z:x\in K,\ |z|\leq R\}$, where $R$
contains the supports of all kernels.
Both maps are therefore continuous in the compact-open topology,
that is, the topology of uniform convergence on compact sets. The same
argument gives this continuity for $T_h$ on exponentially weighted
continuous functions.

For $\boldsymbol v=\min\{\boldsymbol u,\eta\boldsymbol1\}$,
linear comparison gives $Q_s^{\varepsilon,-}\boldsymbol v\leq
C\eta\boldsymbol1$ for $0\leq s\leq1$. The nonlinear loss is
therefore bounded above by $dQ_s^{\varepsilon,-}\boldsymbol v$.
Comparison and the positive exponential series yield
\begin{equation}\label{eq:cooperative-discrete-lower-comparison}
 Q_h^{\varepsilon,-}\boldsymbol u
 \geq e^{-dh}S_\varepsilon(h)\boldsymbol v
 =e^{-(a+d)h}e^{h\mathcal P_\varepsilon}\boldsymbol v
 \geq F_h[\boldsymbol u].
\end{equation}
More precisely, using the positive operator inequality
$T_h\leq e^{-dh}S_\varepsilon(h)$, valid because
$I+h\mathcal P_\varepsilon\leq e^{h\mathcal P_\varepsilon}$,
we have
$\boldsymbol0\leq F_h[\boldsymbol u]
\leq e^{-dh}S_\varepsilon(h)[\eta\boldsymbol1]
\leq C\eta\boldsymbol1\leq\boldsymbol\pi$
for all $\boldsymbol0\leq\boldsymbol u\leq\boldsymbol\pi$.
By positivity of $T_h$ and the definition of $F_h$, we have
$F_h[\boldsymbol u]\leq T_h\boldsymbol u$ for every
$\boldsymbol0\leq\boldsymbol u\leq\boldsymbol\pi$.
More precisely, given any $\theta\in(0,1)$, take the
strictly positive constant vector
$\boldsymbol\epsilon_\theta=\eta\boldsymbol1$.  Then
\[
F_h[\boldsymbol u]
=T_h[\min\{\boldsymbol u,\boldsymbol\epsilon_\theta\}]
\geq (1-\theta)T_h[\min\{\boldsymbol u,
\boldsymbol\epsilon_\theta\}],
\qquad \boldsymbol0\leq\boldsymbol u\leq\boldsymbol\pi.
\]
Thus (A7) and (A8) of \cite[Section~2.2]{DL2015}
hold. For every $\delta>0$ one may take
$\boldsymbol\epsilon_\delta=\eta\boldsymbol1$ in (A8):
if $0<\delta<1$, its conclusion follows with the stronger
factor $1$, while for $\delta\geq1$ its right-hand side is
nonpositive and (A8) holds by positivity of $F_h$.
The constant periodic supersolution $\boldsymbol\pi=M_{\rm sup}\boldsymbol1$
provides (A3), while lattice equivariance and finite propagation
of $F_h$ and $T_h$ establish (A1)--(A2).  No convergence to a
positive equilibrium of $F_h$ is invoked: the stronger assumptions
(A9)--(A10) are not needed for the lower-spreading conclusion
of \cite[Theorem~2.11(ii)]{DL2015}.

We next check the spectral hypothesis for $T_h$. Set
\[
 s_\varepsilon^*=\max_x
 s_{\rm PF}(\underline A^\varepsilon(x)-I).
\]
By Lemma~\ref{lem:approx}, each $\underline{\mathcal L}_\mu^\varepsilon$ has a
positive periodic eigenfunction, with eigenvalue
$\lambda_\varepsilon^-(\mu)>s_\varepsilon^*$.
Write $(T_h)_\mu\boldsymbol f=e^{\mu\cdot x}
T_h[e^{-\mu\cdot x}\boldsymbol f]$ for the exponentially weighted
operator on periodic functions. Its corresponding eigenvalue is
\[
 r_h(\mu)=e^{-(a+d)h}
 \bigl(1+h(a+\lambda_\varepsilon^-(\mu))\bigr).
\]
Its positive eigenfunction has a positive minimum on the torus and is
therefore an order unit: every vector in $X$ is bounded in absolute
value by a multiple of this eigenfunction. Positivity of $(T_h)_\mu$
then implies that $r_h(\mu)$ equals its spectral radius. The convolution part is compact on the periodic space, including
for the exponentially weighted kernels. Modulo the compact operators,
$(T_h)_\mu$ is multiplication by
$e^{-(a+d)h}\{I+h[\underline A^\varepsilon(x)+(a-1)I]\}$.
The pointwise matrices are nonnegative and irreducible. Their Perron
eigenvalues are
$e^{-(a+d)h}[1+h(a+s_{\rm PF}(\underline A^\varepsilon(x)-I))]$.
The essential spectral radius $r_{\rm ess}$, namely the spectral radius
in the quotient by compact operators, is therefore
bounded by their maximum, and hence
\begin{equation}\label{eq:cooperative-discrete-spectral-gap}
 r_{\rm ess}((T_h)_\mu)
 \leq e^{-(a+d)h}\bigl(1+h(a+s_\varepsilon^*)\bigr)
 <r_h(\mu).
\end{equation}
For eventual strong positivity, fix a radius $r>0$ on which every
kernel has a common strictly positive lower bound. On the compact
period torus, a sufficiently high common power of the diagonal
convolution part maps each nontrivial nonnegative component function
to an everywhere positive function: finitely many jumps of length
less than $r$ connect every two points of the torus. The multiplication
matrix $B(x)=I+h[\underline A^\varepsilon(x)+(a-1)I]$ has positive
diagonal and is irreducible for every $x$; its $(m-1)$st power is
entrywise positive. In the nonnegative operator expansion of a
sufficiently high power of $(T_h)_\mu$, retain the term with the
needed convolution steps followed by $m-1$ multiplication steps.
This term is strongly positive, and consequently so is the full
power. The exponent is independent of the input function: the kernel
positivity radius is fixed, finitely many jumps cover the compact
torus, and $m-1$ multiplication steps connect all components at each
point. Continuity and compactness also give a common positive lower
bound for the entries of $B(x)^{m-1}$. Because this power maps every nonzero element of the
periodic positive cone into its interior, the eigenvalue $r_h(\mu)$
is algebraically simple and satisfies
$|\zeta|<r_h(\mu)$ for every other spectral value
$\zeta$ of $(T_h)_\mu$.  This strict spectral inequality,
together with \eqref{eq:cooperative-discrete-spectral-gap},
establishes the spectral gap in the sense of
\cite[Definition~2.6]{DL2015}.
For all sufficiently small $h$, $r_h(0)>1$.

In units of continuous time, define
\[
 f_{h,d}(\mu)=\frac{\log r_h(\mu)}{h},\qquad
 c_{h,d}(e)=\inf_{\rho>0}\frac{f_{h,d}(\rho e)}{\rho},\qquad
 \mathcal W_{h,d}=\bigcap_{e\in\mathbb S^{N-1}}
 \{x:x\cdot e\leq c_{h,d}(e)\}.
\]
We next justify passage to $h=0$ in these speeds. On every bounded
set of weights,
\[
 f_{h,d}(\mu)
 =\frac{\log(1+h(a+\lambda_\varepsilon^-(\mu)))}h-a-d
 \longrightarrow\lambda_\varepsilon^-(\mu)-d
\]
uniformly as $h\downarrow0$. Since $a+\lambda_\varepsilon^-(\mu)>0$,
the convergence is monotone increasing. Positivity of the kernels
on a fixed ball gives constants $\alpha,\beta,D>0$ such that
$\lambda_\varepsilon^-(\mu)\geq\alpha e^{\beta|\mu|}-D$.
For fixed $h_0>0$ this implies
\[
 \liminf_{\rho\to\infty}\inf_{e\in\mathbb S^{N-1}}
 \frac{f_{h_0,d}(\rho e)}{\rho}\geq\frac\beta{h_0}.
\]
Choose $h_0$ so small that $\inf_\mu f_{h_0,d}(\mu)>0$ and
$\beta/h_0>1+\max_e\lambda_\varepsilon^-(e)$.
For $0<h\leq h_0$, monotonicity gives $f_{h,d}\geq f_{h_0,d}$,
whereas testing at $\rho=1$ gives a common upper bound for
$c_{h,d}(e)$. These estimates exclude both $\rho\downarrow0$ and
$\rho\to\infty$, uniformly in $h\leq h_0$ and $e$.
Thus all minimizing weights lie in a fixed compact subinterval of
$(0,\infty)$, and
\[
 c_{h,d}(e)\longrightarrow
 c_{\varepsilon,d}(e):=
 \inf_{\rho>0}\frac{\lambda_\varepsilon^-(\rho e)-d}{\rho}
 \quad\hbox{uniformly in }e.
\]
The same bounds show uniform convergence
$c_{\varepsilon,d}\to c_\varepsilon$ as $d\downarrow0$, where
$c_\varepsilon(e)=\inf_{\rho>0}\lambda_\varepsilon^-(\rho e)/\rho$.
Finally, $0\leq\lambda_L^*(\mu)-\lambda_\varepsilon^-(\mu)
\leq2\varepsilon$, together with \eqref{eq:weighted-moment-lower-bound}
and positivity of the kernels on a fixed ball, give
$c_\varepsilon\to c^*$ uniformly in $e$ as $\varepsilon\downarrow0$.
For the parameters chosen above, $c_{h,d}$ has a positive lower bound
and a finite upper bound on the sphere. Thus $\mathcal W_{h,d}$
contains a ball about zero and is bounded.

We now record precisely the remaining hypotheses needed
for the multidimensional result of Ding and Liang.
For the controlling operator $M=T_h$, positivity and eventual
strong positivity on periodic functions establish (A4).
Compact support of the kernels permits exponential conjugation
for every real weight and establishes (A5).  The argument above
shows that, for each weight, $(T_h)_\mu$ has an algebraically
simple spectral-radius eigenvalue, separated in modulus from the
remaining spectrum, which is (A6); also $r_h(0)>1$ for the chosen
small $h$.  Together with (A1)--(A3) and (A8) verified above,
these are exactly the assumptions used in
\cite[Theorem~2.11(ii)]{DL2015}.  Its statement does not require
compactness of $F_h$, nor does the lower-spreading conclusion
require (A9)--(A10).
The multidimensional spreading theorem of
\cite[Theorem~2.11(ii)]{DL2015}, applied to $F_h$ with $M=T_h$,
gives persistence for its iterates on $nhK$ whenever
$K\Subset\operatorname{int}\mathcal W_{h,d}$.
Here the discrete-time speeds are $h c_{h,d}(e)$.
The theorem requires initial positivity on a sufficiently large fixed
ball. For nonzero nonnegative initial data $\boldsymbol z_0$,
pointwise strict positivity makes $Q_1^{\varepsilon,-}\boldsymbol z_0$
strictly positive on every fixed ball.
Order preservation and $Q_h^{\varepsilon,-}\geq F_h$ give
\begin{equation}\label{eq:cooperative-iterate-comparison}
Q_{1+nh}^{\varepsilon,-}\boldsymbol z_0
=Q_{nh}^{\varepsilon,-}[Q_1^{\varepsilon,-}\boldsymbol z_0]
\geq F_h^n[Q_1^{\varepsilon,-}\boldsymbol z_0].
\end{equation}
This transfers the lower bound to the cooperative solution.
For any $K\Subset\operatorname{int}\mathcal W$, with
\[
 \mathcal W=\bigcap_{e\in\mathbb S^{N-1}}
 \{x:x\cdot e\leq c^*(e)\},
\]
the uniform speed limits allow us to choose successively
$\varepsilon$, $d$, and $h$ such that
$K\Subset\operatorname{int}\mathcal W_{h,d}$.
Coefficient comparison then gives persistence for the original
cooperative system at times $1+nh$. Choose compact convex sets
$K\Subset K'\Subset\operatorname{int}\mathcal W$, with $K'$ containing
a neighborhood of $K$. For the bounded cooperative solution,
positivity of the integral operator and the off-diagonal coupling,
together with the uniform upper bound on the nonlinear loss, imply
$\partial_tz_i\geq-Dz_i$ for a constant $D$ independent of $t,x$;
hence $z_i(t+s,x)\geq e^{-Ds}z_i(t,x)$ whenever $0\leq s\leq h$.
Writing any sufficiently large time as $1+nh+s$ with $0\leq s<h$
and using the discrete-time lower bound transferred by
\eqref{eq:cooperative-iterate-comparison} on $tK'$ proves the desired
lower bound on $tK$. The enlargement from $K$ to $K'$ absorbs the
bounded shift in $t$; in particular no uniform-in-$x$ continuity of
the nonlinear semiflow is assumed. The radial boundary of
$\mathcal W$ is given by \eqref{eq:planar-and-FG-speeds}, which
proves the ray lower bound.

To prove the ray upper bound, let $c>c_{\rm FG}(e)$ and choose $e'\cdot e>0$ with
$c(e'\cdot e)>c^*(e')$, then choose $\rho>0$, $\gamma<\rho c(e'\cdot e)$
and $\boldsymbol\phi\in X^{++}$ such that
$\mathcal L_{\rho e'}\boldsymbol\phi\leq\gamma\boldsymbol\phi$.
For large $C_0$, the function
\[
\boldsymbol Z(t,x)=C_0e^{-\rho x\cdot e'+\gamma t}\boldsymbol\phi(x)
\]
dominates the compactly supported initial function and is a linear
super-solution. Comparison in the weighted space
\eqref{eq:exponentially-weighted-space}, as in
\eqref{eq:weighted-linear-comparison}, gives $\boldsymbol z\leq\boldsymbol Z$, and
$\boldsymbol Z(t,x+cte)\to0$ uniformly for bounded $x$.
This proves assertion~(iv) and completes the proof.
\end{proof}

\subsubsection{Extinction and persistence}

We now apply Lemma~\ref{lem:cooperative-KPP-package} through the upper
and lower comparisons to prove the stated nonlinear results.

\begin{proof}[\textbf{Proof of Theorem~\ref{thm:nonlocal-threshold-dynamics}}]
(i)
Let $\boldsymbol v^+$ solve \eqref{eq:upper-cooperative-KPP} with
$\boldsymbol v^+(0)=\boldsymbol u_0$. By comparison with
\eqref{eq:upper-cooperative-KPP},
$\boldsymbol0\leq\boldsymbol u(t)\leq\boldsymbol v^+(t)$. Since
$\lambda'<0$, Lemma~\ref{lem:cooperative-KPP-package}(i) gives uniform
convergence of $\boldsymbol v^+$ to zero. For completeness, let $\boldsymbol u$ be a nonnegative entire
solution with $\sup_{t,x}\|\boldsymbol u(t,x)\|<\infty$.
The upper coefficient approximation used in Step~2 of the proof of
Lemma~\ref{lem:cooperative-KPP-package} has a positive periodic
eigenfunction $\overline{\boldsymbol\phi}^\varepsilon$ and
negative eigenvalue $\lambda_\varepsilon^+(0)$. Since the
entire solution has a common bound, there exists $C>0$, independent
of $T>0$, such that comparison from time $-T$ gives
\[
\boldsymbol0\leq\boldsymbol u(0,x)
\leq C e^{\lambda_\varepsilon^+(0)T}
\overline{\boldsymbol\phi}^\varepsilon(x).
\]
Letting $T\to\infty$ yields $\boldsymbol u(0,\cdot)=0$.
The same argument after every time translation proves that the entire
solution is identically zero.

\smallskip
\noindent(ii)
By Lemma~\ref{lem:KPP-wellposed-bounded}, there is a time
$T=T(\boldsymbol u_0)$ such that
\[
\|\boldsymbol u(t)\|_\infty\leq M_*+1,
\qquad t\geq T.
\]
Fix $t_0=T+1$. Periodicity is preserved by \eqref{eq:general-nonlocal-KPP}, and
$\boldsymbol u(t_0,\cdot)$ is strictly positive on the compact torus.
Apply \eqref{eq:nonlinear-component-comparison} to the solution with
initial time $T$ and initial bound $M=M_*+1$, and put
$g_*=m\bar qK_{M_*+1}$.
Let $\boldsymbol v^-$ solve
\[
\partial_t\boldsymbol v^-
=\mathcal L\boldsymbol v^--g_*(\boldsymbol v^-)^{\circ(1+p)},\qquad
\boldsymbol v^-(0)=
\left(\min_{1\leq i\leq m}\min_{x\in\mathbb T^N}
u_i(t_0,x)\right)\boldsymbol1.
\]
Then $\boldsymbol v^-(0)\leq\boldsymbol u(t_0)$, and comparison with
\eqref{eq:lower-cooperative-KPP} gives
$\boldsymbol v^-(t)\leq\boldsymbol u(t_0+t)$ for all $t\geq0$.
Because $\lambda'>0$, Lemma~\ref{lem:cooperative-KPP-package}(i) yields
\[
\liminf_{t\to\infty}\min_i\min_{x\in\mathbb T^N}v_i^-(t,x)
\geq\eta_{\rm per}(p,g_*\boldsymbol1)>0.
\]
Thus one may take $\eta_{\rm per}=\eta_{\rm per}(p,g_*\boldsymbol1)$, which is
independent of $\boldsymbol u_0$.

\smallskip
\noindent(iii)
By the eventual bound in Lemma~\ref{lem:KPP-wellposed-bounded},
there is $T=T(\boldsymbol u_0)$ such that
$\|\boldsymbol u(t)\|_\infty\leq M_*+1$ for every $t\geq T$. Set
$t_0=T+1$ and $g_*=m\bar qK_{M_*+1}$. Since
$\boldsymbol u(t_0,x)\gg\boldsymbol0$ everywhere, choose a nonzero
compactly supported continuous function $\boldsymbol w_0$ satisfying
\[
\boldsymbol0\leq\boldsymbol w_0(x)\leq\boldsymbol u(t_0,x),
\qquad x\in\mathbb R^N.
\]
Let $\boldsymbol v^-$ solve
\[
\partial_t\boldsymbol v^-
=\mathcal L\boldsymbol v^--g_*(\boldsymbol v^-)^{\circ(1+p)},
\qquad \boldsymbol v^-(0)=\boldsymbol w_0.
\]
Comparison with \eqref{eq:lower-cooperative-KPP} gives
$\boldsymbol v^-(t)\leq\boldsymbol u(t_0+t)$ for every $t\geq0$.
Since $\lambda^*>0$, Lemma~\ref{lem:cooperative-KPP-package}(ii) gives, for
every $R>0$,
\[
\liminf_{t\to\infty}\min_i\inf_{|x|\leq R}u_i(t,x)
\geq
\eta_{\rm per}(p,g_*\boldsymbol1)>0.
\]
Thus one may take $\eta_{\rm loc}=\eta_{\rm per}(p,g_*\boldsymbol1)$,
independently of $R$ and $\boldsymbol u_0$.
\end{proof}

\begin{proof}[\textbf{Proof of Theorem~\ref{thm:nonlocal-tail-criteria}}]
(i)
Let $M_*$ be the constant in Lemma~\ref{lem:KPP-wellposed-bounded}. Choose $T=T(\boldsymbol u_0)$ with
\[
\|\boldsymbol u(t)\|_\infty\leq M_*+1,
\qquad t\geq T,
\]
and set $t_0=T+1$ and $g_*=m\bar qK_{M_*+1}$. We first show that the solution at time $t_0$ satisfies the required exponential lower bound. The convolution term and the off-diagonal terms of $A$ are nonnegative.
On $[0,t_0]$, the uniform-in-time bound in
Lemma~\ref{lem:KPP-wellposed-bounded} bounds all remaining diagonal
coefficients below by a finite constant $-d_{t_0}$. Hence, componentwise,
\[
\partial_tu_i(t,x)
\geq- d_{t_0}u_i(t,x).
\]
Scalar comparison at each fixed $x$ therefore yields
\[
u_i(t_0,x)\geq e^{-d_{t_0}t_0}u_{0,i}(x)
\geq C^{-1}e^{-d_{t_0}t_0-\xi\cdot x}
\quad\hbox{when }\xi\cdot x\geq B.
\]

Let $\boldsymbol v^-$ be the solution of
\[
\partial_t\boldsymbol v^-=
\mathcal L\boldsymbol v^--g_*(\boldsymbol v^-)^{\circ(1+p)},
\qquad
\boldsymbol v^-(0)=\boldsymbol u(t_0).
\]
The component comparison \eqref{eq:nonlinear-component-comparison},
applied from time $T$ with initial bound $M_*+1$, is valid for all
$t\geq t_0=T+1$. Hence
\eqref{eq:competition-power-control} holds with $g_M=g_*$, and
comparison gives
$\boldsymbol v^-(t)\leq\boldsymbol u(t_0+t)$.
The above exponential lower bound and
Lemma~\ref{lem:cooperative-KPP-package}(iii) imply
\[
\liminf_{t\to\infty}\min_i\inf_{|x|\leq R}v_i^-(t,x)
\geq\eta_{\rm per}(p,g_*\boldsymbol1)>0.
\]
Thus one may take $\eta_{\rm tail}=\eta_{\rm per}(p,g_*\boldsymbol1)$,
which is independent of $\xi$, $R$ and $\boldsymbol u_0$.

\smallskip
\noindent(ii)
Let $\boldsymbol v^+$ solve \eqref{eq:upper-cooperative-KPP}
with $\boldsymbol v^+(0)=\boldsymbol u_0$.
By \eqref{eq:competition-quadratic-control} and comparison,
$\boldsymbol0\leq\boldsymbol u(t)\leq\boldsymbol v^+(t)$.
Lemma~\ref{lem:cooperative-KPP-package}(iii), applied with
$q=1$ and $\boldsymbol b=q_0\boldsymbol1$, gives local uniform
convergence of $\boldsymbol v^+$ to zero, and hence the assertion.
The weighted-space comparison needed for the unbounded exponential
upper function is contained in \eqref{eq:weighted-linear-comparison}.
\end{proof}

\subsubsection{Spreading}

We first verify the geometric formula and then prove the upper and lower
spreading bounds using the cooperative comparison systems.

\begin{proof}[Proof of Proposition~\ref{prop:spreading-geometry}]
Each defining half-space is closed and convex, so the same
holds for $\mathcal W$. The inequality $re\cdot e'\leq c^*(e')$
is automatic when $e\cdot e'\leq0$, since $r\geq0$ and the
planar speeds are positive or infinite. For $e\cdot e'>0$ it is
equivalent to $r\leq c^*(e')/(e\cdot e')$. Taking the intersection
of these constraints gives the asserted radial formula, which is
precisely the Freidlin--G\"artner expression
\eqref{eq:planar-and-FG-speeds}.
\end{proof}

\begin{proof}[\textbf{Proof of Theorem~\ref{thm:nonlocal-spreading}}]
Fix $e\in\mathbb S^{N-1}$ and $R>0$. If $c>c_{\rm FG}(e)$, the
definitions in \eqref{eq:planar-and-FG-speeds} give
$e'\in\mathbb S^{N-1}$ and $\rho>0$ such that
\[
e\cdot e'>0,\qquad \mu:=\rho e'\in\mathcal I,\qquad
\lambda_L^*(\mu)<\rho c(e\cdot e').
\]
Choose $\gamma$ between the last two quantities. By
\eqref{def:lambdaLstar}, applied to the weighted operator \eqref{Q2},
there is $\boldsymbol\phi\in X^{++}$ satisfying
$\mathcal L_\mu\boldsymbol\phi\leq\gamma\boldsymbol\phi$. Since
$\boldsymbol u_0$ has compact support, choose $C>0$ such that
\[
\boldsymbol u_0(x)\leq C e^{-\mu\cdot x}\boldsymbol\phi(x).
\]
The function $\boldsymbol W(t,x)=C e^{\gamma t-\mu\cdot x}\boldsymbol\phi(x)$
is a super-solution of the linear system in $E_\mu$, the weighted space
defined in \eqref{eq:exponentially-weighted-space} with $\xi$
replaced by $\mu$. Compact support gives
$\boldsymbol u_0\in C_{\rm b}\cap E_\mu$. Applying the argument of
\eqref{eq:weighted-linear-comparison}, with the nonnegative competition
term as the loss term, gives
$\boldsymbol u(t)\leq S(t)\boldsymbol u_0\leq\boldsymbol W(t)$. Hence
\[
\max_i\sup_{|x|\leq R}u_i(t,x+cte)
\leq C\|\boldsymbol\phi\|_\infty e^{\rho R}
e^{-[\rho c(e\cdot e')-\gamma]t}\longrightarrow0,
\]
proving \eqref{eq:nonlocal-spreading-outside}.

Choose $T$ using the eventual bound in Lemma~\ref{lem:KPP-wellposed-bounded},
put $t_0=T+1$ and
$g_*=m\bar qK_{M_*+1}$, and choose a nonzero compactly supported
$\boldsymbol w_0\leq\boldsymbol u(t_0)$.  Let $\boldsymbol v^-$ solve
\eqref{eq:lower-cooperative-KPP} with $g_M=g_*$ and initial function
$\boldsymbol w_0$.  Then
$\boldsymbol v^-(t-t_0)\leq\boldsymbol u(t)$ for $t\geq t_0$.
Use the cutoffs $\eta_n$ from the proof of
Theorem~\ref{thm:large-period-limits}(i) and set
$J_{i,n}=\eta_nJ_i$, retaining the loss term $-I$. Denote the corresponding
periodic generalized eigenvalue, planar speed and ray speed by
$\lambda_{L,n}^*(\mu)$, $c_n^*(e')$ and $c_{{\rm FG},n}(e)$.
By \eqref{eq:monotone-kernel-truncation} and the moment lower bound used
in the proof of Theorem~\ref{TH2.1},
\[
\lambda_{L,n}^*(\mu)\uparrow\lambda_L^*(\mu),
\qquad \mu\in\mathbb R^N,
\]
including the value $+\infty$ outside $\mathcal I$. The proof of
Theorem~\ref{TH2.1}(ii) also gives
\[
\min_\mu\lambda_{L,n}^*(\mu)\uparrow\lambda^*>0.
\]
Positivity of the retained kernels on a fixed ball gives constants
$\delta,a,b,d>0$ and $n_0$ such that
\begin{equation}\label{eq:truncated-spreading-coercivity}
\lambda_{L,n}^*(\mu)\geq\delta,\qquad
\lambda_{L,n}^*(\mu)\geq a e^{b|\mu|}-d,
\qquad n\geq n_0,\quad \mu\in\mathbb R^N.
\end{equation}

For fixed $e'$, the speeds $c_n^*(e')$ increase and do not exceed
$c^*(e')$. If their limit is finite, \eqref{eq:truncated-spreading-coercivity} confines their
minimizing weights $\rho_n$ to a compact subinterval of $(0,\infty)$.
Along a subsequence, $\rho_n\to\rho_0>0$. For every fixed $k\geq n_0$,
\[
\lim_{n\to\infty}c_n^*(e')
\geq\frac{\lambda_{L,k}^*(\rho_0e')}{\rho_0}.
\]
Letting $k\to\infty$ gives the reverse inequality. Thus
\[
c_n^*(e')\uparrow c^*(e'),
\qquad e'\in\mathbb S^{N-1},
\]
also when $c^*(e')=+\infty$. The bounds \eqref{eq:truncated-spreading-coercivity} also show that
$c_n^*$ is continuous on $\mathbb S^{N-1}$ for every $n\geq n_0$, and
that $c_n^*(e')\geq c_0>0$ uniformly in $n\geq n_0$ and $e'$.

Fix $0<c<c_{\rm FG}(e)$ and choose $c<c'<c_{\rm FG}(e)$.
Take $0<\alpha<\min\{1,c_0/c'\}$ and set
\[
K_\alpha:=\{e'\in\mathbb S^{N-1}:e\cdot e'\geq\alpha\}.
\]
On this compact set,
\[
\frac{c_n^*(e')}{e\cdot e'}
\uparrow\frac{c^*(e')}{e\cdot e'}\geq c_{\rm FG}(e)>c'.
\]
The increasing open sets
$\{e'\in K_\alpha:c_n^*(e')/(e\cdot e')>c'\}$ cover $K_\alpha$.
Compactness therefore gives a fixed $n$ for which the inequality holds
throughout $K_\alpha$. When $0<e\cdot e'<\alpha$, the quotient is at least
$c_0/\alpha>c'$. Consequently $c_{{\rm FG},n}(e)>c'$.

Let $\boldsymbol z_n$ solve the lower cooperative system with kernels
$J_{i,n}$ and initial function $\boldsymbol w_0$. The comparison \eqref{eq:cooperative-discrete-lower-comparison}
and the spectral-radius bound \eqref{eq:cooperative-discrete-spectral-gap}
for the maps in \eqref{eq:cooperative-discrete-comparison-maps} do not require unit mass,
and $\int_{\mathbb R^N}J_{i,n}(z)\,\mathrm dz\leq1$ preserves the constant
super-solution estimate. Thus Step~5 of the proof of
Lemma~\ref{lem:cooperative-KPP-package} applies to these compactly
supported kernels with their positive masses retained in the convolution
terms. Since $\min_\mu\lambda_{L,n}^*(\mu)>0$, the
inside-spreading result applies to $\boldsymbol z_n$. Kernel ordering
gives $\boldsymbol z_n\leq\boldsymbol v^-$. Hence, for every fixed
$R_1>0$, there is $\beta>0$ such that
\[
v_i^-(s,x+sc'e)\geq\beta,
\qquad |x|\leq R_1,\quad i\in\Lambda,
\]
for all sufficiently large $s$. Use the fixed cube $\Omega$ and the
eigenpair $(\lambda_\Omega,\boldsymbol\phi)$ from Step~4 of the proof of
Lemma~\ref{lem:cooperative-KPP-package}, with $q=p$ and
$\boldsymbol b=g_*\boldsymbol1$. Set
\[
a_\infty=(\lambda_\Omega/g_*)^{1/p},\qquad
a'=\lambda_\Omega a-g_*a^{1+p},\qquad
a(0)=\min\{\beta,a_\infty/2\}.
\]
Choose $R_1>R+\operatorname{diam}\mathcal C+\sup_{y\in\Omega}|y|$.
At time $s_t=(c/c')t$, the preceding lower bound holds on
$B_{R_1}+tce$. Every cell $\ell+\overline{\mathcal C}$ meeting
$B_R+tce$, with $\ell\in\mathbb L$, has
$\Omega+\ell\subset B_{R_1}+tce$. Hence bounded-domain comparison gives,
for all sufficiently large $t$,
\[
u_i(t,x)\geq v_i^-(t-t_0,x)
\geq a(t-t_0-s_t)\phi_i(x-\ell),
\qquad x\in\ell+\overline{\mathcal C}.
\]
Since $t-t_0-s_t\to\infty$ and $a(s)\to a_\infty$, this gives a positive
lower bound independent of $R$, $c$, $e$ and $\boldsymbol u_0$.
For any fixed $0<c<c_{\rm FG}(e)$, the same argument, with no initial time
shift, applies to $\boldsymbol v^-$. Step~2 of the proof of
Lemma~\ref{lem:cooperative-KPP-package} therefore upgrades its lower
bound on the moving balls to $\eta_{\rm per}(p,g_*\boldsymbol1)$.
Using radius $R+ct_0$ for $\boldsymbol v^-$ absorbs the fixed time shift
in the comparison and gives
\[
\eta_{\rm sp}=\eta_{\rm per}(p,g_*\boldsymbol1)>0.
\]
Thus $\eta_{\rm per}$, $\eta_{\rm loc}$, $\eta_{\rm tail}$ and
$\eta_{\rm sp}$ in the three theorems can all be taken to be the minimum
of the positive periodic equilibrium of the same lower cooperative system
with $g_M=g_*$.

For $e\in\mathbb S^{N-1}$, the absence of positive exponential moments
in direction $e$ is expressed by
\begin{equation}
\label{eq:no-directional-exponential-moment}
\int_{\mathbb R^N}J_{i_0}(y)e^{\rho e\cdot y}\,\mathrm dy=+\infty
\quad\text{for some }i_0\text{ and every }\rho>0.
\end{equation}

By \eqref{eq:planar-and-FG-speeds}, $c_{\rm FG}(e)<+\infty$ if and only
if $c^*(e')<+\infty$ for some $e'$ with $e\cdot e'>0$. The latter is
equivalent to $\rho e'\in\mathcal I$ for some $\rho>0$, or, equivalently,
to the existence of $\mu\in\mathcal I$ with $\mu\cdot e>0$.
Finally, the first formula in \eqref{eq:planar-and-FG-speeds} gives
$c^*(e)=+\infty$ if and only if $\rho e\notin\mathcal I$ for every
$\rho>0$. Condition \eqref{eq:no-directional-exponential-moment} implies
this exclusion. Conversely, if it fails, then for each $i\in\Lambda$
there is $\rho_i>0$ such that
\[
\int_{\mathbb R^N}J_i(y)e^{\rho_i e\cdot y}\,\mathrm dy<\infty.
\]
Set $\rho_0=\min_i\rho_i>0$. H\"older's inequality and
$\int_{\mathbb R^N}J_i(y)\,\mathrm dy=1$ give
\[
\int_{\mathbb R^N}J_i(y)e^{\rho_0 e\cdot y}\,\mathrm dy
\leq\left(\int_{\mathbb R^N}J_i(y)e^{\rho_i e\cdot y}\,\mathrm dy
\right)^{\rho_0/\rho_i}<\infty,
\qquad i\in\Lambda.
\]
Thus $\rho_0e\in\mathcal I$, proving the planar equivalence. The second
formula in \eqref{eq:planar-and-FG-speeds} gives
$c_{\rm FG}(e)=+\infty$ if and only if $c^*(e')=+\infty$ for every
$e'\in\mathbb S^{N-1}$ with $e\cdot e'>0$. Applying the planar
equivalence to each such $e'$ proves the stated ray-speed criterion.
\end{proof}

\subsubsection{Period limits}

We first combine the variational speed formulas with the spectral limits
of Sections~\ref{sec:large-period-asymptotics} and~\ref{sec:quantitative-period}.

\begin{proof}[\textbf{Proof of Theorem~\ref{cor:speed-period-limits}}]
For $e\in\mathbb S^{N-1}$, set
\begin{equation*}
F_L(\rho):=\frac{\lambda_L^*(\rho e)}{\rho},
\qquad
F_\infty(\rho):=
\frac{\overline\lambda_\infty(\rho e)}{\rho}
\quad(\rho>0,\ \rho e\in\mathcal I).
\end{equation*}
For $\rho e\notin\mathcal I$, set $F_L(\rho)=F_\infty(\rho)=+\infty$.
By Theorem~\ref{thm:large-period-limits}(i), $F_L(\rho)$ converges to
$F_\infty(\rho)$ as $L\to\infty$ for every $\rho>0$ with
$\rho e\in\mathcal I$. The convergence is locally uniform on
$\{\rho>0:\rho e\in\operatorname{int}\mathcal I\}$.
For the planar limits, we may assume that the admissible ray is nonempty,
since otherwise all corresponding speeds equal $+\infty$.

Assume first that $\min_\mu\overline\lambda_\infty(\mu)>0$.
Theorem~\ref{thm:large-period-limits} gives, for some
$L_0,\delta>0$,
\begin{equation}\label{eq:speed-positive-coercivity-zero}
\lambda_L^*(\rho e)\geq\delta
\quad(L\geq L_0),
\qquad
F_L(\rho)\geq\frac{\delta}{\rho}.
\end{equation}
The lower bound \eqref{eq:weighted-moment-lower-bound} and positivity
of the kernels on a fixed ball give
$\lambda_L^*(\rho e)\geq C_1e^{C_2\rho}-C_3$, with positive constants
independent of $L$. Consequently,
\begin{equation}\label{eq:speed-positive-coercivity-infty}
\lim_{R\to\infty}\inf_{\substack{L\geq L_0\\\rho\geq R}}F_L(\rho)=+\infty.
\end{equation}
Fixing an admissible weight gives a uniform upper bound for
$c_L^*(e)$. The bounds \eqref{eq:speed-positive-coercivity-zero} and
\eqref{eq:speed-positive-coercivity-infty}, together with this upper
bound, confine nearly minimizing weights to a common compact subinterval of $(0,\infty)$.
Choose $L_k\to\infty$ along which $c_{L_k}^*(e)$ converges to its
lower limit, and choose $\rho_k>0$ such that
$F_{L_k}(\rho_k)\leq c_{L_k}^*(e)+1/k$. After passage to a subsequence,
$\rho_k\to\rho_0>0$. Kernel ordering and the locally uniform
compact-kernel limit give, for each fixed $n$,
\[
\liminf_{L\to\infty}c_L^*(e)
=\lim_{k\to\infty}F_{L_k}(\rho_k)
\geq\frac{\overline\lambda_{\infty,n}(\rho_0e)}{\rho_0}.
\]
Letting $n\to\infty$ in
\eqref{eq:effective-constant-truncation-all-weights}, and testing
each fixed admissible $\rho>0$, respectively, yields
\[
\liminf_{L\to\infty}c_L^*(e)\geq F_\infty(\rho_0)
\geq\inf_{\rho>0}F_\infty(\rho),
\qquad
\limsup_{L\to\infty}c_L^*(e)\leq F_\infty(\rho).
\]
Taking the infimum in the last inequality gives
\begin{equation*}
\lim_{L\to\infty}c_L^*(e)
=\lim_{L\to\infty}\inf_{\rho>0}F_L(\rho)
=\inf_{\rho>0}F_\infty(\rho)=c_\infty^*(e).
\end{equation*}

For $L\to0$, assume $\min_\mu\widehat\lambda(\mu)>0$.
Theorem~\ref{thm:small-period-limits} and
\eqref{eq:weighted-moment-lower-bound} give the same uniform bounds,
so nearly minimizing weights again lie in a fixed compact
subinterval of $(0,\infty)$.
Choose $L_k\downarrow0$ along which $c_{L_k}^*(e)$ converges to its
lower limit, and choose $\rho_k>0$ such that
$F_{L_k}(\rho_k)\leq c_{L_k}^*(e)+1/k$. After passage to a subsequence,
$\rho_k\to\rho_0>0$. Denote by $\widehat\lambda_n$ the averaged spectral
quantity associated with the truncated kernels $J_{i,n}$.
For each fixed $n$, the compact-kernel averaging limit is locally
uniform on $\mathbb R^N$. Kernel ordering therefore gives
\[
\liminf_{k\to\infty}F_{L_k}(\rho_k)
\geq\frac{\widehat\lambda_n(\rho_0e)}{\rho_0}.
\]
By the proof of Theorem~\ref{thm:small-period-limits},
$\widehat\lambda_n\uparrow\widehat\lambda$ on $\mathbb R^N$,
including the value $+\infty$ outside $\mathcal I$. Thus
\[
\liminf_{L\to0}c_L^*(e)
\geq\frac{\widehat\lambda(\rho_0e)}{\rho_0}
\geq\widehat c_0^*(e).
\]
Conversely, for every admissible $\rho>0$,
\[
\limsup_{L\to0}c_L^*(e)
\leq\lim_{L\to0}F_L(\rho)
=\frac{\widehat\lambda(\rho e)}{\rho}.
\]
Taking the infimum over $\rho$ gives
\begin{equation*}
\lim_{L\to0}c_L^*(e)
=\inf_{\substack{\rho>0\\\rho e\in\mathcal I}}
\frac{\widehat\lambda(\rho e)}{\rho}=\widehat c_0^*(e).
\end{equation*}

For the ray-speed limits, let $\ell=\overline\lambda_\infty$ when
$L\to\infty$ and $\ell=\widehat\lambda$ when $L\to0$, and assume
$\min_\mu\ell(\mu)>0$. Set
$D_e:=\{\mu\in\mathcal I:\mu\cdot e>0\}$. Writing $\mu=\rho e'$ in
\eqref{eq:planar-and-FG-speeds} gives
\[
c_{{\rm FG},L}(e)
=\inf_{\mu\in D_e}\frac{\lambda_L^*(\mu)}{\mu\cdot e}.
\]
If $D_e=\varnothing$, all the corresponding ray speeds equal
$+\infty$. Otherwise, fixing
one weight in $D_e$ and using its pointwise spectral limit gives
$c_{{\rm FG},L}(e)\leq C$ for $L$ sufficiently near the relevant limit.
Theorems~\ref{thm:large-period-limits} and~\ref{thm:small-period-limits},
\eqref{eq:weighted-moment-lower-bound}, and positivity of the kernels
on a fixed ball give
$\delta,a,b,d>0$ such that, for all such $L$,
\[
\lambda_L^*(\mu)\geq\delta,
\qquad
\lambda_L^*(\mu)\geq ae^{b|\mu|}-d,
\qquad \mu\in\mathbb R^N.
\]
Choose $L_k$ approaching the relevant limit along which
$c_{{\rm FG},L_k}(e)$ approaches its lower limit, and choose
$\mu_k\in D_e$ with
\[
\frac{\lambda_{L_k}^*(\mu_k)}{\mu_k\cdot e}
\leq c_{{\rm FG},L_k}(e)+\frac1k\leq C+1.
\]
The preceding bounds imply that $\mu_k$ is bounded and
$\mu_k\cdot e\geq\delta/(C+1)$. After passage to a subsequence,
$\mu_k\to\mu_0$ with $\mu_0\cdot e>0$.
Let $\ell_n$ be the corresponding limit for the truncated kernels
$J_{i,n}$. Their spectral limits are locally uniform on $\mathbb R^N$,
and the preceding truncation arguments give $\ell_n\uparrow\ell$ on
$\mathbb R^N$, including the value $+\infty$ outside $\mathcal I$.
Kernel ordering therefore yields, for every fixed $n$,
\[
\liminf_{k\to\infty}c_{{\rm FG},L_k}(e)
\geq\frac{\ell_n(\mu_0)}{\mu_0\cdot e}.
\]
Letting $n\to\infty$ forces $\mu_0\in D_e$ and gives the lower bound
\[
\liminf_L c_{{\rm FG},L}(e)
\geq\inf_{\mu\in D_e}\frac{\ell(\mu)}{\mu\cdot e}.
\]
For every fixed
$\mu\in D_e$, the inequality
$c_{{\rm FG},L}(e)\leq\lambda_L^*(\mu)/(\mu\cdot e)$ and pointwise
spectral convergence give
\[
\limsup_L c_{{\rm FG},L}(e)\leq\frac{\ell(\mu)}{\mu\cdot e}.
\]
Taking the infimum over $D_e$ proves the two ray-speed limits.

The condition that every kernel has finite exponential moments of every order is
\begin{equation}
\label{eq:all-exponential-moments}
\int_{\mathbb R^N}J_i(y)e^{\rho|y|}\,\mathrm dy<\infty,
\qquad \rho>0,\quad i=1,\ldots,m.
\end{equation}

Under \eqref{eq:all-exponential-moments} and the respective positivity
assumptions, the bounds \eqref{eq:speed-positive-coercivity-zero} and
\eqref{eq:speed-positive-coercivity-infty} confine the minimizing weights to a common compact subinterval
of $(0,\infty)$, uniformly in $e$. The locally uniform spectral limits give
\begin{align*}
\sup_{e\in\mathbb S^{N-1}}|c_L^*(e)-c_\infty^*(e)|&\longrightarrow0
&& (L\to\infty),\\
\sup_{e\in\mathbb S^{N-1}}|c_L^*(e)-\widehat c_0^*(e)|&\longrightarrow0
&& (L\to0).
\end{align*}
\end{proof}

\subsubsection{Small-period corrections}

We expand the minimizing weights and substitute them into the speed
formulas.

\begin{proof}[\textbf{Proof of Theorem~\ref{thm:Bessel-directional-speed-expansion}}]
By evenness of each kernel, $T_{-\rho}=T_\rho$.
Theorem~\ref{thm:Bessel-refined-expansion} therefore yields,
locally uniformly with two derivatives near $\rho_0$,
\[
 \frac{\lambda_L^*(\pm\rho)}{\rho}
 =F(\rho)+B(\rho)L^{2s}\pm C(\rho)L^{2s+1}
 +Q(\rho)L^{4s}+E(\rho)L^{2s+2}
 +O(L^{4s+1}).
\]
Local uniform spectral convergence near zero and
$\widehat\lambda(0)>0$ give a fixed neighborhood of zero on which
$\lambda_L^*(\pm\rho)\geq\widehat\lambda(0)/2$ for all small $L$.
Thus the speed quotients diverge uniformly as $\rho\downarrow0$.
At the right endpoint, the diagonal comparison
$\lambda_L^*(\pm\rho)\geq M_i(\rho)-C$ for the component with
$\sigma_i=\sigma_{\max}$ forces divergence uniformly in $L$.
The locally uniform small-period limit and the unique limiting
minimum therefore confine all minimizers near $\rho_0$.
Uniform $C^2$ convergence and $F''(\rho_0)>0$ imply strict
convexity there and hence uniqueness for small $L$.
Writing these minimizers as $\rho_L^\pm$, differentiation yields
\[
 \rho_L^\pm=\rho_0-
 \frac{B'(\rho_0)}{F''(\rho_0)}L^{2s}
 +O(L^{2s+1}).
\]
Substitution into the preceding expansion shows that the quadratic
variation of $F$ and linear variation of $B$ combine into
$-[B'(\rho_0)]^2L^{4s}/(2F''(\rho_0))$.
All other optimizer corrections have order at least $L^{4s+1}$;
this proves \eqref{eq:Bessel-refined-speeds}. The terms of orders
$L^{4s}$ and $L^{2s+2}$ are the same for both directions,
so subtraction proves \eqref{eq:Bessel-refined-speed-gap}.
For pointwise symmetric $A$ the averaged spectral operator
is self-adjoint, while $D_sU_s$ is skew-adjoint, and hence
$c_s(\rho)=0$. Theorem~\ref{thm:nonlocal-spreading}, together with
the equality of planar and ray speeds in one dimension, identifies
these spectral speeds with the nonlinear spreading speeds.
\end{proof}

\begin{proof}[Proof of Corollary~\ref{cor:Laplace-minimum-speed}]
The expansion and its first derivative from
Theorem~\ref{thm:Laplace-second-order} give
$\partial_\mu\lambda_L^*(\mu)
=\widehat\lambda'(\mu)+L^2\mathfrak b_2'(\mu)+O(L^3)$.
Taylor expansion at the nondegenerate minimizer and the mean-value
theorem yield the first formula for $\mu_L$; substitution yields
the spectral minimum. For the speed, minimize
$\lambda_L^*(\rho)/\rho$ on $(0,1)$. Local uniform convergence near
zero and $\widehat\lambda(0)>0$ give
$\lambda_L^*(\rho)\geq\widehat\lambda(0)/2$ in a fixed neighborhood
of zero for all small $L$. Thus the quotient diverges uniformly
as $\rho\downarrow0$.
Moreover, the constant positive-vector test and
$M(\rho)=(1-\rho^2)^{-1}\to\infty$ imply
$\lambda_L^*(\rho)\geq M(\rho)-C$ uniformly in $L$,
so the quotient diverges as $\rho\uparrow1$.
The uniform small-period convergence on compact subsets and the
uniqueness of $\rho_0$ localize the minimizers near $\rho_0$.
Apply the same derivative argument to
$f_L(\rho)=\lambda_L^*(\rho)/\rho$.
The last qualification follows from the Taylor remainder.
\end{proof}

\section*{Acknowledgments}
M.-Z. Xin is partially supported by the National Natural
Science Foundation of China (No. 12601334), the China National Postdoctoral Program for
Innovative Talents (BX2026009), and the China Postdoctoral Science Foundation (2026M793394).


\begin{thebibliography}{99}

\bibitem{Bao2020ZAMP}
X. Bao,
Stability of periodic traveling waves for nonlocal dispersal cooperative
systems in space-time periodic habitats,
\emph{Z. Angew. Math. Phys.} 71 (2020), Paper No. 169.

\bibitem{BS} X. Bao and W. Shen, Criteria for the existence of principal eigenvalues of time periodic cooperative linear systems with nonlocal dispersal, \emph{Proc. Amer. Math. Soc.} 145 (2017), 2881--2894.

\bibitem{BaoShenShen2019}
X. Bao, W. Shen, and Z. Shen,
Spreading speeds and traveling waves for space-time periodic nonlocal
dispersal cooperative systems,
\emph{Commun. Pure Appl. Anal.} 18 (2019), 361--396.

\bibitem{BarlesPerthame1987} G. Barles and B. Perthame, Discontinuous solutions of deterministic optimal stopping time problems, \emph{RAIRO Modél. Math. Anal. Numér.} 21 (1987), 557--579.

\bibitem{BentoutVo2025}
S. Bentout and H.-H. Vo,
The eigentheory for nonlocal cooperative-advective system and its role
in the study of free boundary system for directional epidemic models,
preprint, arXiv:2510.10024 (2025).

\bibitem{BentoutVo2026}
S. Bentout and H.-H. Vo,
Long-time dynamics and threshold phenomena for a free-boundary SIS model
with asymmetric kernels in advective periodic environments,
preprint, arXiv:2603.24368 (2026).

\bibitem{Berestycki2016} H. Berestycki, J. Coville, and H.-H. Vo, On the definition and the properties of the principal eigenvalue of some nonlocal operators, \emph{J. Funct. Anal.} 271 (2016), 2701--2751.

\bibitem{BCV2016Persistence}
H. Berestycki, J. Coville, and H.-H. Vo,
Persistence criteria for populations with non-local dispersion,
\emph{J. Math. Biol.} 72 (2016), 1693--1745.

\bibitem{BHN2005}
H. Berestycki, F. Hamel, and N. Nadirashvili,
The speed of propagation for KPP type problems. I: Periodic framework,
\emph{J. Eur. Math. Soc.} 7 (2005), 173--213.

\bibitem{BerestyckiRossi2015}
H. Berestycki and L. Rossi,
Generalizations and properties of the principal eigenvalue of elliptic
operators in unbounded domains,
\emph{Comm. Pure Appl. Math.} 68 (2015), 1014--1065.

\bibitem{Coville2010}
J. Coville,
On a simple criterion for the existence of a principal eigenfunction
of some nonlocal operators,
\emph{J. Differential Equations} 249 (2010), 2921--2953.

\bibitem{Coville2012Harnack}
J. Coville,
Harnack type inequality for positive solution of some integral equation,
\emph{Ann. Mat. Pura Appl.} 191 (2012), 503--528.

\bibitem{C2}
J. Coville,
Singular measure as principal eigenfunction of some nonlocal operators,
\emph{Appl. Math. Lett.} 26 (2013), 831--835.

\bibitem{CovilleDavilaMartinez2008}
J. Coville, J. D\'avila, and S. Mart\'inez,
Existence and uniqueness of solutions to a nonlocal equation with
monostable nonlinearity,
\emph{SIAM J. Math. Anal.} 39 (2008), 1693--1709.

\bibitem{CovilleDavilaMartinez2013}
J. Coville, J. D\'avila, and S. Mart\'inez,
Pulsating fronts for nonlocal dispersion and KPP nonlinearity,
\emph{Ann. Inst. H. Poincar\'e Anal. Non Lin\'eaire}
30 (2013), 179--223.

\bibitem{CrandallIshiiLions1992}
M. G. Crandall, H. Ishii, and P.-L. Lions,
User's guide to viscosity solutions of second order partial differential equations,
\emph{Bull. Amer. Math. Soc. (N.S.)} 27 (1992), 1--67.

\bibitem{CrandallLions1983}
M. G. Crandall and P.-L. Lions,
Viscosity solutions of Hamilton--Jacobi equations,
\emph{Trans. Amer. Math. Soc.} 277 (1983), 1--42.

\bibitem{Deimling1985}
K. Deimling,
\emph{Nonlinear Functional Analysis},
Springer, Berlin, 1985.

\bibitem{DL2015} W. Ding and X. Liang, Principal eigenvalues of generalized convolution operators on the circle and spreading speeds of noncompact evolution systems in periodic media, \emph{SIAM J. Math. Anal.} 47 (2015), 855--896.

\bibitem{DuFangNi2025}
Y. Du, X. Fang, and W. Ni,
Asymptotic limit of the principal eigenvalue of asymmetric nonlocal
diffusion operators and propagation dynamics,
\emph{SIAM J. Math. Anal.}, to appear, arXiv:2503.22062.

\bibitem{DuLiShen2022}
L.-J. Du, W.-T. Li, and W. Shen,
Propagation phenomena for time-space periodic monotone semiflows and
applications to cooperative systems in multi-dimensional media,
\emph{J. Funct. Anal.} 282 (2022), Paper No. 109415.

\bibitem{DuLiXin2026}
L.-J. Du, W.-T. Li, and M.-Z. Xin,
Uniqueness and stability of monostable pulsating fronts for
multi-dimensional reaction-diffusion-advection systems in periodic media,
\emph{J. Math. Pures Appl.} 208 (2026), Paper No. 103843.

\bibitem{ElliottCornell2012}
E. C. Elliott and S. J. Cornell,
Dispersal polymorphism and the speed of biological invasions,
\emph{PLoS ONE} 7 (2012), Paper No. e40496.

\bibitem{Feng2025} Y.-X. Feng, W.-T. Li, Y. Lou, and F.-Y. Yang, Asymptotic behavior of the basic reproduction number for periodic nonlocal dispersal operators and applications, \emph{J. Math. Biol.} 90 (2025), Paper No. 24.

\bibitem{Feng2024} Y.-X. Feng, W.-T. Li, S. Ruan, and M.-Z. Xin, Principal spectral theory of time-periodic nonlocal dispersal cooperative systems and applications, \emph{SIAM J. Math. Anal.} 56 (2024), 4040--4083.

\bibitem{Garnier2011}
J. Garnier,
Accelerating solutions in integro-differential equations,
\emph{SIAM J. Math. Anal.} 43 (2011), 1955--1974.

\bibitem{GartnerFreidlin1979}
J. G\"artner and M. I. Freidlin,
On the propagation of concentration waves in periodic and random media,
\emph{Soviet Math. Dokl.} 20 (1979), 1282--1286.

\bibitem{Girardin2025}
L. Girardin,
Persistence, extinction, and spreading properties of noncooperative
Fisher--KPP systems in space-time periodic media,
\emph{SIAM J. Math. Anal.} 57 (2025), 233--261.

\bibitem{Griette2019}
Q. Griette,
Singular measure traveling waves in an epidemiological model with
continuous phenotypes,
\emph{Trans. Amer. Math. Soc.} 371 (2019), 4411--4458.

\bibitem{GuoLiLiu2026}
H. Guo, W.-T. Li, and B. Liu,
Propagation dynamics for multidimensional nonlocal diffusion equations:
a general Freidlin--G\"artner formula,
preprint, arXiv:2606.13050 (2026).

\bibitem{HamelNadinRoques2011}
F. Hamel, G. Nadin, and L. Roques,
A viscosity solution method for the spreading speed formula in slowly
varying media,
\emph{Indiana Univ. Math. J.} 60 (2011), 1229--1248.

\bibitem{HutsonMartinezMischaikowVickers2003}
V. Hutson, S. Martinez, K. Mischaikow, and G. T. Vickers,
The evolution of dispersal,
\emph{J. Math. Biol.} 47 (2003), 483--517.

\bibitem{KangLiu2025}
H. Kang and S. Liu,
Global dynamics and asymptotic spreading of a diffusive age-structured
model in spatially periodic media,
\emph{Math. Ann.} 393 (2025), 2631--2682.

\bibitem{KangLiu2026} H. Kang and S. Liu, Spreading speeds of diffusive age-structured equations in slowly and rapidly varying media, \emph{J. Differential Equations} 480 (2026), Paper No. 114669.

\bibitem{KLS} C.-Y. Kao, Y. Lou, and W. Shen, Random dispersal vs. nonlocal dispersal, \emph{Discrete Contin. Dyn. Syst.} 26 (2010), 551--596.

\bibitem{Kato1995}
T. Kato,
\emph{Perturbation Theory for Linear Operators},
Classics in Mathematics, Springer, Berlin, 1995.

\bibitem{KotLewisVanDenDriessche1996}
M. Kot, M. A. Lewis, and P. van den Driessche,
Dispersal data and the spread of invading organisms,
\emph{Ecology} 77 (1996), 2027--2042.

\bibitem{LamLou2016}
K.-Y. Lam and Y. Lou,
Asymptotic behavior of the principal eigenvalue for cooperative elliptic
systems and applications,
\emph{J. Dynam. Differential Equations} 28 (2016), 29--48.

\bibitem{LiLiWang2022}
L. Li, W.-T. Li, and M.-X. Wang,
Dynamics for nonlocal diffusion problems with a free boundary,
\emph{J. Differential Equations} 330 (2022), 110--149.

\bibitem{LiWang2025}
L. Li and M.-X. Wang,
Dynamics of a nonlocal diffusive epidemic model with a free boundary,
\emph{Discrete Contin. Dyn. Syst.} 47 (2026), 87--119.

\bibitem{LiXinZhao2025}
W.-T. Li, M.-Z. Xin, and X.-Q. Zhao,
Spatio-temporal dynamics of nonlocal dispersal systems in time-space
periodic habitats,
\emph{J. Differential Equations} 416 (2025), 2000--2042.

\bibitem{Liang2019} X. Liang, L. Zhang, and X.-Q. Zhao, The principal eigenvalue for periodic nonlocal dispersal systems with time delay, \emph{J. Differential Equations} 266 (2019), 2100--2124.

\bibitem{LiangZhao2010}
X. Liang and X.-Q. Zhao,
Spreading speeds and traveling waves for abstract monostable evolution
systems,
\emph{J. Funct. Anal.} 259 (2010), 857--903.

\bibitem{LionsPapanicolaouVaradhan1988}
P.-L. Lions, G. Papanicolaou, and S. R. S. Varadhan,
Homogenization of Hamilton--Jacobi equations,
unpublished manuscript (1988).

\bibitem{Liu2025}
S. Liu,
Asymptotic and global analysis of principal eigenvalues for linear
time-periodic parabolic systems,
\emph{J. Math. Pures Appl.} 203 (2025), Paper No. 103781.

\bibitem{LiuLou2020}
S. Liu and Y. Lou,
A functional approach towards eigenvalue problems associated with
incompressible flow,
\emph{Discrete Contin. Dyn. Syst.} 40 (2020), 3715--3736.

\bibitem{LiuLou2022}
S. Liu and Y. Lou,
Classifying the level set of principal eigenvalue for time-periodic
parabolic operators and applications,
\emph{J. Funct. Anal.} 282 (2022), Paper No. 109338.

\bibitem{LiuLou2024}
S. Liu and Y. Lou,
Monotonicity, asymptotics and level sets for principal eigenvalues
of some elliptic operators with shear flow,
\emph{J. Math. Pures Appl.} 191 (2024), Paper No. 103622.

\bibitem{LiuLou2027}
S. Liu and Y. Lou,
Principal eigenvalues for non-self-adjoint operators: Theory and applications,
\emph{Discrete Contin. Dyn. Syst.} 60 (2027), Paper No. 105528.

\bibitem{Nadin2009Eigenvalue}
G. Nadin,
The principal eigenvalue of a space-time periodic parabolic operator,
\emph{Ann. Mat. Pura Appl.} 188 (2009), 269--295.

\bibitem{PR24} A. Piatnitski and V. Rybalko, Homogenization of nonlocal spectral problems, \emph{Z. Angew. Math. Phys.} 75 (2024), Paper No. 226.

\bibitem{Rossi2017}
L. Rossi,
The Freidlin--G\"artner formula for general reaction terms,
\emph{Adv. Math.} 317 (2017), 267--298.

\bibitem{Shen2015} W. Shen and X. Xie, On principal spectrum points/principal eigenvalues of nonlocal dispersal operators and applications, \emph{Discrete Contin. Dyn. Syst.} 35 (2015), 1665--1696.

\bibitem{SZ} W. Shen and A. Zhang, Spreading speeds for monostable equations with nonlocal dispersal in space periodic habitats, \emph{J. Differential Equations} 249 (2010), 747--795.

\bibitem{Shen2019} Z. Shen and H.-H. Vo, Nonlocal dispersal equations in time-periodic media: principal spectral theory, limiting properties and long-time dynamics, \emph{J. Differential Equations} 267 (2019), 1423--1466.

\bibitem{S20231} Y.-H. Su, W.-T. Li, Y. Lou, and X. Wang, Principal spectral theory for nonlocal systems and applications to stem cell regeneration models, \emph{J. Math. Pures Appl.} 176 (2023), 226--281.

\bibitem{Su2020a} Y.-H. Su, W.-T. Li, Y. Lou, and F.-Y. Yang, The generalised principal eigenvalue of time-periodic nonlocal dispersal operators and applications, \emph{J. Differential Equations} 269 (2020), 4960--4997.

\bibitem{Su2020b} Y.-H. Su, W.-T. Li, and F.-Y. Yang, Asymptotic behaviors for nonlocal diffusion equations about the dispersal spread, \emph{Anal. Appl. (Singap.)} 18 (2020), 585--614.

\bibitem{S20232} Y.-H. Su, X. Wang, and T. Zhang, Principal spectral theory and variational characterizations for cooperative systems with nonlocal and coupled diffusion, \emph{J. Differential Equations} 369 (2023), 94--114.

\bibitem{TaoLiRuanXu2026}
W. Tao, W.-T. Li, S. Ruan, and W.-B. Xu,
A complete classification of spreading speeds for nonlocal dispersal
Fisher--KPP equations in shifting habitats: a Hamilton--Jacobi approach,
\emph{Math. Ann.} 394 (2026), Paper No. 6.

\bibitem{Thu2025} N. V. Thu and H.-H. Vo, Asymptotic behavior of the generalized principal eigenvalue for an inhomogeneous cooperative system with nonlocal dispersal, \emph{Calc. Var. Partial Differential Equations} 64 (2025), Paper No. 125.

\bibitem{TrangVo2026Fractional}
C.-B. Trang and H.-H. Vo,
Spectral theory of fractional cooperative systems and threshold dynamics
in epidemic models, preprint, arXiv:2602.13986 (2026).

\bibitem{Vo2022Neumann}
H.-H. Vo,
Principal spectral theory of time-periodic nonlocal dispersal operators
of Neumann type,
\emph{Math. Nachr.} 295 (2022), 806--826.

\bibitem{Vo2026Advection}
H.-H. Vo,
Principal eigenvalues of nonlocal operators with advection: sharp
scaling limits and spectral phase transitions,
preprint, arXiv:2609.02416 (2026).

\bibitem{VoTrang2025}
H.-H. Vo and C.-B. Trang,
On the existence of forced extinction wave and spreading speed for
delayed population in time-periodic deteriorated environment,
\emph{Int. J. Biomath.} (2025), Paper No. 2550146.

\bibitem{WangWang2026}
M.-X. Wang and Y. Wang,
Properties of the generalized principal eigenvalue for a nonlocal
dispersal operator in one dimension,
\emph{Appl. Math. Lett.} 185 (2027), Paper No. 110124.

\bibitem{WangZhangTime2025}
M.-X. Wang and L. Zhang,
Generalized principal eigenvalue of time-periodic cooperative nonlocal
dispersal systems and applications,
preprint, arXiv:2502.10936 (2025).

\bibitem{WZ-1} M.-X. Wang and L. Zhang, Approximation of the generalized principal eigenvalue of cooperative nonlocal dispersal systems and applications, \emph{J. Math. Pures Appl.} 206 (2026), Paper No. 103835.

\bibitem{Weinberger2002}
H. F. Weinberger,
On spreading speeds and traveling waves for growth and migration models
in a periodic habitat,
\emph{J. Math. Biol.} 45 (2002), 511--548.

\bibitem{WFLS} H. Wu, Y.-X. Feng, W.-T. Li, and J.-W. Sun, Principal eigenvalue theory of nonlocal dispersal cooperative systems in unbounded domains and applications, \emph{J. Differential Equations} 445 (2025), Paper No. 113592.

\bibitem{WuLiSunVo2026}
H. Wu, W.-T. Li, J.-W. Sun, and H.-H. Vo,
Principal spectral theory and asymptotic analysis for time-periodic
nonlocal dispersal cooperative systems and applications,
\emph{J. Differential Equations} 481 (2026), Paper No. 114661.

\bibitem{XinLiBao2023}
M.-Z. Xin, W.-T. Li, and X. Bao,
Propagation dynamics of nonlocal dispersal monostable equations in
time-space periodic habitats,
\emph{J. Differential Equations} 365 (2023), 690--749.

\bibitem{YangLiWang2020}
F.-Y. Yang, W.-T. Li, and J.-B. Wang,
Wave propagation for a class of non-local dispersal non-cooperative
systems,
\emph{Proc. Roy. Soc. Edinburgh Sect. A} 150 (2020),
1965--1997.

\end{thebibliography}
\end{document}